\documentclass[11pt]{article}
\usepackage[utf8]{inputenc}

\title{\textbf{\Large Hyperbolic components and iterated monodromy of\\ polynomial skew-products of $\cmplex^2$}}
\author{Virgile Tapiero}
\date{January 13, 2025}

\usepackage[T1]{fontenc}
\usepackage[english, french]{babel}
\usepackage{hyperref}
\usepackage{amsmath, amssymb, amscd, amsthm}
\usepackage[all]{xy} 
\usepackage{esint}
\usepackage{color}
\usepackage{multirow}
\usepackage{tikz}

\usepackage{cite}

\usepackage{hhline}
\usepackage{nicematrix}

\usepackage{booktabs}  
\usepackage{array}     
\usepackage{amsmath} 

\newtheorem{thm}{Theorem}[section]

\newtheorem{lemme}[thm]{Lemma}
\newtheorem{prop}[thm]{Proposition}
\newtheorem{cor}[thm]{Corollary}
\newtheorem{defn}[thm]{Definition}

\newtheorem{exemple}[thm]{Example}

\newtheorem*{theorem*}{Theorem}
\newtheorem*{theoremA*}{Theorem A}
\newtheorem*{theoremA.I*}{Theorem A.I}
\newtheorem*{theoremA.II*}{Theorem A.II}
\newtheorem*{theoremfrA.I*}{Théorème A.I}
\newtheorem*{theoremfrA.II*}{Théorème A.II}
\newtheorem*{theoremB*}{Theorem B}
\newtheorem*{theoremB.I*}{Theorem B.I}
\newtheorem*{theoremB.II*}{Theorem B.II}
\newtheorem*{theoremfrB.I*}{Théorème B.I}
\newtheorem*{theoremfrB.II*}{Théorème B.II}
\newtheorem*{theoremfrB.I.bis*}{Théorème B.I.bis}
\newtheorem*{propositionC*}{Proposition C}
\newtheorem*{propositionC.I*}{Proposition C.I}
\newtheorem*{corollaryD*}{Corollary D}
\newtheorem*{corollaryD.I*}{Corollary D.I}
\newtheorem*{corollaryfrD.I*}{Corollaire D.I}
\newtheorem*{corollaryE*}{Corollary E}
\newtheorem*{corollaryE.I*}{Corollary E.I}
\newtheorem*{corollaryfrE.I*}{Corollaire E.I}
\newtheorem*{corollaryC*}{Corollary C}

\newcommand{\reels}{\mathbb{R}}
\newcommand{\cmplex}{\mathbb{C}}

\newcommand{\Sk}[2]{\mathrm{Sk}(#1,#2)}

\providecommand{\keywords}[1]
{
  {\small \textbf{\textit{Keywords---}} #1}
}

\providecommand{\acknowledgements}[1]
{
  {\textbf{\textit{Acknowledgements---}} #1}
}

\newenvironment{changemargin}[2]{\begin{list}{}{%
\setlength{\topsep}{0pt}%
\setlength{\leftmargin}{0pt}%
\setlength{\rightmargin}{0pt}%
\setlength{\listparindent}{\parindent}%
\setlength{\itemindent}{\parindent}%
\setlength{\parsep}{0pt plus 1pt}%
\addtolength{\leftmargin}{#1}%
\addtolength{\rightmargin}{#2}%
}\item }{\end{list}}

\usepackage[left=3.20cm,right=3.20cm,top=2cm,bottom=2cm]{geometry}

\usepackage{pdfpages}
\usepackage{graphicx}
\graphicspath{ {./images/} }

\begin{document}
\hypersetup{pdfborder=0 0 0}

\maketitle

\begin{center}
    \begin{changemargin}{0.000005cm}{0.000005cm}
        \begin{otherlanguage}{english}
            \begin{abstract}
            {\footnotesize We study the hyperbolic components of the family $\Sk{p}{d}$ of regular polynomial skew-products of $\cmplex^2$ of degree $d\geq2$, with a fixed base $p\in\cmplex[z]$. Using a homogeneous parametrization of the family, we compute the accumulation set $E$ of the bifurcation locus on the boundary of the parameter space. Then in the case $p(z)=z^d$, we construct a map $\pi_0(\mathcal{D}')\to AB_d$ from the set of unbounded hyperbolic components that do not fully accumulate on $E$, to the set of algebraic braids of degree $d$. This map induces a second surjective map $\pi_0(\mathcal{D}')\to\mathrm{Conj}(\mathfrak{S}_d)$ towards the set of conjugacy classes of permutations on $d$ letters. This article is a continuation in higher degrees of the work of Astorg-Bianchi \cite{AstBian23} in the quadratic case $d=2$, for which they provided a complete classification of the hyperbolic components belonging to $\pi_0(\mathcal{D}')$.
            } 
            \end{abstract}
        \end{otherlanguage}
        \phantom{0}\\
    \end{changemargin}
\end{center}

\section{Introduction}

This article is devoted to study hyperbolic components of some holomorphic families of polynomial mappings of $\cmplex^2$. 

In dimension $1$, the theory of stability within holomorphic families of rational maps is well established since the seminal work of Lyubich \cite{Lyu83} and Ma\~{n}\'{e}-Sad-Sullivan \cite{MSS83}. In such a family, it is known that stability preserves hyperbolicity in the sense that a stable component is either composed only of hyperbolic parameters or none at all. In the first case, the component is called a \textit{hyperbolic component}. Examples of hyperbolic components are given for each degree $d\geq2$ by the shift locus $\mathcal{S}_d$ of $\mathcal{P}_d$, where $\mathcal{P}_d$ is the space of monic polynomial mappings of degree $d$ modulo affine conjugacy. The locus $\mathcal{S}_d$ correspond to parameters whose critical points are all escaping to infinity by iteration. The terminology \textit{shift} comes from the following well known fact: if $p\in\mathcal{S}_d$, then there exists a homeomorphism between the Julia set of $p$ and the Cantor set $\{1,\cdots,d\}^{\mathbb{N}}$ which semi-conjugates $p$ to the one-sided shift on $\{1,\cdots,d\}^{\mathbb{N}}$. We refer to \cite{DeM12,DeMPil17,DeMPil11} for further information and recent investigations on shift loci and more.

The topic of hyperbolic components within holomorphic families of polynomial mappings in higher dimensions has recently been developed by Astorg-Bianchi. They proved in \cite{AstBian23} that stability preserves hyperbolicity for holomorphic families of \textit{regular polynomial skew-products} of $\mathbb{C}^2$, thus giving a sense to the notion of hyperbolic component in this setting. 
In the quadratic case, they classify hyperbolic components analogous to the shift locus $\mathcal{S}_2$. In this article, we study for each degree $d\geq2$ the families of skew-products $\Sk{p}{d}$ and their hyperbolic components analogous to shift loci which we introduce below.\\

A polynomial skew-product of degree $d\geq2$ on $\cmplex^2$ is a polynomial mapping of the form $f(z,w)=(p(z),q(z,w))$, where $p(z)$ and $q(z,w)$ are two polynomial mappings of degree $d$. The map $f$ is said to be \textit{regular} if it extends holomorphically to $\mathbb{P}^2(\cmplex)$. Let us fix $p\in\cmplex[z]$ a polynomial map of degree $d$. The set of all regular skew-products of degree $d$ with a fixed base $p$, modulo conjugacy by affine maps of $\cmplex^2$, is denoted $\Sk{p}{d}$. There exists $\Sk{p}{d}\simeq (f)_{\lambda\in\cmplex^{D_d}}$ a natural parametrization of this family by an affine space $\cmplex^{D_d}$ of dimension $D_d=\frac{1}{2}(d-1)(d+4)$, such that the discriminant $P(\lambda,z):=\mathrm{Disc}_w(q_{\lambda,z}(w))$ is homogeneous of degree $d(d-1)$ in $\lambda\in\cmplex^{D_d}$, see Section \ref{sec:homogeneousparametrization}. We identify the two families $\Sk{p}{d}$ and $(f_{\lambda})_{\lambda\in\cmplex^{D_d}}$ in what follows. 

{\subsection{Main results}\label{sec:mainresults}}
 
We can now state our first result, valid for any $p\in\cmplex[z]$ of degree $d\geq2$. Denote $\mathrm{Bif}$ the support of the \textit{bifurcation current} of $\Sk{p}{d}$, see Section \ref{sec:bifurcation}. Denote $J_p$ the Julia set of $p$.\\

\begin{thm}\label{thm:accumulationBIFalinfini}
    The bifurcation locus $\mathrm{Bif}\subset\cmplex^{D_d}$ accumulates at infinity exactly on $E$ with
    $$E:=\bigcup_{z\in J_p}\left\{[\lambda]\in\mathbb{P}^{D_d-1}_{\infty}:P(\lambda,z)=0\right\}.$$
\end{thm}

Let us now turn our attention to hyperbolic components of $\Sk{p}{d}$. The dynamics of polynomial skew-products have been studied in detail by Jonsson \cite{jon99}. In particular, he studied regular skew-products $f$ of $\cmplex^2$ that are hyperbolic in the sense that $f$ is uniformly expanding on its (small) Julia set, see \eqref{eq:fhyperboliconJ(f)}. Jonsson gave other equivalent conditions for a skew-product to be hyperbolic that are similar to characterizations known for polynomial mappings of one variable, see Theorem \ref{thm:Jonssonequivalenceofhyperbolicity}. In consequence, the connected components of the following set $\mathcal{D}$ are hyperbolic components (when $p$ is itself  hyperbolic):
$$\mathcal{D}:=\left\{\lambda\in\cmplex^{D_d}\ :\ \forall c\in \mathrm{Crit}(q_{\lambda,z}),\ (Q_{\lambda,z}^n(c))_n\ \mathrm{is\ unbounded}\right\},$$
where $Q_{\lambda,z}^{n}$ is the non-autonomous iterates of $q_{\lambda,z}$ defined by:
$$Q_{\lambda,z}^{n}(w):=q_{\lambda,z_{n-1}}\circ\cdots\circ q_{\lambda,z_0}(w),\ z_k:=p^{k}(z).$$
The connected components $\Lambda\in\pi_0(\mathcal{D})$ can be seen as two dimensional counterparts of the shift locus $\mathcal{S}_d$ in dimension $1$ which is connected, see \cite{DeMPil11}. In the sequel, we study the connected components of a subset $\mathcal{D}'$ of $\mathcal{D}$ defined by:
\begin{equation}\label{eq:lensembleD'prime}
    \mathcal{D}':=\{\lambda\in\mathcal{D}\ :\ \exists\ \mathrm{a\ continuous\ path\ inside}\ \mathcal{D}\ \mathrm{joining}\ \lambda\ \mathrm{to}\ \mathbb{P}^{D_d}\backslash E\}.
\end{equation} 
The connected components of $\mathcal{D}'$ are also connected components of $\mathcal{D}$, thus are hyperbolic components when $p$ is hyperbolic. In the case $d=2$, if the Julia set $J_p$ of $p$ is not totally disconnected, Astorg-Bianchi \cite{AstBian23} proved that the set $\pi_0(\mathcal{D}')$ is in one-to-one correspondence with the set 
\begin{equation}\label{eq:shift_p}
    \mathcal{S}_p:=\left\{s:\pi_0(K_p)\to\mathbb{N}:\sum_{U\in\pi_0(K_p)}s(U)\leq 2\right\},
\end{equation}
where $\pi_0(K_p)$ denotes the set of bounded Fatou components of $p$, establishing a classification of the connected components $\Lambda\in\pi_0(\mathcal{D}')$. We give more details about Astorg-Bianchi classification for $p(z)=z^2$ in Section \ref{sec:comparaison}. Our next results, Theorem \ref{thm:algebraicbraidsinJuliasets} and Corollary \ref{cor:abcmap} below, are a continuation in higher degrees $d\geq2$ of this classification in the case $p(z)=z^d$.

The statements of Theorem \ref{thm:algebraicbraidsinJuliasets} and Corollary \ref{cor:abcmap} involve \textit{algebraic braids} of degree $d$. Let us explain briefly what is an algebraic braid in our context, we refer to Section \ref{sec:AlgBraidsandJuliasets} for details. Let $\pi:\mathbb{S}^1\times\cmplex\to\mathbb{S}^1$ be the first projection. A disjoint union $K$ of closed curves in $\mathbb{S}^1\times\cmplex$ such that $\pi:K\to\mathbb{S}^1$ is a smooth unramified covering map of degree $d$, is called a \textit{closed braid}. If moreover $K$ has the form $K=\{q(z,w)=0\}\cap(\mathbb{S}^1\times\cmplex)$, where $q\in\cmplex[z,w]$, we say that $K$ is an \textit{algebraic braid} of degree $d$. 

It is known that the set of ambient isotopy classes in $\mathbb{S}^1\times\cmplex$ (for isotopies equal to $\mathrm{Id}$ on the base) of closed braids of degree $d$, are in one-to-one correspondence with $\mathrm{Conj}(B_d)$, the set of conjugacy classes of the braid group $B_d$ of Artin (see Section \ref{sec:AlgBraidsandJuliasets} for details). We denote $AB_d$ the set of conjugacy classes of $B_d$ which correspond to classes of algebraic braids of degree $d$. By abuse of language, the elements of $AB_d$ are also called algebraic braids. If $\mathcal{K}\in AB_d$ corresponds to the class of an algebraic braid $\{q=0\}\cap(\mathbb{S}^1\times\cmplex)$, we denote $\mathrm{ab}(q):=\mathcal{K}$, and if $q=q_{\lambda}$ for some $\lambda\in\cmplex^{D_d}$, we denote $\mathrm{ab}(f_{\lambda}):=\mathcal{K}$.

Given an algebraic braid of the form $\{q=0\}\cap (\mathbb{S}^1\times\cmplex)$, the monodromy of the braid above $\mathbb{S}^1$ induces a permutation $S_q$ on the roots $\{w\in\cmplex:q(1,w)=0\}$. By labeling these roots with numbers between $1$ and $d$, we can identify the permutation $S_q$ with an element of $\mathfrak{S}_d$, the group of permutations of $\{1,\cdots,d\}$. Moreover, the way the roots are labeled does not affect the conjugacy class of  $S_q$ in $\mathfrak{S}_d$, thus we have a well defined map:
\begin{equation}\label{eq:themapc:AB_dtoCB_d}
    \mathrm{c}:\left|
            \begin{array}{ll}
              AB_d & \longrightarrow\ \mathrm{Conj}(\mathfrak{S}_d)\\
              \textcolor{white}{} & \textcolor{white}{}  \textcolor{white}{}\\
              \mathrm{ab}(q)   & \longmapsto\ S_q\ \mathrm{mod.\ conj.}\\
            \end{array}
        \right.
\end{equation}

In the following statement $J(f_{\lambda})$ is the (small) Julia set of $f_{\lambda}$, see Section \ref{sec:partitionCBD_hyperbolicity}. To avoid any confusion, let us also mention that in this statement an isotopy class of a compact subset $C\subset\mathbb{S}^1\times\cmplex$ is the class obtained by identifying any compact subset $C'\subset\mathbb{S}^1\times\cmplex$ to $C$ if there exists an ambiant isotopy of $\mathbb{S}^1\times\cmplex$, equal to $\mathrm{Id}$ on the base, deforming $C$ into $C'$, see \eqref{eq:equivrelaambientisotopy}.
\begin{thm}\label{thm:algebraicbraidsinJuliasets}
    Let $\Lambda\in\pi_0(\mathcal{D}')$ with $\mathcal{D}'\subset \Sk{z^d}{d}$ defined by \eqref{eq:lensembleD'prime}. 
    \begin{enumerate}
        \item[\textbullet] For each $\lambda\in\Lambda$, the connected components of $J(f_{\lambda})$ are continuous closed loops winding at least one time above $\mathbb{S}^1$ in $\mathbb{S}^1\times\cmplex$.
        \item[\textbullet] Let $\lambda_0\in\Lambda$ such that $[\lambda_0]\notin E$. If the norm of $\lambda$ if sufficiently large, then the set $\{q_{\lambda_0}=0\}\cap(\mathbb{S}^1\times\cmplex)$ is an algebraic braid of degree $d$. Assume it is now the case.
        \item[\textbullet] For each $\lambda\in\Lambda$, $f_{\lambda}$ admits exactly $d$ different fixed points in $J(f_{\lambda})\cap\pi^{-1}(1)$, and the ambient isotopy class $\mathrm{ab}(f_{\lambda})$ of the union of the connected components in $J(f_{\lambda})$ of these fixed points coincide with the isotopy  class $\mathrm{ab}(f_{\lambda_0}):$ $\mathrm{ab}(f_{\lambda})=\mathrm{ab}(f_{\lambda_0})$.
    \end{enumerate}
\end{thm}
This theorem shows that for any parameter $\lambda$ in a component $\Lambda\in\pi_0(\mathcal{D}')$, the Julia set $J(f_{\lambda})$ contains an algebraic braid of degree $d$ which does not vary as $\lambda$ varies in $\Lambda$.

\begin{cor}\label{cor:abcmap}
    Let $d\geq2$ and let $\Lambda\in\pi_0(\mathcal{D}')$ with $\mathcal{D}'\subset \Sk{z^d}{d}$ defined by \eqref{eq:lensembleD'prime}. Then the following map is well defined:
    \begin{equation}\label{eq:themapab:pi_0(D')_dtoAB_d}
        \mathrm{ab}:\left|
            \begin{array}{ll}
              \pi_0(\mathcal{D}') & \longrightarrow\ AB_d\\
              \textcolor{white}{} & \textcolor{white}{}  \textcolor{white}{}\\
              \Lambda\ni\lambda   & \longmapsto\ \mathrm{ab}(f_{\lambda})\\
            \end{array}
        \right.
    \end{equation}
    Moreover, the composition $\mathrm{cab}:=\mathrm{c}\circ\mathrm{ab}$ with the map \eqref{eq:themapc:AB_dtoCB_d} gives a surjective map:
    \begin{equation}\label{eq:themapab:AB_dtoCB_d}
        \mathrm{cab}:\left|
            \begin{array}{ll}
              \pi_0(\mathcal{D}') & \overset{\mathrm{ab}\ }{\longrightarrow}\ AB_d \ \ \overset{\mathrm{c}\ }{\longrightarrow}\ \mathrm{Conj}(\mathfrak{S}_d)\\
              \textcolor{white}{} & \textcolor{white}{}  \textcolor{white}{}\\
              \Lambda\ni\lambda   & \longmapsto\ \mathrm{ab}(f_{\lambda}) \longmapsto\ S_{q_{\lambda}}\ \mathrm{mod.\ conj.}    \\
            \end{array}
        \right.
    \end{equation}
\end{cor}

The proof of Theorem \ref{thm:algebraicbraidsinJuliasets} relies on the computation of the \textit{iterated monodromy} of $f_{\lambda_0}$ on a pre-image tree, with $\lambda_0\in\mathcal{D}'$ being an \textit{admissible parameter}, see Definition \ref{defn:admissileparameters}. If $\lambda_0\in\Lambda$ with $[\lambda_0]\notin E$ is of sufficiently large norm, $\lambda_0$ is admissible and $\{q_{\lambda_0}=0\}\cap(\mathbb{S}^1\times\cmplex)$ is an algebraic braid of degree $d$, obtaining the second item of Theorem \ref{thm:algebraicbraidsinJuliasets}, see Proposition \ref{prop:lambdaisadmissible}. Roughly speaking, $\lambda_0$ admissible means that $f=f_{\lambda_0}$ is a partial self-cover of $\mathbb{S}^1\times\mathbb{D}_R$, for some $R>0$, and $(f^{-n}(\mathbb{S}^1\times\mathbb{D}_R))_{n\geq1}$ is a nested sequence such that $\bigcap_{n\geq0}f^{-n}(\mathbb{S}^1\times\mathbb{D}_R)=J(f)$. In this case, we can associate to $f$ (and a base point $t_0=(1,0)$) a group $\mathrm{IMG}(f,t_0)$ called the iterated monodromy group of $f$. $\mathrm{IMG}(f,t_0)$ acts on the pre-image tree $\mathcal{X}$ defined by:
$$\mathcal{X}:=\bigcup_{n=0}^{+\infty}f^{-n}(t_0)\supset X:=\bigcup_{n=0}^{+\infty}f^{-n}(t_0)\cap\pi^{-1}(1).$$ 
The group $\mathrm{IMG}(f,t_0)$ is defined as the quotient of the fundamental group $\pi_1(\mathbb{S}^1\times\mathbb{D}_R,t_0)$ by the kernel of the monodromy action of $f$ on $\mathcal{X}$. A key point is that the action of $\mathrm{IMG}(f,t_0)$ on the tree $X\subset\mathcal{X}$ defined above can be computed explicitly. Let us specify here that $\mathrm{IMG}(f,t_0)$ does not preserve $X$, the elements of $\mathrm{IMG}(f,t_0)$ which preserve the set $X$ are the elements of the following collection of subgroups of $\pi_1(\mathbb{S}^1\times\mathbb{D}_R,t_0)$:
\begin{equation}\label{eq:collectionsousgroupes}
    \mathbb{A}=\{\mathbb{A}_{d^n},n\in\mathbb{N}\}\hookrightarrow\mathrm{IMG}(f,t_0),
\end{equation}
where $\mathbb{A}_{d^n}$ is the subgroup generating by the loop $t\in[0,1]\mapsto(e^{2i\pi td^n},0)$. For each $n\geq0$, the set $\mathbb{A}_{d^n}$ acts on $\bigcup_{k=0}^{n}f^{-k}(t_0)\cap\pi^{-1}(1)$, we compute explicitly this action in Section \ref{sec:computationmonodromyactionofAonX}. We use the action of $\mathbb{A}$ on $X$ to analyze the connected components of the Julia set $J(f)$, allowing to obtain the first item of Theorem \ref{thm:algebraicbraidsinJuliasets}, see Section \ref{sec:topologicaldescriptionofJ(f)}.

For proving the third item of Theorem \ref{thm:algebraicbraidsinJuliasets}, we proceed in two steps. First, we show that the isotopy class $\mathcal{C}(f_{\lambda_0})$ of the union of the connected components of the fixed points in $J(f_{\lambda_0})\cap\pi^{-1}(1)$ is equal to the class $\mathrm{ab}(f_{\lambda_0})$, see Theorem \ref{thm:Kisotopicf^-1(S^1times0)}. Second, we use Theorem \ref{thm:isotopyoverJ_p} to construct for each $\lambda\in\Lambda$ an isotopy between $\mathcal{C}(f_{\lambda_0})$ and $\mathcal{C}(f_{\lambda})$ (the union of the connected components of the fixed points in $J(f_{\lambda})\cap\pi^{-1}(1)$) proving $\mathrm{ab}(f_{\lambda_0})=\mathcal{C}(f_{\lambda_0})=\mathcal{C}(f_{\lambda})=:\mathrm{ab}(f_{\lambda})$, see Lemma \ref{lemma:braidab(f_lambda)constantonLambda}.

An other application of the iterated monodromy action of $\mathbb{A}$ on $X$, is that we can show that $J(f_{\lambda_0})$ is homeomorphic to a topological space (a suspension) that involves only the permutation $S_{q_{\lambda_0}}$ of the algebraic braid $\mathrm{ab}(q_{\lambda_0})$ of $f_{\lambda_0}$. It is a topological result of independent interest.

\begin{thm}
    Let $\Lambda\in\pi_0(\mathcal{D}')$ and let $\lambda_0\in\Lambda$ be an admissible element in the sense of $\mathrm{Defintion\ \ref{defn:admissileparameters}}$. Let $S:=S_{q_{\lambda_0}}\in\mathfrak{S}_d$ be the permutation given by the map $\mathrm{\eqref{eq:themapc:AB_dtoCB_d}}$. We denote $h:(a_n)_{n\geq1}\mapsto\left(S^{\circ d^{n-1}}(a_n)\right)_{n\geq1}$ which is a homeomorphism of $\{1,\cdots,d\}^{\mathbb{N}^*}$. Then there exists a homeomorphism:
    $$J(f_{\lambda_0})\overset{\simeq}{\longrightarrow}\mathcal{S}_{\Lambda}:=\left([0,1]\times\{1,\cdots,d\}^{\mathbb{N}^*}\right)\left/\displaystyle{\phantom{\int}\!\!\!\!\!\!(0,x)\sim_h(1,y)\Leftrightarrow h(x)=y}\right. .$$
    Moreover, for any  other parameter $\lambda\in\Lambda$, $J(f_{\lambda})$ is also homeomorphic to $\mathcal{S}_{\Lambda}$. 
\end{thm}

Iterated monodromy groups have been used in complex dynamics by different authors. In the article \cite{BartNekra06}, Bartholdi-Nekrashevych solved a problem posed by Hubbard in dimension $1$ (the \textit{twisted rabbit problem}) by using some iterated monodromy groups. Nekrashevych introduced later a general definition of iterated monodromy groups and gave applications to complex dynamics, we refer to \cite{Nekra11}. 

In practice, computing iterated monodromy groups can be difficult, and there are very few examples of them in complex dynamics in higher dimensions. The ones in \cite{BelchKoch09, Nekra12, Bowman21} are the only ones we know. To the best of our knowledge, this is the first time iterated monodromy has been used to study holomorphic families in higher dimensions.

\subsection{The quadratic case}

For the quadratic family $f_{\lambda}=(z^2,w^2+a^2z^2+b^2z+c^2),$ with $\lambda=(a,b,c)\in\cmplex^{D_2}=\cmplex^3$, the situation is much easier to understand than for higher degrees $d>2$. Indeed, a braid made of two strings is equivalent to the datum of an integer $n\in\mathbb{N}$ indicating how many times the two strings are braided together, with a sign $\pm$ to indicate whether the strings are braided in a direct or indirect sense. It is convenient to represent these braids by using a representative diagram as follows:
$$
\tikzset{every picture/.style={line width=0.75pt}}    
\begin{tikzpicture}[x=0.75pt,y=0.75pt,yscale=-1,xscale=1] 
\draw    (7,162) -- (7,213) ;
\draw    (147,162) -- (147,212) ;
\draw    (7,203) .. controls (36,202) and (9,173) .. (37,173) ; 
\draw    (46,203) .. controls (75,202) and (48,173) .. (76,173) ;
\draw    (86,203) .. controls (115,202) and (88,173) .. (116,173) ;
\draw    (117,202) .. controls (146,201) and (119,172) .. (147,172) ;
\draw [color={rgb, 255:red, 0; green, 0; blue, 0 }  ,draw opacity=1 ]   (7,172) .. controls (16,172) and (14,179) .. (19,185) ;
\draw    (26,193) .. controls (34,203) and (35,203) .. (46,203) ;
\draw [color={rgb, 255:red, 0; green, 0; blue, 0 }  ,draw opacity=1 ]   (37,173) .. controls (57,173) and (52,177) .. (56,185) ;
\draw    (66,193) .. controls (74,203) and (75,203) .. (86,203) ; 
\draw [color={rgb, 255:red, 0; green, 0; blue, 0 }  ,draw opacity=1 ]   (76,173) .. controls (96,173) and (91,177) .. (95,185) ;
\draw [color={rgb, 255:red, 0; green, 0; blue, 0 }  ,draw opacity=1 ]   (104,194) .. controls (109,201) and (110,202) .. (117,202) ;
\draw [color={rgb, 255:red, 0; green, 0; blue, 0 }  ,draw opacity=1 ]   (116,173) .. controls (121,174) and (125,175) .. (130,180) ;
\draw [color={rgb, 255:red, 0; green, 0; blue, 0 }  ,draw opacity=1 ]   (136,185) .. controls (140,193) and (135,204) .. (147,202) ;
\draw (-1,223) node [anchor=north west][inner sep=0.75pt]   [align=left] {\ \ \ \ \ \ \ \ \ \ \ $\sigma_1^{+4}\in B_2$};
\draw (-44,196) node [anchor=north west][inner sep=0.75pt]   [align=left] { \ \ \ \ \ \ \ \ $1$ \ \ \ \ \ \ \ };
\draw (110,196) node [anchor=north west][inner sep=0.75pt]   [align=left] { \ \ \ \ \ \ \ \ $1$ \ \ \ \ \ \ \ };
\draw (-44,165) node [anchor=north west][inner sep=0.75pt]   [align=left] { \ \ \ \ \ \ \ \ $2$ \ \ \ \ \ \ \ };
\draw (110,165) node [anchor=north west][inner sep=0.75pt]   [align=left] { \ \ \ \ \ \ \ \ $2$ \ \ \ \ \ \ \ };
\end{tikzpicture}
$$
Here, $\sigma_1$ is the generator of the Artin group $B_2$, see \eqref{eq:GroupeArtin}. This diagram represents a braid made of two strings which are braided together $n=4$ times in a direct way. It is encoded by the word $\sigma_1^{+4}$ in the group $B_2$.

Astorg-Bianchi \cite{AstBian23} proved that a braid $K_{\lambda}=f_{\lambda}^{-1}(\mathbb{S}^1\times\{0\})$ (for adapted parameters $\lambda\in\mathcal{D}'$) is composed of two disjoint closed curves or only one closed curve. In the first case, they also computed the winding number $\mathrm{Wind}(K_{\lambda})$ of the two curves which is equal to is $0$ or $1$ (see Lemma \ref{lemma:lemmeAB}). The word in $B_2$ associate to the braid $K_{\lambda}$ is equal to $\sigma_1^{2\mathrm{Wind}(K_{\lambda})}$ if $K_{\lambda}$ is disconnected, and is equal to  $\sigma_1^{0}$ if $K_{\lambda}$ is connected, see Corollary \ref{eq:abford=2}. If $K_{\lambda}$ is disconnected, we have 
$$\mathrm{Wind}(K_{\lambda})=\frac{\#\{z\in\mathbb{D}:a^2z^2+b^2z+c^2=0\}}{2},\ \mathrm{where}\ \lambda=(a,b,c).$$
 
Thus, braids do not necessarily need to be involved in understanding the combinatorial structure of the connected components of $\pi_0(\mathcal{D}')$, it is sufficient to compute $s:\lambda=(a,b,c)\mapsto\#\{z\in\mathbb{D}:a^2z^2+b^2z+c^2=0\}$. We observe that the datum of $\lambda\mapsto s(\lambda)$ is equivalent to the explicit computation of the set $\mathcal{S}_p$ for $p(z)=z^2$, where $\mathcal{S}_p$ is defined by \eqref{eq:shift_p}. We recall these results of Astorg-Bianchi in detail in Section \ref{sec:comparaison}, and provide a comparison between their classification and the results obtained in this article in the same section.

\subsection{Examples of degree \texorpdfstring{$d=3$}{TEXT}}

In what follows $f(z,w)=(z^3,q_z(w))$ is a regular skew-product of degree $d=3$ for which there exists $\lambda\in\cmplex^{D_3}=\cmplex^7$ such that $f=f_{\lambda}$, see Section \ref{sec:homogeneousparametrization} for the definition of $f_{\lambda}$. We give three different examples below where $f$ is hyperbolic and such that $K:=\{q=0\}\cap(\mathbb{S}^1\times\cmplex)$ is an algebraic braid of degree $3$. We also give an explicit computation of the braid $K$.

As mentioned in Section \ref{sec:mainresults}, there exists a bijective map from the set of ambient isotopy classes of closed braids of degree $3$ towards the set $\mathrm{Conj}(B_3)$ of conjugacy classes of $B_3$. Thus, we can associate to $K$ a word $\sigma(K)$ in $B_3$ representing the class of $K$ in $\mathrm{Conj}(B_3)$, we do it for each example, and we also provide a diagram associated to $\sigma(K)$. It then allows us to deduce the permutation $S_q$ (up to conjugacy) associated to the monodromy of $K$ above $\mathbb{S}^1$ i.e. the permutation given by the map \eqref{eq:themapc:AB_dtoCB_d}.

For computing explicitly $K$, we proceed as follows. Observe that $K=f^{-1}(\mathbb{S}^1\times\{0\})$, thus $K$ can be determined by finding the different lifts of the curve $\gamma_d(t)=(e^{2i\pi dt},0)$ under $f$ (with $d=3$) which contain the roots $\{(1,w)\ : q_1(w)=0\}$ above $z=1$. Indeed, this works since $f$ is a partial self-cover on $\mathbb{S}^1\times\mathbb{D}_R$ for some $R>0$, we choose $\lambda\in\cmplex^7$ of sufficiently large norm to ensure this property. Further details can be found in the proof of Proposition \ref{lemma:particularexamples}, where the same idea is used to compute braids and their associated permutations.

\begin{exemple} 
    $\mathrm{Take}\ \lambda=(0,0,0,-2,0,0,0):$
    $$q_z(w)=(w-2z)(w-2e^{i\pi/3}z)(w-2e^{2i\pi/3}z)=w^3-8z^3.$$
\end{exemple}
\noindent The roots above $z=1$ are $\{q_1(w)=0\}=\{2,2e^{i\pi/3},2e^{2i\pi/3}\}=:\{x_1,x_2,x_3\}$ thus
$$K=\alpha_1([0,1])\sqcup\alpha_2([0,1])\sqcup\alpha_3([0,1]),\ \mathrm{with}\ \alpha_k(t):=(e^{2i\pi t},x_ke^{2i\pi t}),$$
where $\alpha_1$, $\alpha_2$ and $\alpha_3$ are the three lift of $\gamma_3$ by $f$. Thus, $\sigma(K)$ and its diagram are given by:
$$
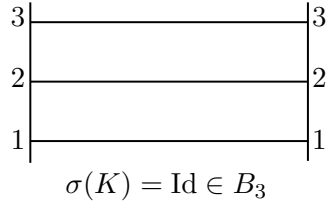
\tikzset{every picture/.style={line width=0.75pt}}        
\begin{tikzpicture}[x=0.75pt,y=0.75pt,yscale=-1,xscale=1]
\draw    (100,120) -- (100,201) ;
\draw    (240,120) -- (240,200) ;
\draw    (100,130) -- (240,130) ;
\draw    (100,160) -- (240,160) ;
\draw    (100,190) -- (240,190) ;
\draw (92,204) node [anchor=north west][inner sep=0.75pt]   [align=left] { \ \ \ \ \ $\sigma(K)=\mathrm{Id}\in B_3$ \ \ \ \ \ \ \ };
\draw (50,183) node [anchor=north west][inner sep=0.75pt]   [align=left] { \ \ \ \ \ \ \ \ $1$ \ \ \ \ \ \ \ };
\draw (202,183) node [anchor=north west][inner sep=0.75pt]   [align=left] { \ \ \ \ \ \ \ \ $1$ \ \ \ \ \ \ \ };
\draw (50,152) node [anchor=north west][inner sep=0.75pt]   [align=left] { \ \ \ \ \ \ \ \ $2$ \ \ \ \ \ \ \ };
\draw (202,152) node [anchor=north west][inner sep=0.75pt]   [align=left] { \ \ \ \ \ \ \ \ $2$ \ \ \ \ \ \ \ };
\draw (50,121) node [anchor=north west][inner sep=0.75pt]   [align=left] { \ \ \ \ \ \ \ \ 3 \ \ \ \ \ \ \ };
\draw (202,121) node [anchor=north west][inner sep=0.75pt]   [align=left] { \ \ \ \ \ \ \ \ 3 \ \ \ \ \ \ \ };
\end{tikzpicture}
$$
The permutation associated to this braid is $S_q=\mathrm{Id}$.

\begin{exemple}
    $\mathrm{Take}\ \lambda=(0,-2,0,0,0,0,0):$ $q_z(w)=w^3-8z.$
\end{exemple}
\noindent The roots above $z=1$ are $\{q_1(w)=0\}=\{2e^{2i\pi \frac{k-1}{3}},\ k=1,2,3\}=:\{x_1,x_2,x_3\}$. Then $K$ is this time given by a single curve:
$$K=\alpha_1([0,3])=\alpha_2([0,3])=\alpha_3([0,3]),\ \mathrm{with}\ \alpha_k(t):=(e^{2i\pi t},x_ke^{2i\pi \frac{t}{3}}).$$ 
Thus, $\sigma(K)\in B_2$ and its diagram are given by:
$$
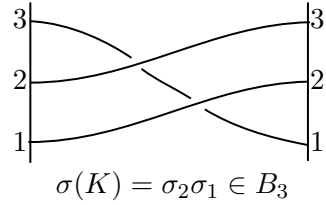
\tikzset{every picture/.style={line width=0.75pt}}      
\begin{tikzpicture}[x=0.75pt,y=0.75pt,yscale=-1,xscale=1]
\draw    (140,160) -- (140,241) ;
\draw    (280,160) -- (280,240) ;
\draw    (191,188) .. controls (180,180.5) and (161,170.5) .. (140,169.5) ;
\draw    (280,232) .. controls (259,228) and (243,223) .. (228,214) ;
\draw    (140,200.5) .. controls (192,200.5) and (231,171) .. (281,170) ;
\draw    (139,230.5) .. controls (191,230.5) and (230,201) .. (280,200) ;
\draw    (196,193.5) .. controls (204,199.5) and (213,203.5) .. (221,208.5) ;
\draw (132,244) node [anchor=north west][inner sep=0.75pt]   [align=left] { \ \ \ \ $\sigma(K)=\sigma_{2} \sigma_{1}\in B_3$ \ \ \ \ \ \ \ };
\draw (91,224) node [anchor=north west][inner sep=0.75pt]   [align=left] { \ \ \ \ \ \ \ \ $1$ \ \ \ \ \ \ \ };
\draw (241,224) node [anchor=north west][inner sep=0.75pt]   [align=left] { \ \ \ \ \ \ \ \ $1$ \ \ \ \ \ \ \ };
\draw (91,193) node [anchor=north west][inner sep=0.75pt]   [align=left] { \ \ \ \ \ \ \ \ $2$ \ \ \ \ \ \ \ };
\draw (241,193) node [anchor=north west][inner sep=0.75pt]   [align=left] { \ \ \ \ \ \ \ \ $2$ \ \ \ \ \ \ \ };
\draw (91,162) node [anchor=north west][inner sep=0.75pt]   [align=left] { \ \ \ \ \ \ \ \ $3$ \ \ \ \ \ \ \ };
\draw (241,162) node [anchor=north west][inner sep=0.75pt]   [align=left] { \ \ \ \ \ \ \ \ $3$ \ \ \ \ \ \ \ };
\end{tikzpicture}
$$
The permutation associated to this braid is $S_q=(1,2,3)$.

\begin{exemple}
    $\mathrm{Take}\ \lambda=(0,0,0,0,0,4i,0):$ $q_z(w)=w^3-16zw=w(w^2-16z)$.
\end{exemple}
\noindent 
The roots above $z=1$ are $\{q_1(w)=0\}=\{0,4,-4\}=:\{x_1,x_2,x_3\}$ and $K$ is given by:
$$K=\alpha([0,1])\sqcup\beta([0,2]),\ \mathrm{with}\ \alpha(t):=(e^{2i\pi t},0)\ \mathrm{and}\ \beta(t)=(e^{2i\pi t}, x_ke^{2i\pi \frac{t}{2}}),$$
where $k\in\{2,3\}$ is chosen arbitrarily. Thus, $\sigma(K)$ and its diagram are given by:
$$
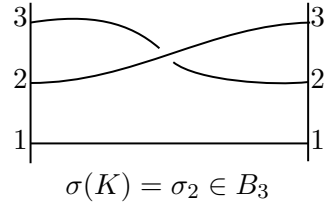
\tikzset{every picture/.style={line width=0.75pt}} 
\begin{tikzpicture}[x=0.75pt,y=0.75pt,yscale=-1,xscale=1]
\draw    (81,100) -- (81,181) ;
\draw    (221,100) -- (221,180) ;
\draw    (81,140.5) .. controls (133,140.5) and (172,111) .. (222,110) ;
\draw    (81,110) .. controls (112,105) and (131,110) .. (146,123) ;
\draw    (153,130) .. controls (165,139.5) and (204,142) .. (221,140) ;
\draw    (81,171) -- (221,171) ;
\draw (73,184) node [anchor=north west][inner sep=0.75pt]   [align=left] { \ \ \ \ \ $\sigma(K)=\sigma_2\in B_3$ \ \ \ \ \ \ \ };
\draw (32,163) node [anchor=north west][inner sep=0.75pt]   [align=left] { \ \ \ \ \ \ \ \ $1$ \ \ \ \ \ \ \ };
\draw (182,163) node [anchor=north west][inner sep=0.75pt]   [align=left] { \ \ \ \ \ \ \ \ $1$ \ \ \ \ \ \ \ };
\draw (32,132) node [anchor=north west][inner sep=0.75pt]   [align=left] { \ \ \ \ \ \ \ \ $2$ \ \ \ \ \ \ \ };
\draw (182,132) node [anchor=north west][inner sep=0.75pt]   [align=left] { \ \ \ \ \ \ \ \ $2$ \ \ \ \ \ \ \ };
\draw (32,101) node [anchor=north west][inner sep=0.75pt]   [align=left] { \ \ \ \ \ \ \ \ $3$ \ \ \ \ \ \ \ };
\draw (182,101) node [anchor=north west][inner sep=0.75pt]   [align=left] { \ \ \ \ \ \ \ \ $3$ \ \ \ \ \ \ \ };
\end{tikzpicture}
$$
The permutation associated to this braid is $S_q=(2,3)$.\\

\noindent\acknowledgements{The author thanks Matthieu Astorg for sharing his time through fruitful discussions, and for introducing him to the subject of skew-products and hyperbolic components. This research work was partially supported by the ANR Project ANR PADAWAN /ANR-21-CE40-0012-01, and by the French Italian University and Campus France through the Galileo program, under the project {From rational to transcendental: complex dynamics and parameter spaces}.}

\section{The parameter space \texorpdfstring{$\Sk{p}{d}$}{TEXT}}

{\subsection{Homogeneous parametrization}\label{sec:homogeneousparametrization}}

For each $\lambda=(a_{0,0},\cdots,a_{0,d},\cdots,a_{d-2,0},\cdots,a_{d-2,2})\in\cmplex^{D_d}$, where $D_d=\frac{1}{2}(d-1)(d+4)$, we denote $f_{\lambda}(z,w)=(p(z),q_{\lambda,z}(w))$, where $q_{\lambda,z}$ is defined by:
\begin{equation}\label{eq:parametricformoff_version_intro}
    q_{\lambda,z}(w):=w^d+\sum_{j=0}^{d-2}\left(\sum_{k=0}^{d-j}{(a_{j,k})}^{d-j}z^k\right)w^{j}.
\end{equation} 
Every regular skew-product of degree $d$, whose first component is equal to $p$, is conjugated to an element of $(f_{\lambda})_{\lambda\in\cmplex^{D_d}}$ via an affine map of $\cmplex^2$ of the form $A(z,w)=(z,B(z,w))$. Thus, by definition, for each element $f$ of the family $\Sk{p}{d}$ there exists $\lambda\in\cmplex^{D_d}$ such that $f=f_{\lambda}$. Throughout all the article, we identify the two families to the affine space $\cmplex^{D_d}$:
$$\Sk{p}{d}\simeq\{f_{\lambda},\ \lambda\in\cmplex^{D_d}\}\simeq\cmplex^{D_d}.$$

The mappings $q_{\lambda,z}(w)$ and $q_{\lambda,z}'(w)$ are homogeneous in $(\lambda,w)$ of degree $d$ and $d-1$ respectively. Thus their discriminant is homogeneous of degree $d(d-1)$:
\begin{equation}\label{eq:discriminantofq_lambdaz_version_intro}
         P(\lambda,z):=\mathrm{Res}_w(q_{\lambda,z}(w);q_{\lambda,z}'(w))=\prod_{q_{\lambda,z}'(c)=0}q_{\lambda,z}(c)\ \in\cmplex[\lambda,z].
\end{equation}
In particular, $\lambda\mapsto P(\lambda,z)$ is a map compatible with the projective compactification of $\cmplex^{D_d}$:
\begin{equation}\label{eq:compactificationoftheparameterspace_version_intro}
     \mathbb{P}^{D_d}:=\cmplex^{D_d}\sqcup \mathbb{P}_{\infty}^{D_d-1},\ \mathrm{with}\ \mathbb{P}_{\infty}^{D_d-1}=\left\{[\lambda],\ \lambda\in\cmplex^{D_d}\backslash\{0\}\right\}.
\end{equation}
We then denote:
\begin{equation}\label{defn:EetlesE_z_version_intro}      E:=\bigcup_{z\in J_p}E_z\ \mathrm{and}\ E_z:=\left\{[\lambda]\in\mathbb{P}^{D_d-1}_{\infty}:P(\lambda,z)=0\right\}.
\end{equation}

\subsection{The bifurcation locus \texorpdfstring{$\mathrm{Bif}$}{TEXT}}\label{sec:bifurcation}

Recall that the standard definition of the bifurcation current of a holomorphic family $(f_{\lambda})_{\lambda\in M}$ of endomorphisms of $\mathbb{P}^2$ of degree $d\geq2$
is the current $T_{\mathrm{bif}}:=dd^c\mathcal{L}$, where
the function $\mathcal{L}:\lambda\in M\mapsto \int_{\mathbb{P}^2}\mathrm{Log}|\mathrm{Jac}f_{\lambda}|\ \mathrm{d}\mu_\lambda$ is the Lyapunov function, and where $\mu_{\lambda}$ stands for the equilibrium measure of $f_{\lambda}$, we refer to \cite{BB07,BB09,BBD18}. In the context of the family $\Sk{p}{d}$, Astorg-Bianchi \cite{AstBian23} have shown Theorem \ref{thm:T_bif=intT_bif_zmu_p(z)} below which gives a more explicit form of $T_{\mathrm{bif}}$ using the vertical Green functions defined by \eqref{eq:fonctiondeGreenverticale}. For any $z\in J_p$ and any $\lambda\in\Sk{p}{d}$, the \textit{vertical Green function} $G_{\lambda,z}(w)$ of $q_{\lambda,z}$ is defined by:
\begin{equation}\label{eq:fonctiondeGreenverticale}
    G_{\lambda,z}(w)=\lim_{n\to+\infty}\frac{1}{d^n}\mathrm{Log}^+|Q^{n}_{\lambda,z}(w)|.
\end{equation}
The map $G_{\lambda,z}(w)$ is continuous in the three variables $(\lambda,z,w)\in\cmplex^{D_d}\times\cmplex^2$ and is plurisubharmonic in $\lambda$. In particular, $T_{\mathrm{bif}_z}:=dd^c_{\lambda}\left(\sum_{c\in C_{\lambda,z}}G_{\lambda,z}(c)\right)$ is a well defined positive current, where $C_{\lambda,z}:=\mathrm{Crit}(q_{\lambda,z})$ is the critical set of $q_{\lambda,z}$. We denote $\mathrm{Bif}_z:=\mathrm{Supp}(T_{\mathrm{bif}_z})$. Recall that $J_p$ is the support of the equilibrium measure $\mu_p$ of $p$.

\begin{thm}[Astorg-Bianchi {\cite[Thm. 3.3]{AstBian23}}]\label{thm:T_bif=intT_bif_zmu_p(z)}
    $T_{\mathrm{bif}}$ and its locus $\mathrm{Bif}:=\mathrm{Supp}(T_{\mathrm{bif}})$ satisfy:
    $$T_{\mathrm{bif}}=\int_{J_p}T_{\mathrm{bif}_z}\ \mathrm{d}\mu_p(z)\ ,\ \mathrm{Bif}=\overline{\bigcup_{z\in J_p}\mathrm{Bif}_z}\ \mathrm{and}\ \mathrm{Bif}_z=\partial\mathcal{B}_z,$$
    where $\mathcal{B}_z:=\{\lambda\in\cmplex^{D_d},\ \exists c\in C_{\lambda,z}:(Q^n_{\lambda,z}(c))_n\ \mathrm{is\ bounded}\}$. 
\end{thm}
The following lemma is also due to Astorg-Bianchi {\cite[Lemma 2.8]{AstBian23}}, we use it to show that the accumulation at infinity of $\mathrm{Bif}$ contains the set $E$ defined by \eqref{defn:EetlesE_z_version_intro}, see Section \ref{sec:proofofthm1.2}. 

\begin{lemme}\label{lemma:variationofB_zandBif_z}
    For each compact subset $M$, $z\mapsto\mathcal{B}_z\cap M$ is upper semi-continuous.
\end{lemme}

\subsection{Partition of the parameter space and hyperbolicity}\label{sec:partitionCBD_hyperbolicity}

Recall that $C_{\lambda,z}:=\mathrm{Crit}(q_{\lambda,z})$.
Following Astorg-Bianchi \cite{AstBian23} we introduce:

\begin{defn}\label{defn:CBD} 
    We denote for $z\in J_p$ and $\lambda\in\cmplex^{D_d}$ : 
    \begin{enumerate}
        \item[\textbullet] $K_z(f_{\lambda}):=\{w\in\cmplex\ :\ (Q^n_{\lambda,z}(c))_n\ \mathrm{is\ bounded}\}$: the filled Julia set above $z$.
        \item[\textbullet] $J_z(f_{\lambda}):=\partial K_z(f_{\lambda})$: the Julia set above $z$. 
    \end{enumerate}
    We then distinguish the following subsets of the parameter space:
    \begin{enumerate}
        \item[\textbullet] $\mathcal{C}:=\bigcap_{z\in J_p}\{\lambda\in\cmplex^{D_d} :\ \forall c\in C_{\lambda,z},\ c\in K_z(f_{\lambda})\}\subset\bigcap_{z\in J_p}\mathcal{B}_z$.
        \item[\textbullet] $\mathcal{D}:=\bigcap_{z\in J_p}\mathcal{D}_z$ with $\mathcal{D}_z:=\{\lambda\in\cmplex^{D_d}\ :\ \forall c\in C_{\lambda,z},\ G_{\lambda,z}(c)>0\}$.
        \item[\textbullet] $\mathcal{D}':=\{\lambda\in\mathcal{D} :\mathrm{there\ exists\ a\ continuous\ path\ inside}\ \mathcal{D}\ \mathrm{joining}\ \lambda\ \mathrm{to}\ \mathbb{P}^{D_d-1}_{\infty}\backslash E\}.$
    \end{enumerate}
    
\end{defn}
Note that a sequence $(Q^n_{\lambda,z}(w))_n$, with $z\in J_p$ and $w\in\cmplex$, is bounded if and only if $G_{\lambda,z}(w)=0$, where $G_{\lambda,z}(w)$ is defined by \eqref{eq:fonctiondeGreenverticale}. The set $\mathcal{C}$ is known (when $J_p$ is connected) as the \textit{connected locus} of the family, this terminology was introduced by Jonsson \cite{jon99}. It is compact as shown in Corollary \ref{cor:locusCiscompact}. The set $\mathcal{D}$ is open by continuity of the vertical Green function. We should also mention that for a regular skew-product $f(z,w)=(p(z),q(z,w))$ of degree $d\geq2$, we have \cite[Corollary 4.4]{jon99}:
\begin{equation}\label{eq:definitionJuliasetSK}
    J(f)=\overline{\bigcup_{z\in J_p}\{z\}\times J_z(f)},
\end{equation}
where $J(f)$ denotes the support of the equilibrium measure $\mu_f$ of $f$, it is called the (small) Julia set of $f$. 

In this article we are interested in \textit{hyperbolic skew-products}, Jonsson studied them in \cite{jon99}. We recall that a regular skew-product $f(z,w)=(p(z),q(z,w))$ of degree $d\geq2$ is hyperbolic if there exists two constants $c>0$ and $k>1$ such that for all $n\geq0$:
\begin{equation}\label{eq:fhyperboliconJ(f)}
    ||d_pf^n||\geq ck^n,\ \forall p\in J(f).
\end{equation}
We mention that the condition \eqref{eq:fhyperboliconJ(f)} requires that $p$ is itself hyperbolic. Jonsson showed in \cite{jon99} that when $p$ is hyperbolic, Condition \eqref{eq:fhyperboliconJ(f)} is actually equivalent to the following condition: there exists $K>1$ and $C>0$ such that for any $z\in J_p$ and for any $n\geq0$:
\begin{equation}\label{eq:fverticallyexpandingonJ(f)}
    |(Q_{z}^n)'(w)|\geq CK^n,\ \forall w\in J_z(f).
\end{equation}
If this second condition is satisfied we say that $f$ is \textit{vertically expanding} (over $J_p$), a terminology introduced by Jonsson \cite{jon99}. In this case, $p$ is not required to be hyperbolic, and for any $p$ of degree $d$ there exists $q$ such that $f(z,w)=(p(z),q(z,w))$ is a regular skew-product of degree $d$ which is vertically expanding. The notion of vertically expanding skew-products is thus more general. We also mention that Jonsson also showed that Condition \eqref{eq:fverticallyexpandingonJ(f)} implies :
\begin{equation}\label{eq:fhypimpliesJ=cupJ_z}
    J(f) = \bigcup_{z\in J_p}\{z\}\times J_z(f).
\end{equation}

It is well known that a polynomial map of one complex variable of degree $d\geq2$ is hyperbolic, if and only if, its postcritical set does not accumulate on its Julia set. Theorem \ref{thm:Jonssonequivalenceofhyperbolicity} below is the counterpart of this property for polynomial skew-products:

\begin{thm}[Jonsson {\cite{jon99}}]\label{thm:Jonssonequivalenceofhyperbolicity} 
Let $f\in\Sk{p}{d}$. Then the following points are equivalent:
    \begin{enumerate}
        \item $f$ is vertically expanding over $J_p$.
        \item $\overline{\bigcup_{n\geq1}f^n(C_\lambda)}\cap J(f) = \emptyset$, where $C_\lambda:=\bigcup_{z\in J_p}\{z\}\times C_{\lambda,z}$.
    \end{enumerate}
Moreover, $f$ is hyperbolic if and only if $p$ is hyperbolic and \textit{1.} or \textit{2.} is satisfied.
\end{thm}

It is known that if all the critical points of a polynomial mapping $q(w)$ of degree $d\geq2$ escape to infinity by iteration, then $q(w)$ is hyperbolic. This property is transferred to skew-products as follows. It can be shown by elementary arguments that the condition $\overline{\bigcup_{n\geq1}f^n(C_\lambda)}\cap J(f) = \emptyset$ is verified when $(Q^n_{\lambda,z}(c))_{n}$ escapes to infinity for any $(z,c)\in C_{\lambda}$, i.e. $\lambda\in\mathcal{D}$. So, Theorem \ref{thm:Jonssonequivalenceofhyperbolicity} implies the following:

\begin{cor}\label{cor:elementsofDarehyperbolic}
    If $\lambda\in\mathcal{D}$ then $f_{\lambda}$ is vertically expanding. If furthermore $p$ is hyperbolic then $f_{\lambda}$ is hyperbolic. 
\end{cor}

This corollary implies that the connected components of $\mathcal{D}$ are only composed of vertically expanding (or hyperbolic) elements, so we can call these sets {vertically expanding components} (or {hyperbolic components}) of the family $\Sk{p}{d}$. We mention that a connected component of $\mathcal{D}'$ is also a connected component of $\mathcal{D}$, thus the connected components of $\mathcal{D}'$ are also vertically expanding components (or hyperbolic components) of the family $\Sk{p}{d}$.  
 
As mention in introduction, Astorg-Bianchi \cite{AstBian23} proved the {structural stability} of Julia sets inside stability components of the family $\Sk{p}{d}$, generalizing the results of Lyubich \cite{Lyu83} and Ma\~{n}\'{e}-Sad-Sullivan \cite{MSS83} who proved the structural stability of Julia sets inside stability components of families of rational maps. More precisely, Astorg-Bianchi proved in {\cite[Thm 5.1, Lem. 5.4]{AstBian23}} the following:

\begin{thm}\label{thm:isotopyoverJ_p}
    Let $\Lambda$ be a stability component of $\Sk{p}{d}$. If there exists $\lambda_0\in\Lambda$ such that $f_{\lambda_0}$ is vertically expanding (resp. hyperbolic), then $f_{\lambda}$ is vertically expanding (resp. hyperbolic) for every $\lambda\in\Lambda$. 
    
    In this case, for a sufficiently small neighborhood $U\subset\Lambda$ of $\lambda_0$, there exists a unique continuous map $h:U\times J(f_{\lambda_0})\to \cmplex^2$ of the form $h_\lambda(z,w)=(z,k_{\lambda}(z,w))$ such that for every $\lambda\in U$, $h_{\lambda}:J(f_{\lambda_0})\to J(f_{\lambda})$ is a homoemorphism and $h_{\lambda}\circ f_{\lambda_0}=f_{\lambda}\circ h_{\lambda}$ on $J(f_{\lambda_0})$.
\end{thm}

It is not difficult to see that for any $z\in J_p$, $\lambda\mapsto \sum_{c\in C_{\lambda,z}}G_{\lambda,z}(c)$ is pluriharmonic on $\mathcal{D}$, thus it implies $\mathcal{D}\subset\cmplex^{D_d}\backslash\mathrm{Bif}$ according to Theorem \ref{thm:T_bif=intT_bif_zmu_p(z)}. So, the connected components of $\mathcal{D}$ are components of stability of the family $\Sk{p}{d}$. Thus this theorem can be applied to the connected components $\Lambda\in\pi_0(\mathcal{D}')$, a fact that we will use in Sec. \ref{sec:proofofcorollary1.3} to prove Corollary \ref{cor:abcmap}.

\section{Accumulation of the bifurcation locus at infinity}

\subsection{A coherent norm on the parameter space}

Since $J_p$ is not a finite set, we can define the following norm $R_{\bullet}$ on $\cmplex^{D_d}$:
\begin{equation}\label{eq:defnitiondeR_lambda}
    \lambda=(a_{0,0},\cdots,a_{d-2,2})\in\cmplex^{D_d}\mapsto R_\lambda:=\sup_{z\in J_p}\ \sum_{j=0}^{d-2}\ \left(\ \sum_{k=0}^{d-j}|a_{j,k}|^{(d-j)}|z|^k\ \right)^{1/(d-j)}.
\end{equation}

\begin{lemme}\label{lemma:rayoncritique}
    Let $d\geq 2$ and let $\alpha>1$. We denote $R_{\lambda}(\alpha):=R_{\lambda}^{\alpha}$.
    Then there exists $R_{\star}=R_{\star}(\alpha)>1$ such that for any $\lambda\in\cmplex^{D_d}$ satisfying $R_{\lambda}\geq R_{\star}$ we have:
    \begin{enumerate}{
        \item For all $(z_0,w_0,n_0)\in J_p\times\cmplex\times\mathbb{N}$ we have:
        $$|Q_{\lambda,z_0}^{n_0}(w_0)|\geq R_{\lambda}(\alpha)\Rightarrow|Q_{\lambda,z_0}^{n}(w_0)|\geq 2^{n-n_0}|Q_{\lambda,z_0}^{n_0}(w_0)|, \forall n\geq n_0.$$
        \item { $\exists\ 1<R_{\lambda}'(\alpha)<R_{\lambda}(\alpha)$ s.t. for any $z\in J_p:$ $Q_{\lambda,z}^{-(n+1)}(\mathbb{D}_{R_{\lambda}(\alpha)})\subset Q_{\lambda,z}^{-n}(\mathbb{D}_{R_{\lambda}'(\alpha)}),$ $\forall n\geq0$.}
        }
    \end{enumerate}
\end{lemme}
\noindent\textbf{\underline{Proof$:$}} 

\noindent\textit{1.} Denote $(z_n,w_n):=f_{\lambda}^n(z_0,w_0)$ and let us assume $|w_{n_0}|\geq R_{\lambda}^{\alpha}\geq1$ for some $n_0$. By invariance of $J_p$ and according to Definition \eqref{eq:defnitiondeR_lambda} we have that $|w_{n_0+1}|\geq |w_{n_0}|^d-\sum_{j=0}^{d-2}R_{\lambda}^{d-j}|w_{n_0}|^j$. Using $|w_0|\geq R^{\alpha}_{\lambda}$ we deduce: 
\begin{equation*}
    |w_{n_0+1}|\geq |w_{n_0}|^{d-1}\left[|w_{n_0}|-\sum_{j=0}^{d-2}\frac{R_{\lambda}^{d-j}}{|w_{n_0}|^{d-1-j}}\right]\geq |w_{n_0}|^{d-1}R_{\lambda}^{\alpha}\left[1-\frac{1-R_{\lambda}^{-(\alpha-1)(d-1)}}{R_{\lambda}^{2(\alpha-1)}(1-R_{\lambda}^{-(\alpha-1)})}\right].
\end{equation*}
Since $\alpha>1$, if $R_{\star}>1$ is large enough we have $R_{\lambda}^{\alpha}\geq 4$ and the term $[1-(\cdots)]$ here is greater than $1/2$. So, we have $|w_{n_0+1}|\geq |w_{n_0}|^{d-1}\frac{1}{2}R_{\lambda}^{\alpha}\geq 2|w_{n_0}|$ and thus $|w_{n_0+1}|\geq R_{\lambda}^\alpha$. Similarly, we prove by induction that $|w_{n}|\geq 2^{n-n_0}|w_{n_0}|,\ n\geq n_0$.\\ 

\noindent\textit{2.} Assume by contradiction that the statement is false. Then there exists $n(k)\geq0$ (which may not be a subsequence) and there exists $(w_k)_{k\geq1}$ and $(z_k)_{k\geq1}\in J_p^{\mathbb{N}^*}$ s.t. $|Q^{n(k)+1}_{\lambda,z_k}(w_k)|< R_{\lambda}^{\alpha}$ and $R_{\lambda}^{\alpha}-1/k\leq|Q^{n(k)}_{\lambda,z_k}(w_k)|$. According to the first item, we also have $|Q^{n(k)}_{\lambda,z_k}(w_{k})|<R_{\lambda}^{\alpha}$ for $k\geq1$. Observe that $|Q^{n(k)+1}_{\lambda,z_k}(w_k)|\geq (R_{\lambda}^{\alpha}-1/k)^d-R_{\lambda}^d\sum_{j=0}^{d-2}R_{\lambda}^{(\alpha-1)j}$. Using also $|Q^{n(k)+1}_{\lambda,z_k}(w_k)|<R_{\lambda}^{\alpha}$, when $k$ tends to infinity we obtain:
\begin{equation}\label{eq:dernierinegalitedeLemme22}
    \frac{1}{R_{\lambda}^{\alpha(d-1)}}\geq \left(1-\frac{1}{R_{\lambda}^{d(\alpha-1)}}\cdot\frac{1-R_{\lambda}^{(\alpha-1)(d-1)}}{1-R_{\lambda}^{\alpha-1}}\right)\ \underset{R_{\lambda}\to+\infty}{\sim}\ \left(1-\frac{1}{R_{\lambda}^{2(\alpha-1)}}\right).
\end{equation}
So if $R_{\star}=R_{\star}(\alpha)$ have been previously chosen sufficiently large, then for any $R_{\lambda}\geq R_{\star}$ the left and-side of \eqref{eq:dernierinegalitedeLemme22} is less than $0.1$ and its right hand-side is greater than $0.9$, thus we get a contradiction and the result follows.\qed

\subsection{Proof of Theorem \ref{thm:accumulationBIFalinfini}}\label{sec:proofofthm1.2}

Recall that the support of the bifurcation current is equal to $\mathrm{Bif}=\overline{\bigcup_{z\in J_p}\mathrm{Bif}_z}$, see Theorem \ref{thm:T_bif=intT_bif_zmu_p(z)}.
Our goal is to prove that $\mathrm{Bif}$ accumulates at infinity exactly on $E$ defined by \eqref{defn:EetlesE_z_version_intro}. To do so, we use the fact that $\mathrm{Bif}_z=\partial\mathcal{B}_z$ where $\mathcal{B}_z$ is the set of parameters $\lambda\in\Sk{p}{d}$ such that $(Q^{n}_{\lambda,z}(c))_{n\geq0}$ is bounded for some $c\in C_{\lambda,z}$, see again Theorem \ref{thm:T_bif=intT_bif_zmu_p(z)}. For $A\subset\cmplex^{D_d}$ let us denote $\mathrm{Acc}^{\infty}(A)\subset\mathbb{P}^{D_d-1}_{\infty}$ the cluster set at infinity of $A$. We will prove the following theorem, which is the counterpart of {\cite[Thm. 4.1]{AstBian23}} in the quadratic case. It implies in particular Theorem \ref{thm:accumulationBIFalinfini}.

\begin{thm}\label{thm:Accumulationatinfinity}
    For all $z\in J_p$, $\mathrm{Acc}^{\infty}(\mathcal{B}_z)=\mathrm{Acc}^{\infty}(\mathrm{Bif}_z)=E_z$. We have $\mathrm{Acc}^{\infty}(\mathrm{Bif})=E$.
\end{thm}

For sake of completeness, we deduce the following corollary as in \cite{AstBian23}.

\begin{cor}\label{cor:locusCiscompact}
    The locus $\mathcal{C}$ in $\mathrm{Definition\ \ref{defn:CBD}}$ is compact. 
\end{cor}
\noindent\textbf{\underline{Proof$:$}} Observe that $\mathcal{C}$ is closed by continuity of the Green function $G_{\lambda,z}(w)$, so we only need to prove that $\mathcal{C}$ is bounded. 
Assume that there exists $[\lambda]\in\mathbb{P}^{D_d-1}_{\infty}$ such that $[\lambda]\in\mathrm{Acc}^{\infty}(\mathcal{C})$. Let $z_1,\cdots,z_N$ be $N$ complex numbers distinct two by two in $J_p$, where $N>(\mathrm{degree\ of}\ z\mapsto P(\lambda,z))$. By definition we have $\mathcal{C}\subset \mathcal{B}_{z_1}\cap\cdots\cap\mathcal{B}_{z_N}$, thus $[\lambda]\in\mathrm{Acc}^{\infty}(\mathcal{C})\subset E_{z_1}\cap\cdots\cap E_{z_N}$ according to Theorem \ref{thm:Accumulationatinfinity}. So we have $P(\lambda,z_j)=0$ for $j=1,\cdots,N$, which implies that the coefficients of $P(\lambda,\cdot)$ are all equal to $0$, contradiction since $\lambda\neq0$.\qed\\

To prove Theorem \ref{thm:Accumulationatinfinity} we use the following proposition, proved below.

\begin{prop}\label{prop:tools} 
    The following points hold:
    \begin{enumerate}
        \item One has $\mathrm{Acc}^{\infty}(\mathcal{B}_z)\subset E_z$ for all $z\in J_p$.
        \item Let $z\in J_p$ such that $p^n(z)=z$ for some $n\geq1$. We denote the following set:
        $$\mathrm{Per}(n,z):=\left\{\lambda\in\cmplex^{D_d},\ \exists c\in C_{\lambda,z}\ :\ Q_{\lambda,z}^n(c)-c=0\right\}.$$
        It is an algebraic hypersurface of $\cmplex^{D_d}$ and we denote $\mathcal{P}(n,z):=\mathrm{Acc}^{\infty}(\mathrm{Per}(n,z))$.
        Then we have $\mathcal{P}(n,z)=\mathrm{Acc}^{\infty}(\mathcal{B}_z)= E_z$. 
        \item The sets $(\mathcal{B}_z)_{z\in J_p}$ and $\cup_{z\in J_p}\mathcal{B}_z$ are closed and $\mathrm{Acc}^{\infty}\left(\cup_{z\in J_p}\mathcal{B}_z\right)\subset E$.\\
    \end{enumerate}
\end{prop}
\noindent\textbf{\underline{Proof of Theorem \ref{thm:Accumulationatinfinity}$:$}} 
The assertion $\mathrm{Acc}^{\infty}(\partial\mathcal{B}_z)=E_z$, for all $z\in J_p$, implies the equality $\mathrm{Acc}^{\infty}(\mathrm{Bif})=E$. Indeed, this assertion implies $E\subset\mathrm{Bif}$, and the third item of Proposition \ref{prop:tools} implies the other inclusion. It remains to prove $\mathrm{Acc}^{\infty}(\partial\mathcal{B}_z)=\mathrm{Acc}^{\infty}(\mathcal{B}_z)=E_z$ for $z\in J_p$. 

Let us fix $z\in J_p$ and let us explain first that $E_z\subset\mathrm{Acc}^{\infty}(\mathcal{B}_z)$. By density of periodic points of $p$ in $J_p$, there exists a sequence of periodic points $z_n\in J_p$ such that $\lim_n z_n=z$. According to the second item of Proposition \ref{prop:tools}, $\mathrm{Acc}^{\infty}(\mathcal{B}_{z_n})=E_{z_n}$. Let us fix a point $P\in\mathbb{P}_{\infty}^{D_d-1}\backslash E$, we denote $\mathcal{P}$ the pencil of lines passing through $P$ in $\mathbb{P}^{D_d}$ that are not contained inside $\mathbb{P}^{D_d-1}_{\infty}$. Since $B=\cup_{z\in J_p}\mathcal{B}_z$ is a closed set such that $\mathrm{Acc}^{\infty}(B)\subset E$ (by the third item of Proposition \ref{prop:tools}), we should observe that $B\cap D$ is compact and that $\mathcal{B}_{z_n}\cap D=\mathcal{B}_{z_n}\cap B\cap D$ (and similarly replacing $z_n$ by $z$), for each line $D\in\mathcal{P}$. The compactness of $B\cap D$ then implies by Lemma \ref{lemma:variationofB_zandBif_z}:
\begin{equation}\label{eq:limsupdroite}
    \limsup_{n\to+\infty}\mathcal{B}_{z_n}\cap D\subset \mathcal{B}_z\cap D.
\end{equation}
Since this is true for any complex line $D\in\mathcal{P}$, we can deduce that:
\begin{equation}\label{eq:limsupAcc(B_z_n)subsetAcc(B_z)}
    \limsup_{n\to+\infty}\mathrm{Acc}^{\infty}(\mathcal{B}_{z_n})\subset\mathrm{Acc}^{\infty}(\mathcal{B}_z).
\end{equation}
Let us explain briefly how one can obtain \eqref{eq:limsupAcc(B_z_n)subsetAcc(B_z)} from \eqref{eq:limsupdroite}: consider a element $x\in\mathbb{P}^{D_d-1}_{\infty}$ such that $x=\lim_jx_j$ with $x_j\in\mathrm{Acc}^{\infty}(\mathcal{B}_{z_{n_j}})$ such that $(n_j)_{j\geq1}$ is a subsequence. Denote $Y_j:=\mathrm{Per}(N_j,z_{n_j})$, where $N_j$ satisfies $p^{N_j}(z_{n_j})=z_{n_j}$. Let $(D_m)_{m\geq1}$ be a family of lines of $\mathcal{P}$ such that $\limsup_m D_m\subset\mathbb{P}^{D_d-1}_{\infty}$. For each $m\geq1$, let $(y^m_j)_{j\geq1}$ be a sequence such that for each $j\geq1$, $y^m_j\in Y_j\cap D_m$, it exists since $Y_j$ is an algebraic hypersurface of $\cmplex^{D_d}$ such that $\mathrm{Acc}^{\infty}(Y_j)\subset E$ does not contain the center of the pencil $\mathcal{P}$. For each $m\geq1$, since $B\cap D_m$ is compact, we can extract a convergent subsequence from $(y^m_{j})_{j\geq1}$ such that the limit belongs to $\mathcal{B}_z$ by \eqref{eq:limsupdroite}. Then
we can select a collection of subsequences $(\varphi_k)_{k\geq1}$ such that for each $m\geq 1$, there exists $y^m\in\mathcal{B}_z$ satisfying 
$d_{\mathbb{P}^{D_d}}(y^m_{\varphi_1\circ\cdots\circ\varphi_m(j)},y^m)\leq 2^{-m}$ for all $j\geq1$. Then consider a subsequence $(\psi(m))_{m\geq1}$ such that for all $m\geq1$ and for all $1\leq j\leq m$:
$$d_{\mathbb{P}^{D_d}}\left(y^{\psi(m)}_{\varphi_1\circ\cdots\circ\varphi_{\psi(m)}(j)},x_{\varphi_1\circ\cdots\circ\varphi_{\psi(m)}(j)}\right)\leq 2^{-m}.$$
Observe that for each $m\geq1$, $d_{\mathbb{P}^{D_d}}(y^{\psi(m)},x)\leq 2^{-m+1}+d_{\mathbb{P}^{D_d}}(x_{\varphi_1\circ\cdots\circ\varphi_{\psi(m)}(m)},x)$. Since $\lim_j x_j=x$, we deduce that $\lim_my^{\psi(m)}=x$ proving $x\in\mathrm{Acc}^{\infty}(\mathcal{B}_z)$ and thus \eqref{eq:limsupAcc(B_z_n)subsetAcc(B_z)}.

Since $E_{z_n}=\mathrm{Acc}^{\infty}(\mathcal{B}_{z_n})$, Equation \eqref{eq:limsupAcc(B_z_n)subsetAcc(B_z)} implies $\limsup_nE_{z_n}\subset\mathrm{Acc}^{\infty}(\mathcal{B}_z)$. One can check that $\limsup_nE_{z_n}=E_z$, so we get $E_z\subset\mathrm{Acc}^{\infty}(\mathcal{B}_z)$. It implies $E_z=\mathrm{Acc}^{\infty}(\mathcal{B}_z)$ according to the first item of Proposition \ref{prop:tools}. To conclude, it remains to  prove that $\mathrm{Acc}^{\infty}(\partial\mathcal{B}_z)=E_z$, but since $\mathcal{B}_z$ is a closed subset of $\cmplex^{D_d}$ which accumulates at infinity exactly on $E_z$ which has empty interior, the boundary $\partial\mathcal{B}_z$ must also accumulates exactly on $E_z$.\qed\\

\newpage 

\noindent\textbf{\underline{Proof of Proposition \ref{prop:tools}$:$}} Let $1<\alpha<d$, $R_{\star}=R_{\star}(\alpha)$ be given by Lemma \ref{lemma:rayoncritique}.

\noindent\textit{1.} Let $[\lambda_0]$ be an element of $\mathbb{P}^{D_d-1}_{\infty}\backslash E_z$. There exists $\mathbb{B}$ a small ball centered at $\lambda_0$ in $\cmplex^{D_d}$ such that for all $\lambda\in 2\mathbb{B}$: $q_{\lambda,z}(c)\neq 0$ for all $c\in C_{\lambda,z}$. 
Recall that $(\lambda,w)\mapsto q_{\lambda,z}'(w)$ is homogeneous of degree $(d-1)$, thus we deduce that $tC_{\lambda,z}=C_{t\lambda,z}$ for any $t\in \cmplex$ (recall that $C_{\lambda,z}$ is the critical set of $q_{\lambda,z}$). Recall also that $(\lambda,w)\mapsto q_{\lambda,z}(w)$ is homegeneus of degree $d$. So, since $d-\alpha>0$, there exists $C>0$ and $r>0$ such that for all $\lambda\in\mathbb{B}$ and for all $|t|\geq r$, $R_{t\lambda}\geq R_{\star}$ and:
$$R_{t\lambda}^{\alpha}=|t|^{\alpha}R_{\lambda}^{\alpha}\leq C|t|^d\leq \inf_{c\in C_{\lambda,z}}|q_{t\lambda,z}(tc)|=\inf_{c\in C_{t\lambda,z}}|q_{t\lambda,z}(c)|,$$
By using Lemma \ref{lemma:rayoncritique}, we deduce that $t\lambda\not\in\mathcal{B}_z$ for all $\lambda\in\mathbb{B}$ and for all $|t|\geq r$. It implies that $[\lambda_0]\not\in \mathrm{Acc}^{\infty}(B_z)$. The result follows.\\

\noindent\textit{2.} If $R_{n,z}(\lambda):=\mathrm{Res}_w(Q^{n}_{\lambda,z}(w)-w,q'_{\lambda,z}(w))\in\cmplex[\lambda]$ then $\mathrm{Per}(n,z)=\{\lambda:R_{n,z}(\lambda)=0\}$, thus this set is an algebraic hypersurface of $\cmplex^{D_d}$. Let $T$ be a complex coordinate used to compactify $\cmplex^{D_d}$ into $\mathbb{P}^{D_d}$ i.e. the elements of $\mathbb{P}^{D_d}$ can be written $[\lambda_1:\cdots:\lambda_{D_d}:T]$ and $\mathbb{P}^{D_d-1}_{\infty}=\{T=0\}$. In particular we have:
$$\mathcal{P}(n,z) = \left\{[\lambda:T]\in\mathbb{P}^{D_d},\ \exists c\in C_{\lambda,z}\ :\ Q^{n}_{\lambda,z}(c,T)-cT^{d^n-1}=0\right\}\cap \left\{T=0\right\},$$
where $(\lambda,w,T)\mapsto Q^n_{\lambda,z}(w,T)$ is the homogeneization of $(\lambda,w)\mapsto Q^{n}_{\lambda,z}(w)$ into a homogeneous polynomial of degree $d^n$. Observe now that
$$\mathcal{P}(n,z)=\left\{[\lambda]\in\mathbb{P}^{D_d-1}_{\infty},\ \exists c\in C_{\lambda,z}\ :\ Q^{n}_{\lambda,z}(c,0)=0\right\}.$$
We are going to compute explicitly $Q^{n}_{\lambda,z}(w,0)$ for each $[\lambda]$ and each $w\in\cmplex$. 

Let us denote for each $k\geq0$, $W_k(\lambda,w):=Q^{k}_{\lambda,z}(w)\in\cmplex[\lambda,w]$. $W_k(\lambda,w)$ can be homogenized into a homogeneous polynomial $W_k(\lambda,w,T)$ of degree $d^k$ for every $k\geq1$, and moreover we have a recurrent formula (where $z_k:=p^{k}(z)$):
\begin{equation}\label{eq:formulerecurrencepolynomehomogeneise}
    W_{k+1}(\lambda,w,T)=W_k(\lambda,w,T)^d + \sum_{j=0}^{d-2}A_{\lambda,j}(z_k)\cdot\left(W_k(\lambda,w,T)\right)^j\cdot T^{(d^k-1)(d-j)},\ k\geq1,
\end{equation}
where $A_{\lambda,j}(z):=\sum_{k=0}^{d-j}a_{j,k}^{d-j}z^k$ with $\lambda=(a_{j,k})$. For $k=1$ we have $W_1(\lambda,w)=w^d+\sum_{j=0}^{d-2}A_{\lambda,j}(z)w^j$ which is homogeneous in $(\lambda,w)$ of degree $d$, and $W_2(\lambda,w)=W_1(\lambda,w)^d + \sum_{j=0}^{d-2}A_{\lambda,j}(z_1)W_1(\lambda,w)^j$ can thus be homogenized by $W_2(\lambda,w,T)=W_1(\lambda,w)^d + \sum_{j=0}^{d-2}A_{\lambda,j}(z_1)\cdot W_1(\lambda,w)^j\cdot T^{(d-1)(d-j)}$. One can prove similarly by induction on $k$ the recurrent formula \eqref{eq:formulerecurrencepolynomehomogeneise}. This recurrent formula implies in particular (taking $T=0$):
$$W_{k+1}(\lambda,w,0)=W_k(\lambda,w,0)^d,\ k\geq1.$$
We thus deduce $W_n(\lambda,w,0)=W_1(\lambda,w)^{d^{n-1}}$. Observe that $W_n(\lambda,w,0)=Q^n_{\lambda,z}(w,0)$ and that $W_1(\lambda,w)=q_{\lambda,z}(w)$, thus finally $Q^n_{\lambda,z}(w,0)=(q_{\lambda,z}(w))^{d^{n-1}}$. In consequence we have
\begin{align*}
    \mathcal{P}(n,z) & = \left\{[\lambda]\in\mathbb{P}^{D_d-1}_{\infty},\ \exists c\in C_{\lambda,z}\ :\ (q_{\lambda,z}(c))^{d^{n-1}}=0\right\}=E_z.
\end{align*}
At last, since $p^n(z)=z$, one can observe that for each $\lambda$ the vertical Green function $G_{\lambda,z}$ of $f$, and the vertical Green function $\Tilde{G}_{\lambda,z}$ of $f^n$, are actually equal $G_{\lambda,z}=\Tilde{G}_{\lambda,z}$. So, we deduce that $\mathrm{Per}(n,z)\subset\mathcal{B}_z$ and thus $E_z=\mathcal{P}(n,z)\subset\mathrm{Acc}^{\infty}(\mathcal{B}_z)$. The first item allows to conclude.\\

\noindent\textit{3.} Let us denote $B:=\cup_{z\in J_p}\mathcal{B}_z$ and let $\lambda\in\overline{B}$: $\lambda=\lim_n\lambda_n$, where $\lambda_n\in B$. There exists a sequence $(z_n)_n$ of elements of $J_p$ and a sequence $(c_n)_n$ such that $c_n\in C_{\lambda_n,z_n}$ and such that $G_{\lambda_n,z_n}(c_n)=0$. Up to extract a subsequence, we can assume that $z_n\to z$ with $z\in J_p$. Then $\cup_{n\geq0}C_{\lambda_n,z_n}$ must be bounded in $\cmplex$, thus $c_n\to c\in\cmplex$, up to a subsequence. The continuity of the Green function implies $G_{\lambda,z}(c)=0$, and the continuity of $(\lambda,w,z)\mapsto q'_{\lambda,z}(w)$ implies that $c\in C_{\lambda,z}$. Finally $\lambda\in\mathcal{B}_z$ by definition of $\mathcal{B}_z$ and thus $B$ is closed. By similar arguments we show that $\mathcal{B}_z$ is closed for each $z\in J_p$.

Let us prove now that $\mathrm{Acc}^{\infty}(B)\subset E$. Let $(\lambda_n)_n$ being a sequence of elements of $B$ which tends to some $[\lambda]\in\mathbb{P}^{D_d-1}_{\infty}$. For each $n$ there exists $z_n\in J_p$ and $c_n\in C_{\lambda_n,z_n}$ such that $G_{\lambda_n,z_n}(c_n)=0$. Then denoting $\underline{\lambda}:=\frac{\lambda}{||\lambda||_{\cmplex^2}}$ and $\underline{\lambda}_n:=\frac{\lambda_n}{||\lambda_n||_{\cmplex^2}}$, we have (up to a subsequence) $\underline{\lambda}_n\to \xi\underline{\lambda}$ in $\cmplex^{D_d}$, for some $\xi\in\mathbb{S}^1$. Without loss of generality we can assume that $\xi=1$. Observe that $\underline{c}_n:=\frac{c_n}{||\lambda||_{\cmplex^2}}$ satisfies $\underline{c}_n\in C_{\underline{\lambda}_n,z_n}$. Up to extract a subsequence, we can assume that $z_n\to z\in J_p$. Then $\bigcup_{n\geq0}C_{\underline{\lambda}_n,z_n}$ is bounded, so up to extract again a subsequence, we can also assume that $\underline{c}_n\to\underline{c}$ for some $\underline{c}\in\cmplex$. We denote $c:=||\lambda||_{\cmplex^2}\times \underline{c}$. We have $q_{\lambda,z}'(c)=||\lambda||_{\cmplex^2}^{d-1}\times q'_{\underline{\lambda},z}(\underline{c})=||\lambda||_{\cmplex^2}^{d-1}\lim_n q'_{\underline{\lambda}_n,z_n}(\underline{c}_n)=0$, since $c_n\in C_{\lambda_n,z_n}$ for all $n\geq0$. We have thus proved that $c\in C_{\lambda,z}$.

We prove now that $q_{\lambda,z}(c)=0$. To do so we use that $G_{\lambda_n,z_n}(c_n)=0$ which implies $|q_{\lambda_n,z_n}(c_n)|\leq R_{\lambda_n}^{\alpha}$ if $R_{\lambda_n}\geq R_{\star}$, according to Lemma \ref{lemma:rayoncritique}. We can assume $R_{\lambda_n}\geq R_{\star}$ for any $n\geq0$, since $R_{\lambda_n}\to+\infty$. By homogeneity we deduce $|q_{\underline{\lambda_n},z_n}(\underline{c}_n)|\leq R_{\underline{\lambda_n}}^{\alpha}\times \frac{||\lambda_n||_{\cmplex^2}^{\alpha}}{||\lambda_n||_{\cmplex^2}^d}\to 0$ since $d>\alpha$. 
By passing to the limit we obtain by continuity of $q$ that $q_{\underline{\lambda},z}(\underline{c})=0$, and thus $q_{\lambda,z}(c)=0$ by homogeneity. So, $P(\lambda,z)=0$ i.e. $[\lambda]\in E_z$ and we have $\mathrm{Acc}^{\infty}(B)\subset E$.\qed

\section{Hyperbolic components and partial self-coverings}

Our aim in what follows, is to show that every $\Lambda\in\pi_0(\mathcal{D}')$ contains an element $\lambda\in\Lambda$ such that $f_{\lambda}$ is a partial self-cover on $\mathcal{M}=J_p\times\cmplex$ and $\bigcap_{n\geq0}f_{\lambda}^{-n}(\mathcal{M})=J(f_{\lambda})$, such a parameter is called \textit{admissible}. This will be used in Section \ref{sec:topologicaldescriptionofJ(f)} to establish a topological description of the Julia set $J(f_{\lambda})$ by computing the iterated monodromy of $f_{\lambda}$ on a pre-image tree inside $\mathcal{M}$ (with $p(z)=z^d$). 

\subsection{Admissible parameters and the Poincaré metric}\label{sec:admissibleparameters}
 
We denote $\pi:J_p\times\cmplex\to J_p$ the first projection, and if $\mathcal{M}\subset J_p\times\cmplex$, we denote $\mathcal{M}(z):=\mathcal{M}\cap\pi^{-1}(z)$.

\begin{defn}\label{defn:admissileparameters}
    A parameter $\lambda\in\mathcal{D}$ is admissible if there exists $R>0$ such that:
    \begin{enumerate}
        \item Denoting $\mathcal{M}:=J_p\times\mathbb{D}_R$ and $\mathcal{M}_n:=f^{-n}_{\lambda}(\mathcal{M})$, we have $\mathcal{M}_{n+1}\Subset\mathcal{M}_n$ for each $n\geq0$.
        \item $f_{\lambda}:\mathcal{M}_1\to\mathcal{M}$ is a $d^2-$fold covering map.
        \item $J(f_{\lambda})=\bigcap_{n\geq0}\mathcal{M}_n$.
    \end{enumerate}
    We will also say that $\lambda$ is $R-$admissible.
\end{defn}

\begin{prop}\label{prop:propertiesofpartialselfcovers}
    Let $\lambda\in\mathcal{D}$ and assume it is $R-$admissible for some $R>0$. Let $\mathcal{M}=J_p\times\mathbb{D}_R$ and $\mathcal{M}_n=f^{-n}_{\lambda}(\mathcal{M})$. 
    \begin{enumerate}
        \item For each $n\geq0$, $f_{\lambda}^n:\mathcal{M}_n\to\mathcal{M}$ is a $d^{2n}-$fold covering map. 
        \item For each $n\geq0$ and each pair $(z',z)$ such that $p^n(z')=z$, the map $f_{\lambda}^n:\mathcal{M}_n(z')\to\mathcal{M}(z)$ is a $d^n-$fold covering map.
    \end{enumerate}
\end{prop}
\noindent\textbf{\underline{Proof$:$}} The first item is a direct consequence of the fact that $\mathcal{M}_1\subset\mathcal{M}$ and that $f:\mathcal{M}_1\to\mathcal{M}$ is a $d-$fold covering map. Then observe that $\mathcal{M}$ does not contain any critical value of $f^n|_{\mathcal{M}}$, thus for $z\in J_p$ we have $\#(f^{-n}(\mathcal{M})\cap\pi^{-1}(z))=d^n$. The second item then follows.\qed\\

In the sequel of the article, it is convenient to use the Poincaré (or hyperbolic) metric $d\rho_{R}(w)=\frac{2R}{R^2-|w|^2}|dw|$ on $\mathbb{D}_R$. Let us fix $\lambda\in\mathcal{D}$ and let us assume that $\lambda$ is $R-$admissible for some $R>0$, we denote again $\mathcal{M}=J_p\times\mathbb{D}_R$ and $\mathcal{M}_n=f_{\lambda}^{-n}(\mathcal{M})$. For simplicity, we drop the index $\lambda$ in the rest of the paragraph. 

Since $\mathcal{M}_1\Subset\mathcal{M}$ there exists a disc $D$ centered at $0$ such that $\mathcal{M}_1(z)\subset\{z\}\times D\subset\{z\}\times\overline{D}\subset\mathcal{M}(z)$ for any $z\in J_p$. Let $d_{\mathcal{M}}$ be the metric on $\mathcal{M}$ defined by:
\begin{equation}\label{eq:defnofd_M}
    d_{\mathcal{M}}\left((z_1,w_1),(z_2,w_2)\right):=d_{J_p}(z_1,z_2) + \rho_R(w_1,w_2),
\end{equation}
where $d_{J_p}$ is the standard metric of $\cmplex$ restricted to $J_p$. Observe that the metric $d_{\mathcal{M}}$ coincide with $\rho_R$ in each fiber of $\pi$ over $J_p$, thus there exists a constant $C>0$ such that for any $(z,w_1),(z,w_2)\in J_p\times\overline{D}$ (see for instance \cite[Thm. 4.3]{CarlGam93}):
\begin{equation}\label{eq:equivalenced_Met||.||surD}
    \frac{1}{C}|w_1-w_2|\leq d_{\mathcal{M}}\left((z,w_1),(z,w_2)\right)\leq C|w_1-w_2|.
\end{equation}
The interest of such a metric $d_{\mathcal{M}}$ is the following. Let us denote $\Gamma(z)$ the set of global sections of the covering map $\{z\}\times q_z^{-1}(\mathbb{D}_R)\to\{z\}\times\mathbb{D}_R$ defined by $(z,w)\mapsto (z,q_z(w))$. Since $q_z^{-1}(\mathbb{D}_R)\Subset\mathbb{D}_R$ the elements $g\in\Gamma(z)$ contract the Poincaré metric $d\rho_R$ (Schwarz-Pick Theorem). In consequence, for each $g\in\Gamma(z)$ we have $k(g)<1$, where $k(g)$ is the Lipchitz constant:
$$k(g):=\mathrm{Lip}\left[g:\left(\{z\}\times\overline{D},d_{\mathcal{M}}\right)\to\left(\{z\}\times\overline{D},d_{\mathcal{M}}\right)\right].$$
Using the fact that $f$ is holomorphic and has a structure of finite covering on $\mathcal{M}=J_p\times\mathbb{D}_R$, we can cover $J_p$ 
\begin{equation}\label{eq:recouvrementduJuliaJ_p}
    J_p=\bigcup_{j=1}^NI_j
\end{equation}
by a finite number of open sets $I_j$ (if $J_p=\mathbb{S}^1$, $I_j$ is a small interval) such that the following result holds: 
\begin{prop}\label{lemma:contractiondanslesfibres}
    Let $\lambda\in\mathcal{D}$ and assume it is $R-$admissible. Using the same notation than above, there exists a family $\underset{\ \ \ \ \ 1\leq k\leq d}{(\mathcal{G}_{jk})_{1\leq j\leq N}}$ satisfying:
    \begin{enumerate}
        \item For each $jk$, $\mathcal{G}_{jk}:I_j\times\mathbb{D}_R\to\mathcal{M}_1$ is a section of the covering $f_{\lambda}:\mathcal{M}_1\to\mathcal{M}$.
        \item For each $jk$, there exists $l_{jk}\in]0,1[$ such that for each $(z_1,w_1),(z_2,w_2)\in I_j\times\overline{D}$:
        $$d_{\mathcal{M}}\left(\mathcal{G}_{jk}(z_1,w_1),\mathcal{G}_{jk}(z_2,w_2)\right)\leq l_{jk}d_{\mathcal{M}}((z_1,w_1),(z_2,w_2)).$$
        \item Let $z\in J_p$ and let $j$ such that $p(z)\in I_j$. Then for each $g\in\Gamma(z)$, there is a unique $k$ such that $\widehat{g}_z(w)=\widehat{\mathcal{G}}_{jk}(p(z),w)$ for all $w\in\mathbb{D}_R$. The mappings $\widehat{g}_z$ and $\widehat{\mathcal{G}}_{jk}$ are defined by the following relations for all $z'\in I_j$ and for all $w\in\mathbb{D}_R$:
        $$g(z,w)=:(z,\widehat{g}_z(w))\ \mathrm{and}\ \mathcal{G}_{jk}(z',w)=:\left(\mathcal{G}_{jk}^1(z'),\widehat{\mathcal{G}}_{jk}(z',w)\right).$$
    \end{enumerate}
    In particular we have:
    \begin{equation}\label{eq:constantedeLipchitzglobale}
        k(f_{\lambda}):=\sup_{z\in J_p}\max_{g\in\Gamma(z)}k(g)\ \leq \max\left\{l_{jk},\ 1\leq j\leq N,\ 1\leq k\leq d\right\} <\ 1.
    \end{equation}
\end{prop}

\subsection{Existence of admissible parameters in \texorpdfstring{$\Lambda\in\pi_0(\mathcal{D}')$}{TEXT}}

Let us fix $\Lambda\in\pi_0(\mathcal{D}')$. Our goal is to prove the existence of an admissible parameter in $\Lambda$. To do so, we first prove that the postcritical set of $f_{\lambda}$, for adapted parameters $\lambda\in\Lambda$, does not accumulate on a neighborhood $\mathcal{M}$ of $J(f_{\lambda})$. We recall that the set $E$ is defined by \eqref{defn:EetlesE_z_version_intro}. We fix an arbitrary choice of $\alpha\in]1,d[$ and let $R_{\star}=R_{\star}(\alpha)$ and $R_{\lambda}(\alpha)$ be as in Lemma \ref{lemma:rayoncritique}. The postcritical set of a map $f_{\lambda}$ is defined by:
\begin{equation}\label{eq:postcriticalsetoff_lambda}
    P_{f_{\lambda}}:=\bigcup_{n\geq1}f^n_{\lambda}(\mathrm{Crit}\ f_{\lambda})
\end{equation}
If we assume that $J_p\cap\mathrm{Crit}(p)=\emptyset$, we observe that for any $\lambda\in\cmplex^{D_d}$:
\begin{equation}\label{eq:postcritic_inter_JxC}
    (J_p\times\cmplex)\cap P_{f_{\lambda}}=\bigcup_{n\geq1}f_{\lambda}^n(C_{\lambda})\ \mathrm{with}\ C_{\lambda}=\bigcup_{z\in J_p}\{z\}\times C_{\lambda,z}.
\end{equation}

\begin{lemme}\label{lemma:criticalpointsescape}
    Assume that $J_p\cap\mathrm{Crit}(p)=\emptyset$, and let us fix $\Lambda\in\pi_0(\mathcal{D}')$. Let $\lambda_0\in\Lambda$ such that $R_{\lambda_0}\geq R_{\star}$ and assume that $[\lambda_0]\not\in E$. 
    \begin{enumerate}
        \item For any $t\in\cmplex$ of modulus $|t|\geq 1$ we have $t\lambda_0\in\Lambda$. 
        \item There exists $C=C(\lambda_0)>0$ such that for any $t\in \cmplex$ and any $z\in J_p$: 
        $$\frac{1}{C}|t|^d\leq \min_{c\in C_{t\lambda_0,z}}|q_{t\lambda_0,z}(c)|\leq \max_{c\in C_{t\lambda_0,z}}|q_{t\lambda_0,z}(c)|\leq C|t|^d.$$
        \item There exists $r(\lambda_0)>1$ such that for any $|t|\geq r(\lambda_0)$:
        $$\left(J_p\times\mathbb{D}_{R_{t\lambda_0}(\alpha)}\right)\cap \overline{P_{f_{t\lambda_0}}}=\emptyset.$$
    \end{enumerate}
\end{lemme}
\noindent\underline{\textbf{{Proof$:$}}} \ \\
\noindent{\textit{1.}} Let $|t|\geq1$, $z\in J_p$ and $c\in C_{\lambda,z}$. Let us prove that for each $n\in\mathbb{N}$ the following property holds:
\begin{equation}\label{eq:propP(N)}
    \left(|Q^{n}_{\lambda_0,z}(c)|\geq R_{\lambda_0}^{\alpha}\right)\Longrightarrow \left(|Q^{n}_{t\lambda_0,z}(tc)|\geq R_{t\lambda_0}^{\alpha}\right).
\end{equation}
For $n\in\mathbb{N}$, a direct computation using that $(\lambda,w)\mapsto q_{\lambda,z}(w)$ is homogeneous of degree $d$ shows:
$$Q^{n}_{t\lambda_0,z}(tc) = t^dq_{\lambda_0,z_{n-1}}\left(t^{d-1}q_{\lambda_0,z_{n-2}}\left(\cdots t^{d-1}q_{\lambda_0,z}(c)\cdots\right)\right),$$
where $z_k:=p^{k}(z)$. So, if we assume that $|Q^{n}_{t\lambda_0,z}(tc)|< R_{t\lambda_0}^{\alpha}=t^{\alpha}R_{\lambda_0}^{\alpha}$ then Lemma \ref{lemma:rayoncritique} ensures successively that (using $|t|\geq1$ and $\alpha<d$):
$$\left|q_{\lambda_0,z_{n-2}}\left(\cdots t^{d-1}q_{\lambda_0,z}(c)\cdots\right)\right|<R_{\lambda_0}^{\alpha}\ ,\ \cdots\ ,\ |q_{\lambda_0,z}(c)|<R_{\lambda_0}^{\alpha}\ ,\ |c|<R_{\lambda_0}^{\alpha}.$$
Thus \eqref{eq:propP(N)} is true by contraposition for each $n\in\mathbb{N}$.

Property \eqref{eq:propP(N)} shows that, for any $z\in J_p$ and for any $c\in C_{\lambda_0,z}$, if the sequence $(Q^n_{\lambda_0,z}(c))_{n}$ is not bounded then the sequence $(Q^n_{t\lambda_0,z}(tc))_{n}$ is also not bounded (use again Lemma \ref{lemma:rayoncritique}). It proves that if $\lambda_0\in\mathcal{D}$ then $t\lambda_0\in\mathcal{D}$ by definition of $\mathcal{D}$. At last, since this is true for any $|t|\geq1$, the element $t\lambda_0$ belongs to the same connected component of $\lambda_0$ inside $\mathcal{D}$, thus in particular $t\lambda_0\in\Lambda$ since connected components of $\mathcal{D}'$ are also connected components of $\mathcal{D}$.\\ 

\noindent\textit{2.} By compactness of $J_p$ and by a continuity argument for the roots of polynomial mappings, one can show that $C_{\lambda_0}$ is compact, similarly as we proved that the set $\cup_{z\in J_p}\mathcal{B}_z$ is closed in Proposition \ref{prop:tools}. Thus the continuous map $(z,c)\in C_{\lambda_0}\mapsto |q_{\lambda_0,z}(c)|$ attains its minimum and its maximum, and moreover, for any $z\in J_p$, $c\in C_{\lambda_0,z}$, we have $|q_{\lambda_0,z}(c)|>0$ since $[\lambda_0]\not\in E$, see \eqref{eq:discriminantofq_lambdaz_version_intro} and \eqref{defn:EetlesE_z_version_intro}. Thus there exists a constant $C>0$ such that for any $z\in J_p$ (use that $(\lambda,w)\mapsto q_{\lambda,z}(w)$ is homogeneous of degree $d$):
\begin{equation}\label{eq:minvaleurscritiques}
    \frac{1}{C}|t|^d\leq \min_{c\in C_{\lambda_0,z}}|q_{t\lambda_0,z}(tc)|\leq \max_{c\in C_{\lambda_0,z}}|q_{t\lambda_0,z}(tc)|\leq C|t|^d.
\end{equation}
We conclude by observing that $C_{t\lambda_0,z}=t\times C_{\lambda_0,z}$.\\

\noindent\textit{3.} Let $r(\lambda_0)>1$ be large enough such that for any $|t|\geq r(\lambda_0)$ we have $C^{-1}|t|^d\geq \textcolor{black}{2}|t|^\alpha R_{\lambda_0}^{\alpha}$. It exists since $\alpha<d$. Then for any $|t|\geq r(\lambda_0)$ and for any $z\in J_p$, using \eqref{eq:minvaleurscritiques} we deduce:
$$\textcolor{black}{2}R_{t\lambda_0}^{\alpha}\leq \frac{1}{C}|t|^{d}\leq \min_{c\in C_{\lambda_0,z}}|q_{t\lambda_0,z}(tc)|.$$
According to Lemma \ref{lemma:rayoncritique}, we deduce that for any $n\geq 1$, for any $z\in J_p$ and for any $c\in C_{t\lambda_0,z}$:
$$|Q^{n}_{t\lambda_0,z}(c)|\geq 2R_{t\lambda_0}^{\alpha}.$$ 
Using \eqref{eq:postcritic_inter_JxC} we deduce finally that $(J_p\times\cmplex)\cap P_{f_{t\lambda_0}}\subset\cmplex^{2}\backslash \mathbb{D}_{R_{t\lambda_0}^{\alpha}}$ and the result follows.\qed\\

We can now prove the existence of admissible parameters in $\Lambda$ under the assumption $J_p\cap\mathrm{Crit}(p)=\emptyset$:

\begin{prop}\label{prop:lambdaisadmissible}
    Assume $J_p\cap\mathrm{Crit}(p)=\emptyset$, and let us fix $\Lambda\in\pi_0(\mathcal{D}')$. Let $\lambda_0\in\Lambda$ such that $[\lambda_0]\notin E$ and $R_{\lambda_0}\geq R_{\star}$. Let $r(\lambda_0)$ be as in $\mathrm{Lemma\ \ref{lemma:criticalpointsescape}}$. Then $t\lambda_0$ is $R_{t\lambda_0}(\alpha)-$admissible for every $|t|\geq r(\lambda_0)$. 
\end{prop}
\noindent\textbf{\underline{Proof$:$}} Let $|t|\geq r(\lambda_0)$. According to the first item of Lemma \ref{lemma:criticalpointsescape}, $\lambda:=t\lambda_0\in\Lambda$. Let us denote $\mathcal{M}=J_p\times\mathbb{D}_{R_{\lambda}(\alpha)}$ and $\mathcal{M}_n=f_{\lambda}^{-n}(\mathcal{M})$. According to the second item of Lemma \ref{lemma:rayoncritique}, we have for each $n\geq0$:
\begin{align*}
    \mathcal{M}_{n+1} = \bigcup_{z\in J_p} \{z\}\times Q_{\lambda,z}^{-(n+1)}(\mathbb{D}_{R_{\lambda}(\alpha)}) & \subset \bigcup_{z\in J_p}\{z\}\times Q_{\lambda,z}^{-n}(\mathbb{D}_{R_{\lambda}'(\alpha)})\\
               & \Subset\bigcup_{z\in J_p}\{z\}\times Q_{\lambda,z}^{-n}(\mathbb{D}_{R_{\lambda}(\alpha)})=\mathcal{M}_{n}.
\end{align*}
Since $f_{\lambda}$ is a polynomial mapping of topological degree $d^2$, the map $f_{\lambda}:\cmplex^2\backslash f_\lambda^{-1}\left(\overline{P_{f_{\lambda}}}\right)\to\cmplex^2\backslash\overline{P_{f_{\lambda}}}$ is a $d^2-$fold covering map, and according to Lemma \ref{lemma:criticalpointsescape} we have $\mathcal{M}\cap\overline{P_{f_{\lambda}}}=\emptyset$, thus by restriction $f_{\lambda}:\mathcal{M}_1\to\mathcal{M}$ is a $d^2-$fold covering map. 

It remains to prove that $J(f_{\lambda})=\bigcap_{n\geq0}\mathcal{M}_n$. By using Lemma \ref{lemma:rayoncritique}, observe that $\bigcap_{n\geq0}\mathcal{M}_n=\bigcup_{z\in J_p}\{z\}\times K_z(f_{\lambda})$. So it sufficient to prove that:
$$J(f_{\lambda})=\bigcup_{z\in J_p}\{z\}\times K_z(f_{\lambda}).$$
The justification of this last equality remains on the fact that $f_{\lambda}$ is vertically expanding (because $\lambda\in\mathcal{D}$), it is implicitly contained in the article \cite{jon99} of Jonsson, let us explain why. We have $J(f)=\bigcup_{z\in J_p}\{z\}\times J_z(f_{\lambda})$ by \eqref{eq:fhypimpliesJ=cupJ_z} and $J_z(f_{\lambda})=\partial K_{z}(f_{\lambda})$, thus it is sufficient to prove that ${\overset{\circ}{K}}_{z}(f_{\lambda})=\emptyset$. We argue by contradiction, assume that there exists $z\in J_p$ such that ${\overset{\circ}{K}}_{z}(f_{\lambda})$ is not empty and let $F\subset\{z\}\times {\overset{\circ}{K}}_{z}(f_{\lambda})$ be a connected component. By \cite[Cor. 3.6]{jon99} there exists $(n_k)_k$ a subsequence such that 
$$\lim_{k\to+\infty}\mathrm{dist}\left(f_{\lambda}^{n_k}(F),P_{f_{\lambda}}\right)=0.$$
Observe that $f_{\lambda}^{n_k}(F)\subset\overline{\mathcal{M}}_1$ for each $k$. So we deduce that $\mathrm{dist}\left(\overline{\mathcal{M}}_1,\overline{P}_{f_{\lambda}}\right)=0$, but $\mathcal{M}_1\Subset\mathcal{M}$ and by Lemma \ref{lemma:criticalpointsescape} we have $\mathcal{M}\cap \overline{P}_{f_{\lambda}}=\emptyset$, contradiction.\qed

\section{Iterated monodromy of admissible parameters}\label{sec:IMG}

We assume in this section that $p(z)=z^d$ and thus $J_p=\mathbb{S}^1$. Let us fix $\Lambda\in\mathcal{D}'$ and let $\lambda\in\Lambda$ be a $R-$admissible parameter for some $R>0$ in the sense of Definition \ref{defn:admissileparameters}, it exists according to Proposition \ref{prop:lambdaisadmissible}. We denote $f=f_{\lambda}$, $\mathcal{M}=\mathbb{S}^1\times\mathbb{D}_R$ and $\mathcal{M}_n=f^{-n}(\mathcal{M})$. We also denote $\pi:\mathbb{S}^1\times\cmplex\to\mathbb{S}^1$ the first projection, and $\mathcal{M}_n(z)=\mathcal{M}_n\cap\pi^{-1}(z)$. In particular, we denote $f_1:\mathcal{M}_1(1)\to\mathcal{M}(1)$ the $d-$fold covering induced by the restriction of $f:\mathcal{M}_1\to\mathcal{M}$ on $\mathcal{M}_1(1)$. 

\subsection{Iterated monodromy}\label{sec:IteratedMonodromy}

Let us defined the following two sets which are pre-image trees of $f$ and $f_1$: 
$$\mathcal{X}:=\bigcup_{n=0}^{+\infty}f^{-n}(t_0)\ \supset\ X:=\bigcup_{n=0}^{+\infty}f^{-n}_1(t_0),$$ 
where $t_0:=(1,0)$ is a base point in $\mathcal{M}$. We also denote $\mathcal{X}_n=f^{-n}(t_0)$ and $X_n=f^{-n}_1(t_0)$. Our aim in this section is to introduce and to compute explicitly an action of (a subset of) $\pi_1(\mathcal{M},t_0)$ on $X$. It is well known that $\pi_1(\mathcal{M},t_0)$ acts by monodromy on $\mathcal{X}$, let us recall briefly this fact. One can define the group $\mathrm{Aut}(\mathcal{X})$ of automorphisms of $\mathcal{X}$ as the subgroup of bijections $g:\mathcal{X}\to\mathcal{X}$ such that $g\in\mathfrak{S}(\mathcal{X}_n)$ and $g\circ f= f\circ g$ on $\mathcal{X}_n$ for each $n\geq1$. The monodromy action of $\pi_1(\mathcal{M},t_0)$ on $\mathcal{X}$ is then define by the group homomorphism:
\begin{equation}\label{eq:morphismedegroupesX}
    \phi:
    \left\{
        \begin{array}{ll}
          \pi_1(\mathcal{M},t_0)  & \longrightarrow\ \mathrm{Aut}(\mathcal{X})\\
          \ \ \ \ \ \mathrm{[}\alpha\mathrm{]} & \longmapsto\ \left(x\mapsto f^{-n}(\alpha)[x],\ \mathrm{if}\ x\in\ \mathcal{X}_n\right)\\
        \end{array}
    \right.
\end{equation}
where $f^{-n}(\alpha)[x]:=f^{-n}(\alpha)_x(1)$ and $f^{-n}(\alpha)_x:[0,1]\to\mathcal{M}$ is the lift starting at $x$ of $\alpha$ by the covering $f^n:\mathcal{M}_n\to\mathcal{M}$. Note that $\pi_1(\mathcal{M},t_0)\simeq\pi_1(\mathbb{S}^1,1)\simeq\mathbb{Z}$ is abelian, thus \eqref{eq:morphismedegroupesX} is a morphism and not an antimorphism. In what follows, it is convenient to fix an explicit isomorphism between $\pi_1(\mathcal{M},t_0)$ and $\mathbb{Z}$:
\begin{equation}\label{eq:groupepi_1ducercle}
    \begin{array}{ll}
          \mathbb{Z} & \longrightarrow\ \pi_1(\mathcal{M},t_0)\\
                    n & \longmapsto\ \mathrm{[}\gamma_n\mathrm{]}\\
    \end{array}
\end{equation}
where $\gamma_n:[0,1]\to\mathcal{M}$ is the path defined by $\gamma_n(t):=(e^{2i\pi t n},0)$ for every $t\in[0,1]$.

The group generated by the loops of $\pi_1(\mathcal{M},t_0)$ acting non trivially on $\mathcal{X}$ is called the \textit{iterated monodromy group} of the covering $f:\mathcal{M}_1\to\mathcal{M}$ at $t_0$, and it is denoted:
$$\mathrm{IMG}(f,t_0):=\frac{\pi_1(\mathcal{M},t_0)}{\mathrm{Ker}\ \phi}.$$
This group has been introduced by V. Nekrashevych, we refer to the survey \cite{Nekra11}. We mention that this group does not depend on the choice of $t_0\in\mathcal{M}$ (see \cite{Nekra11}) but our choice of $t_0=(1,0)$ is natural in our setting. 

Actually, the set of real interest in our work is the set $X\subset\mathcal{X}$, but the action of $\mathrm{IMG}(f,t_0)$ on $\mathcal{X}$ does not preserve $X$. So, we consider below a collection of subgroups, denoted $\mathbb{A}$, composed of all the elements that preserve certain levels $X_n$ of the tree $X$. 

\begin{lemme}
    Let $n\in\mathbb{N}$, let $x\in X_n$ and let $[\alpha]=[\gamma_m]\in\pi_1(\mathcal{M},t_0)$. Then $f^{-n}(\alpha)[x]\in X_n$ if and only if $m\in d^n\mathbb{Z}$. In particular, the group $\mathbb{A}_{d^n}:=\{[\gamma_{d^nm}],\ m\in\mathbb{Z}\}$ acts on $X_k$ for every $0\leq k\leq n$.
\end{lemme}
\noindent\textbf{\underline{Proof$:$}} We can assume $\alpha=\gamma_m$. Let us defined $y:=f^{-n}(\gamma_{m})[x]$. The path $f^{-n}(\gamma_{m})_{x}$ can be written $f^{-n}(\gamma_{m})_{x}(t)=(Z_t,W_t)$, where $(Z_t,W_t)$ is a path inside $\mathcal{M}$. We have $p^n(Z_t)=\pi\circ f^n(Z_t,W_t)=\pi\circ\gamma_{m}(t)=e^{2i\pi tm}=p^n\circ\pi\circ\gamma_{m/d^n}(t)$,
and $Z_0=1$ since $(Z_0,W_0)=x\in\mathcal{M}(1)$. By uniqueness of lifts by the $d^n-$fold covering map $p^n:\mathbb{S}^1\to\mathbb{S}^1$, we deduce that $Z_t=e^{2i\pi tm/{d^n}}$. Then observe that $y=(Z_1,W_1)$ by construction of $y\in\mathcal{X}_n$. Thus $y\in X_n$ if and only if $Z_1=1$, that is equivalent to say that $m/d^n\in\mathbb{Z}$.\qed\\

Let us now define the collection $\mathbb{A}$ of subgroups of $\pi_1(\mathcal{M},t_0)$ by:
\begin{equation}\label{eq:filtrationA}
    \mathbb{A}:=\left\{\mathbb{A}_{d^n},\ n\in\mathbb{N}\right\}.
\end{equation}
This set is exactly the elements of $\mathrm{IMG}(\mathcal{M},t_0)$ which preserve some levels of the tree $X$. More precisely, the group $\mathbb{A}_{d^n}\in\mathbb{A}$ acts by monodromy on the finite tree $\bigcup_{k=0}^nX_k$ for each $n\geq0$. By abuse of language, we will say that $\mathbb{A}$ acts on the tree $X$ by monodromy. Our task is now to give an explicit formula of the monodromy action of $\mathbb{A}$ on $X$. To do so, we follow Nekrashevych work's \cite{Nekra11} which furnish the possibility to compute this monodromy action by using a recurrent formula, see Proposition \ref{prop:formulederecurrence} and Corollary \ref{cor:lafomulederecurrencedeNekrashevych}. This is the major tool used in Section \ref{sec:topologicaldescriptionofJ(f)} to give a topological description of $J(f)$. 

We denote $\mathcal{A}:=\{1,\cdots,d\}$ the alphabet on $d$ letters. It generates a free monoid $\mathcal{A}^*$ formed of finite words $a_n\cdots a_1$ (written from right to left) and equipped with a natural concatenation low $(a_n\cdots a_1)\cdot(b_m\cdots b_1)=a_n\cdots a_1 b_m\cdots b_1$, we also add $\emptyset$ in $\mathcal{A}^*$ to have a neutral element for the concatenation law. We denote $\mathcal{A}_n$ the set of words of length $n$, with $\mathcal{A}_0=\{\emptyset\}$. 
Let us arbitrarily label $x_1,\cdots,x_d$ the elements of $X_1$, and let us choose smooth paths $l_j:[0,1]\to\mathcal{M}(1)$ joining $t_0$ to $x_j$ for $1\leq j\leq d$.

\begin{defn}\label{defn:codageL}
    We denote $L:\mathcal{A}^*\to X$ the map defined inductively on each level $\mathcal{A}_n,\ n\geq0$, by the recurrent formulas: 
    \begin{equation}\label{eq:recurrentformulapourL}
    \left\{
        \begin{array}{ll}
          n=0:     & L(\emptyset)=t_0\\
          n=1:     & L(a)=x_a\\
          n\geq 2: & L(a_n\cdots a_1)=f^{-(n-1)}(l_{a_n})[L(a_{n-1}\cdots a_1)]\\
        \end{array}
    \right.
    \end{equation}
\end{defn}
\noindent Observe that $f_1^{-(n-1)}(l_{a_n})_{L(a_{n-1}\cdots a_1)}$ is also a lift by $f^{n-1}$ of $l_{a_n}$ starting at $L(a_{n-1}\cdots a_1)$. So, by uniqueness of lifts by $f^{n-1}:\mathcal{M}_{n-1}\to\mathcal{M}$ starting at the same point, we have $f^{-(n-1)}(l_{a_n})[L(a_{n-1}\cdots a_1)]=f_1^{-(n-1)}(l_{a_n})[L(a_{n-1}\cdots a_1)]\in X$. In particular, it guarantees that $L$ takes its values in $X$ and not in $\mathcal{X}\backslash X$. 

The uniqueness in the path lifting property of covering mappings allows to prove:

\begin{lemme}\label{lemma:treeisomorphismLieLcode}
    The map $L:\mathcal{A}^*\to X$ satisfies:
    \begin{enumerate} 
        \item For $a\in\mathcal{A}_1$ we have $f(L(a))=L(\emptyset)=t_0$.
        \item For $n\geq 2$ and for $a_n\cdots a_1\in\mathcal{A}_n$ we have $f(L(a_n\cdots a_1))=L(a_n\cdots a_{2})$. 
        \item $L:\mathcal{A}_n\to X_n$ is bijective for any $n\geq0$, and thus $L:\mathcal{A}^*\to X$ is bijective.
    \end{enumerate}
\end{lemme}

\begin{defn}\label{defn:notationspratiques}
    For $n\geq1$ and for $[\alpha]\in\mathbb{A}_{d^n}$ we denote:
    \begin{enumerate}
        \item For $a_n\cdots a_1\in \mathcal{A}_n$, we put $(a_n\cdots a_1)^{\alpha}:=L^{-1}\left(f^{-n}(\alpha)[L(a_n\cdots a_1)]\right)$. 
        \item For $x\in X_n$, we put $x^{\alpha}:=f^{-n}(\alpha)[x]$.
        \item[3.] We put for each $a\in\mathcal{A}_1$, $\langle\alpha,a\rangle:=l_a\alpha_a l_b^{-1}$, where $b:=a^{\alpha}$ and $\alpha_a:=f^{-1}(\alpha)_{x_a}$. Observe that $\langle\alpha,a\rangle$ is a loop based at $t_0$ inside $\mathcal{M}$.
    \end{enumerate}
\end{defn}
\noindent With these notations we have for any $v\in\mathcal{A}^*$:
\begin{equation}\label{eq:formuleL(v^alpha)=L(v)^{alpha}}
    L(v^{\alpha})=L(v)^\alpha.
\end{equation}
We can then deduce a recurrent formula for the monodromy, see \cite[Proposition 2.1]{Nekra11}.

\begin{prop}\label{prop:formulederecurrence}
    For any $n\geq2$, $a_n\cdots a_1\in \mathcal{A}_n$ and $[\alpha]\in\mathbb{A}_{d^n}$ we have:
    $$(a_n\cdots a_1)^{\alpha}=a_n^{\alpha}\cdot (a_{n-1}\cdots a_{1})^{\langle\alpha,a_n\rangle}.$$
\end{prop}

\subsection{Computation of the action of \texorpdfstring{$\mathbb{A}$}{TEXT} on \texorpdfstring{$X$}{TEXT}}\label{sec:computationmonodromyactionofAonX}

\begin{defn}\label{defn:definitiondeS}
    Let $S\in\mathfrak{S}_d$ be the permutation of $\mathcal{A}=\{1,\cdots,d\}$ defined by:
    $$S(a):=a^{\gamma_d},\ \forall a\in\mathcal{A}.$$
\end{defn}

\begin{lemme}\label{lemma:toursdanslabase}
    The permutation $S\in \mathfrak{S}_d$ satisfies that for any $m\in\mathbb{Z}$:
    \begin{enumerate}
        \item For $n\geq0$ and for $a\in\mathcal{A}_1$, we have $a^{\gamma_{d^{n+1}m}}=S^{d^nm}(a)$. 
        \item For $n\geq1$ and for $a\in \mathcal{A}_1$, we have $[\langle\gamma_{d^nm}, a\rangle] = [\gamma_{d^{n-1}m}]$ in $\pi_1(\mathcal{M},t_0)$ (see $\mathrm{Definition\ \ref{defn:notationspratiques}}$).
    \end{enumerate}
\end{lemme}
\noindent\textbf{\underline{Proof$:$}} Since $[\gamma_{d^{n+1}m}]=[\gamma_d]^{d^{n}m}$, we deduce that $a^{\gamma_{d^{n+1}m}}=S^{d^nm}(a)$. Let us prove now the second item. Let us defined $b:=a^{\gamma_{d^nm}}$ with $m\in\mathbb{Z}$. By Definition \ref{defn:notationspratiques} we have $\langle\gamma_{d^nm},a\rangle=l_a\cdot(f^{-1}(\gamma_{d^nm})_{L(a)})\cdot l_{b}^{-1}$. The path $f^{-1}(\gamma_{d^nm})_{L(a)}$ can be written $f^{-1}(\gamma_{d^nm})_{L(a)}(t)=(Z_t,W_t)$, where $(Z_t,W_t)$ is a path inside $\mathcal{M}$. We have $p(Z_t)=\pi\circ f(Z_t,W_t)=\pi\circ\gamma_{d^nm}(t)=e^{2i\pi td^nm}=p\circ\pi\circ\gamma_{d^{n-1}m}(t)$ and $Z_0=1$ since $(Z_0,W_0)=L(a)\in \mathcal{M}(1)$. By uniqueness of lifts we deduce that $Z_t=e^{2i\pi td^{n-1}m}$. Then observe that $l_a(t)=(1,w_a(t))$ and similarly $l_b(t)=(1,w_{b}(t))$, thus we have explicitly:
$$l_a\cdot(f^{-1}(\gamma_{d^nm})_{L(a)})\cdot l_{b}^{-1}(t)=\left(\left(C_1\cdot Z_{\bullet}\cdot C_1\right)(t),\left(w_{a}\cdot W_{\bullet}\cdot{w_{b}}^{-1}\right)(t)\right),\ \forall t\in[0,1],$$
where $C_1$ is the constant path equal to $1$, and where $(w_{a}\cdot W_{\bullet}\cdot w_{b}^{-1})$ is a loop centered at $0$ homotopic inside $\mathbb{D}_{R}$ to $C_0$ the constant path equal to $0$. Let $(h_t)_{t\in [0,1]}$ (resp. $(k_t)_{t\in[0,1]}$) be a homotopy fixing $1$ (resp. $0$) and deforming the loop $C_1\cdot Z_{\bullet}\cdot C_1$ inside $\mathbb{S}^1$ (resp. $w_a\cdot W_{\bullet}\cdot w_b^{-1}$ inside $\mathbb{D}_{R}$) into the loop $Z_{\bullet}$ (resp. $C_0$). Then $H_t(s):=(h_t(s),k_t(s))$ is a homotopy fixing $t_0=(1,0)$ and deforming the loop $l_a\cdot(f^{-1}(\gamma_{d^nm})_{L(a)})\cdot l_{b}^{-1}$ into the loop $(Z_{\bullet},0)=\gamma_{d^{n-1}m}$ inside $\mathcal{M}$. We have then proved that $[\langle\gamma_{d^nm}, a\rangle] = [\gamma_{d^{n-1}m}]$ in $\pi_1(\mathcal{M},t_0)$.\qed\\

\begin{cor}\label{cor:lafomulederecurrencedeNekrashevych}
    We have for any $n\geq1$ and $m\in\mathbb{Z}$, and for any $a_n\cdots a_1\in\mathcal{A}_n$:
    $$(a_n\cdots a_1)^{\gamma_{d^{n}m}}=S^{d^{n-1}m}(a_n)\cdots S^{d^0m}(a_1).$$
\end{cor}
\noindent\textbf{\underline{Proof$:$}} By Proposition \ref{prop:formulederecurrence} we have $(a_n\cdots a_1)^{\gamma_{d^{n}m}}=a_n^{\gamma_{d^{n}m}}(a_{n-1}\cdots a_1)^{\langle\gamma_{d^{n}m},a_n\rangle}$ and, according to Lemma \ref{lemma:toursdanslabase}, we have $(a_{n-1}\cdots a_1)^{\langle\gamma_{d^{n}m},a_n\rangle}=(a_{n-1}\cdots a_1)^{\gamma_{d^{n-1}m}}$. Thus we have 
$$(a_n\cdots a_1)^{\gamma_{d^{n}m}}=a_n^{\gamma_{d^{n}m}}(a_{n-1}\cdots a_1)^{\gamma_{d^{n-1}m}}.$$ 
So by induction on $n\geq1$ we deduce that
$(a_n\cdots a_1)^{\gamma_{d^{n}m}}=a_n^{\gamma_{d^{n}m}}a_{n-1}^{\gamma_{d^{n-1}m}}\cdots a_1^{\gamma_{d^{n-(n-1)}m}}$. The first item of Lemma \ref{lemma:toursdanslabase} then allows to conclude.\qed

\section{Topological description of \texorpdfstring{$J(f)$}{TEXT}}\label{sec:topologicaldescriptionofJ(f)}

The present section is composed of three theorems: Theorem \ref{thm:CompConnexesJ(f)}, Theorem \ref{thm:Returnmapisahomeo} and Theorem \ref{thm:Kisotopicf^-1(S^1times0)}. We present the content of these theorems in the heading of this section, the subsections below are devoted to prove them.

Let $\Lambda\in\pi_0(\mathcal{D}')$ and let $\lambda\in\Lambda$ be a $R-$admissible parameter, we use the same notations and objects than the one introduced in Section \ref{sec:IMG}, and we drop the index $\lambda$ in $f_{\lambda}=f$. Recall that $J(f)=\bigcap_{n\in\mathbb{N}}\mathcal{M}_n$ with $\mathcal{M}_n=f^{-n}(\mathcal{M})$ and $\mathcal{M}=\mathbb{S}^1\times\mathbb{D}_R$. Putting $t_0=(1,0)$ we have introduced $X=\bigcup_{n=0}^{+\infty}f_1^{-n}(t_0)$ and $\mathbb{A}=\{\mathbb{A}_{d^n},n\in\mathbb{N}\}$ (defined in \eqref{eq:filtrationA}) acting on $X$ level by level i.e. $\mathbb{A}_{d^n}$ acts on $f^{-n}_1(t_0)$ for each $n$. In particular, the action of $\mathbb{A}_d$ on $f_1^{-1}(t_0)$ gives a permutation $S\in\mathfrak{S}_d$, see Definition \ref{defn:definitiondeS}, which depends on the labeling $\{x_1,\cdots,x_d\}=f^{-1}_1(t_0)$.\\ 

We focus in this section on the Cantor set $\mathcal{C}_{\mathcal{A}}$ associated to the alphabet $\mathcal{A}=\{1,\cdots,d\}$ defined as follows:
\begin{equation}\label{eq:CantorsetC_A}
    \mathcal{C}_{\mathcal{A}}:=\{a_1\cdots a_n\cdots\ |\ a_i\in\mathcal{A}\}\simeq \mathcal{A}^{\mathbb{N}}.
\end{equation}
Let  $\sigma:\mathcal{C}_{\mathcal{A}}\to\mathcal{C}_{\mathcal{A}}$ be the map defined by the one-sided shift:
\begin{equation}\label{eq:onesidedshiftsigma}
    \sigma(a_1a_2\cdots a_n\cdots):=(a_2\cdots a_{n+1}\cdots).
\end{equation}
We also consider a homeomorphism $h:\mathcal{C}_{\mathcal{A}}\to\mathcal{C}_{\mathcal{A}}$ defined by:
\begin{equation}\label{eq:homeoh}
    h(a_1a_2a_3\cdots a_n\cdots):=\left(S(a_1)S^d(a_2)S^{d^2}(a_3)\cdots S^{d^{n-1}}(a_n)\cdots\right).
\end{equation}
The Julia set inside $\pi^{-1}(1)$ is equal to $\{1\}\times J_1(f)=\bigcap_{n\in\mathbb{N}} f_1^{-n}(\mathcal{M}(1))$, and since $q_1$ belongs to the shift locus $\mathcal{S}_d$, this space is homeomorphic to $\mathcal{C}_\mathcal{A}$ via a homeorphism $\omega:\{1\}\times J_1(f)\to \mathcal{C}_{\mathcal{A}}$ that conjugates $f_1$ to $\sigma$. We recall these facts in Proposition \ref{prop:lecodageomega} below. In particular the points $p_j:=\omega^{-1}(j\cdots j\cdots),\ j\in\{1,\cdots,d\},$ are exactly the fixed points of $f$ contained in $\pi^{-1}(1)$. The first result we obtain is the following:

\begin{thm}\label{thm:CompConnexesJ(f)}
    For each $j\in\{1,\cdots,d\}$ the connected component $\mathcal{C}_j$ of $p_j$ inside $J(f)$ is equal to $\mathcal{C}_j=\mathcal{P}_j([0,1])$, where $\mathcal{P}_j:[0,1]\to J(f)$ is a continuous loop starting at $p_j$. The path $\mathcal{C}_j$ winds exactly $m_j$ times above $\mathbb{S}^1$, where $m_j$ is the order of $\omega(p_j)$ under $h$. Two such curves $\mathcal{C}_i$ and $\mathcal{C}_j$ are equal if and only if $\omega(p_i)$ and $\omega(p_j)$ are in the same orbit of $h$. 

    More generally, any connected component of $J(f)$ is parameterised by a loop winding at least one time above $\mathbb{S}^1$.
\end{thm}

By winding a certain number of times above $\mathbb{S}^1$ we understand the following: a continuous loop $\gamma:[0,1]\to\mathcal{M}$ winds $m$ times above $\mathbb{S}^1$ if $\gamma$ has the form
$$\gamma(t)=(e^{2i\pi tm},W(t))$$
and if the points $W(k/m),\ k\in\{0,\cdots,m-1\},$ are distinct two by two.\\ 

It is then natural to consider return maps for each $t\in[0,1]$ using Theorem \ref{thm:CompConnexesJ(f)}. For each $z\in J_p$, denote $J_z:=J_z(f)$ for simplicity. Then one can consider for each $t\in[0,1]$ a map (with $z=e^{2i\pi t}$):
$$\mathcal{P}_t:\{z\}\times J_z\to \{1\}\times J_1$$
which is defined as follows: if $p=(z,w)\in \{z\}\times J_z$ then by Theorem \ref{thm:CompConnexesJ(f)} there exists $\mathcal{P}_p:[0,1]\to J(f)$ a continuous loop starting at $p$ which parametrize (in the trigonometric sense) the connected component of $(z,w)$ in $J(f)$, so we can define 
\begin{equation}\label{eq:mapsP_t}
    \mathcal{P}_t(p):=\mathcal{P}_p\left(\min\{s\in I_t\ :\ \mathcal{P}_p(s)\in\pi^{-1}(1)\}\right),
\end{equation}
where $I_t=[0,1]$ if $t<1$ and $I_t=]0,1]$ if $t=1$. These return maps furnish a map:
\begin{equation}\label{eq:Returnmap}
    \mathcal{P}:\left\{
        \begin{array}{ll}
          J(f)  & \longrightarrow\ \mathcal{S}_{\mathcal{A}}\\
          p=(e^{2i\pi t},w) & \longmapsto\ \left(t,\omega(\mathcal{P}_t(p))\right),\ t<1
        \end{array}
    \right.
\end{equation}
where $\mathcal{S}_{\mathcal{A}}$ is the suspension defined by the quotient:
\begin{equation}\label{eq:suspensionS_A}
    \mathcal{S}_{\mathcal{A}}:=\left([0,1]\times\mathcal{C}_{\mathcal{A}}\right)\left/\displaystyle{\phantom{\int}\!\!\!\!\!\!(0,x)\sim_h(0,y)\Leftrightarrow h(x)=y}\right. .
\end{equation}
The second result we obtained is then:

\begin{thm}\label{thm:Returnmapisahomeo}
    The map \eqref{eq:Returnmap} is a homeomorphism.
\end{thm}

At last, we come to ambient isotopies. Let $\sim$ be the equivalent relation on compact subsets of $\mathbb{S}^1\times\cmplex$ defined by:
\begin{align}
    C\sim C' \Longleftrightarrow\ & \exists\mathrm{ambient\ isotopy}\ \left(h_t:\mathbb{S}^1\times\cmplex\overset{\simeq}{\longrightarrow}\mathbb{S}^1\times\cmplex\right)_{t\in[0,1]}\ \notag\\ 
                                 & \mathrm{deforming}\ C\ \mathrm{into}\ C'\ \mathrm{of\ the\ form}\ h_t(z,w)=(z,k_t(z,w)). \label{eq:equivrelaambientisotopy}
\end{align}
If $C\sim C'$ we say that $C$ and $C'$ are ambient isotopic in $\mathbb{S}^1\times\cmplex$. Denote $\mathcal{K}$ the associated isotopy class of $\bigcup_{j=1}^{d}\mathcal{C}_j$ (where $\mathcal{C}_j$ is the connected component of $p_j$ given by Theorem \ref{thm:CompConnexesJ(f)}). The third result we obtained is then:

\begin{thm}\label{thm:Kisotopicf^-1(S^1times0)}
    The class $\mathcal{K}$ is equal to the class of $f^{-1}(\mathbb{S}^1\times\{0\})$ for the relation \eqref{eq:equivrelaambientisotopy}. 
\end{thm}

\subsection{Description of connected components of fixed points}\label{sec:topologicaldescription}

Recall that $J(f)=\bigcup_{z\in\mathbb{S}^1}\{z\}\times J_z$, where $J_z=J_z(f)$. We begin by recalling that $J_1\subset \mathbb{D}_{R}$ is a Cantor set on $d$ letters in Proposition \ref{prop:lecodageomega} below. Since $w\mapsto q_1(w)$ belongs to the shift locus $\mathcal{S}_d$, we mention that the content of Proposition \ref{prop:lecodageomega} is classic and well known, so we prefer to avoid the proof. 

Let us recall some notations and properties of the Poincaré metric $\rho_R$ introduced in Section \ref{sec:admissibleparameters}. We fix a disc $D$ centered at $0$ such that $\mathcal{M}_1(z)\subset\{z\}\times D\Subset\mathcal{M}(z)$ for each $z\in\mathbb{S}^1$, and we define a metric $d_{\mathcal{M}}$ on $\mathcal{M}$ such that $d_{\mathcal{M}}$ coincide with the Poincaré metric $\rho_R$ on $\mathbb{D}_R$ on each fiber of $\pi$, see \eqref{eq:defnofd_M}. Let $(g_j)_{1\leq j\leq d}$ be the inverse branches of $f_1:\mathcal{M}_1(1)\to\mathcal{M}(1)$ such that $g_j(t_0)=x_j$ for each $1\leq j\leq d$. Each inverse branch $g_j\in\Gamma(1)$ is contracting for the metric $d_{\mathcal{M}}$ by Proposition \ref{lemma:contractiondanslesfibres}: 
\begin{equation}\label{eq:g_jcontracte}
        \forall x,y\in\{1\}\times\overline{D},\ d_{\mathcal{M}}(g_j(x),g_j(y))\leq k(f)\times d_{\mathcal{M}}(x,y)\leq\rho k(f),
\end{equation}
where $\rho=\sup\{d_{\mathcal{M}}(a,b),\ a,b\in\{1\}\times\overline{D}\}$. Recall also that $d_{\mathcal{M}}$ is equivalent to the euclidean metric on $\{1\}\times\overline{D}$, see \eqref{eq:equivalenced_Met||.||surD}. Then it is a classical fact that these inverse branches generate an Iterated Function System (IFS) such that:
$$\{1\}\times J_1=\bigcap_{n=0}^{+\infty}f_1^{-n}(\mathcal{M}(1))=\bigcap_{n=1}^{+\infty}\bigcup_{1\leq j_1,\cdots,j_n\leq d}g_{j_1}\circ\cdots\circ g_{j_n}\left(\{1\}\times\mathbb{D}_R\right).$$
This IFS is associated with a natural coding map $\omega:\{1\}\times J_1\to\mathcal{C}_{\mathcal{A}}$ which satisfies the content of the following proposition:

\begin{prop}\label{prop:lecodageomega} There exists a map $\omega:\{1\}\times J_1\to\mathcal{C}_{\mathcal{A}}$ such that:
    \begin{enumerate}
        \item If $\omega(x)=(a_1\cdots a_n\cdots)$ then $x=g_{a_1}\circ\cdots\circ g_{a_n}(f^n(x)), \forall n\geq1$.
        \item $\omega:\{1\}\times J_1\to\mathcal{C}_{\mathcal{A}}$ is a homeomorphism such that $\omega\circ f_1=\sigma\circ\omega$ on $\mathcal{C}_{\mathcal{A}}$ and:
        \begin{equation}\label{eq:inversedeomega}
        \omega^{-1}:
            \left\{
            \begin{array}{ll}
                \mathcal{C}_{\mathcal{A}} & \longrightarrow \{1\}\times J_1\\
                (a_1\cdots a_n\cdots)     & \longmapsto     \lim_{n}g_{a_1}\circ\cdots\circ g_{a_n}(t_0).
            \end{array}
            \right.
        \end{equation}
        \item $f_1:\mathcal{M}_1(1)\to\mathcal{M}(1)$ has exactly $d$ fixed points $p_1,\cdots,p_d$ and they satisfy:
        $$\omega(p_j)=(j,\cdots,j\cdots),\ j\in\{1,\cdots,d\}.$$
    \end{enumerate}
\end{prop}
By construction of $L:\mathcal{A}^*\to X$ using the monodromy of $f$, if $(a_n)_{n\geq1}\in\mathcal{C}_{\mathcal{A}}$ then it is not difficult to observe that for each $n$ we have $L(a_n\cdots a_1)=g_{a_1}\circ\cdots\circ g_{a_n}(t_0)$. Therefore, by using the third item of Proposition \ref{prop:lecodageomega} and the contraction property \eqref{eq:g_jcontracte} of the inverse branches, we deduce the following lemma:

\begin{lemme}\label{lemma:Lemma1}
    Let $x\in\{1\}\times J_1$ and let $(a_1\cdots a_n\cdots)=\omega(x)$. Then we have for each $n\geq 1$:
    $$L(a_n\cdots a_1)=g_{a_1}\circ\cdots\circ g_{a_n}(t_0)\ \mathrm{and}\ x=g_{a_1}\circ\cdots\circ g_{a_n}\circ f^n(x),$$
    thus we deduce:
    $$\lim_{n\to+\infty} L(a_n\cdots a_1)=x.$$
\end{lemme}
Recall that $S\in\mathfrak{S}_d$ is defined in Definition \ref{defn:definitiondeS} and that $h$ is defined by \eqref{eq:homeoh}. Observe that for each $j\in\{1,\cdots,d\}$ we have:
\begin{equation}\label{eq:defnofm_j}
    \#\{h^{n}(\omega(p_j)),\ n\geq0\}=\#\{S^n(j),\ n\geq0\}=:m_j.
\end{equation}

\begin{lemme}\label{lemma:nestedintersectionandconnectedcomponent}
    Let $j\in\{1,\cdots,d\}$. For each $n\geq1$, we denote:
    $$\mathcal{P}_{j,n}(t):=f^{-n}(\gamma_{d^nm_j})_{L(nj)}(t),\ \forall t\in[0,1],$$
    with $L(nj)=L(j\cdots j)$ where the integer $j$ is repeated $n$ times. For each $n\geq1$, the path $\mathcal{P}_{j,n}$ is contained in the connected component $\mathcal{M}_{n,j}$ of $\mathcal{M}_n$ which contains $L(nj)$. The sequence $(\mathcal{M}_{n,j})_{n\geq1}$ satisfies the following points:
    \begin{enumerate}
        \item $(\mathcal{M}_{n,j})_{n\geq1}$ is decreasing and $\mathcal{M}_j:=\bigcap_{n\in\mathbb{N}^*}\mathcal{M}_{n,j}$ is connected and contains $p_j$. 
        \item The set $\mathcal{M}_{j}$ is a connected component of $ J(f)$.
    \end{enumerate}
\end{lemme}
\noindent\textbf{\underline{Proof$:$}} 
Since $\mathcal{P}_{j,n}(0)=L(nj)\in\mathcal{M}_{n,j}$ and since the path $\mathcal{P}_{j,n}$ is by construction contained inside $\mathcal{M}_n$, we can conclude that it is actually contained inside the connected component $\mathcal{M}_{n,j}$. Let us now prove the two items.\\ 

\noindent\textit{1.} According to Lemma \ref{lemma:Lemma1}, $L((n+1)j)=g_j\circ L(nj)$ thus $L((n+1)j)$ and $L(nj)$ must be contained inside the same connected component of $\mathcal{M}_n(1)$ and so inside the same connected component of $\mathcal{M}_n$. We then deduce that $\mathcal{M}_{n+1,j}\subset\mathcal{M}_{n,j}$. Since $\mathcal{M}_{n+1}\Subset\mathcal{M}_n$ we immediately deduce $\mathcal{M}_{n+1,j}\Subset\mathcal{M}_{n,j}$. 

For every $n\geq1$, using Lemma \ref{lemma:Lemma1} we can observe that $p_j$ belongs to the same connected component of $L(nj)$ inside $\mathcal{M}_n(1)$, thus it belongs to the same connected component of $L(nj)$ inside $\mathcal{M}_n$, that is $\mathcal{M}_{n,j}$. This shows that $p_j\in\bigcap_{n\geq1}\mathcal{M}_{n,j}$. At last, $(\overline{\mathcal{M}}_{n,j})_{n\geq1}$ is a nested sequence of connected compact subsets (such that $\overline{\mathcal{M}}_{n+1,j}\subset\mathcal{M}_{n,j}$) as proved just below, thus $\bigcap_{n\geq1}\mathcal{M}_{n,j}=\bigcap_{n\geq1}\overline{\mathcal{M}}_{n,j}$ is connected.\\

\noindent\textit{2.} It is straightforward to see that $\mathcal{M}_j\subset\bigcap_{n\in\mathbb{N}}f^{-n}(\mathcal{M})=J(f)$. Let us then prove that $\mathcal{M}_j$ is a connected component of $J(f)$. Let $\mathcal{C}_j$ be the connected component of $J(f)$ containing the point $p_j$. Then $\mathcal{M}_j\subset\mathcal{C}_j$ according to the preceding point. For $n\geq1$, we also have $\mathcal{C}_j\subset f^{-n}(\mathcal{M})=\mathcal{M}_n$, and thus $\mathcal{C}_j\subset \mathcal{M}_{n,j}$ since $p_j\in\mathcal{M}_{n,j}$ as explained above. At last $\mathcal{C}_j=\mathcal{M}_j$.\qed\\

We can then prove the following result which implies the first part of Theorem \ref{thm:CompConnexesJ(f)}:

\begin{prop}\label{prop:Lemma2}
    Let $j\in\{1,\cdots,d\}$. For each $n\geq1$, the path $\mathcal{P}_{j,n}$ has the following form:
    $$\mathcal{P}_{j,n}(t)=\left(e^{2i\pi t m_j}, W_{j,n}(t)\right),\ \forall t\in[0,1],$$
    where $W_{j,n}:[0,1]\to\mathbb{D}_R$ is a continuous function.
    Moreover, the sequences $(W_{j,n})_{n\geq1}$ and $(\mathcal{P}_{j,n})_{n\geq1}$ converges uniformly on $[0,1]$ to a continuous function $W_j:[0,1]\to \mathbb{D}_{R}$ and to a path $\mathcal{P}_{j}:[0,1]\to J(f)$ such that:
    \begin{enumerate}
        \item $\mathcal{P}_{j}(t)=(e^{2i\pi tm_j},W_j(t))$ for any $t\in[0,1]$.
        \item The points $\mathcal{P}_j(k/m_{j}),\ k\in\{0,\cdots,m_{j}-1\},$ are distinct two by two.
        \item $\mathcal{P}_{j}(0)=\mathcal{P}_{j}(1)=p_j$.
    \end{enumerate}
    In other words, $\mathcal{P}_{j}$ is a loop based at $p_j$ winding $m_j$ times above $\mathbb{S}^1$.
    \begin{enumerate}
        \item[4.] $\mathcal{C}_j:=\mathcal{P}_j([0,1])$ is the connected component of $J(f)$ containing $p_j$.
    \end{enumerate}
\end{prop}
\noindent\textbf{\underline{Proof$:$}}\ \\
\noindent{\textit{1.}} Let us fix $t\in[0,1]$. We denote $w_k:=W_{j,k}(t)$ and $z_k:=e^{2 i\pi t d^k m_j}$. For any $z\in\mathbb{S}^1$, let $\mathcal{M}_{1,j}(z):=\mathcal{M}_{1,j}\cap\pi^{-1}(z)$, where $\mathcal{M}_{1,j}$ is the connected component of $\mathcal{M}_1$ which contains $L(j)$. Since $L(nj)=g_{j}^{\circ n}(1,0)$, $\mathcal{M}_{1,j}$ is also the connected component of $L(nj)$ in $\mathcal{M}_1$. 

Let us prove the following property:
\begin{equation}\label{eq:etrememevoisinage}
    \forall n\geq1,\forall l\geq n,\forall k\in\{0,\cdots,n\},\ (z_{k},Q_{z_0}^{k}(w_l))\in\mathcal{M}_{1,j}(z_{k}).
\end{equation}
Let $n\geq1$, $l\geq n$ and $k\in\{0,\cdots,n\}$ be fixed. Then observe that $(z_{k},Q_{z_0}^{k}(w_l))=f^{k}(z_0,w_l)=f^{k}\circ\mathcal{P}_{j,l}(t)=f^{-(l-k)}(\gamma_{d^{l}m_j})_{L((l-k)j)}(t)$, by using the definition of $\mathcal{P}_{j,l}(t)$. From the equality $(z_{k},Q_{z_0}^{k}(w_l))=f^{-(l-k)}(\gamma_{d^{l}m_j})_{L((l-k)j)}(t)$ we deduce that $(z_{k},Q^{k}_{z_0}(w_l))$ must belongs to the connected component of $L((l-k)j)$ in $\mathcal{M}_1$, that is $\mathcal{M}_{1,j}$. So we have $(z_{k},Q_{z_0}^{k}(w_l))\in\mathcal{M}_{1,j}(z_{k})$ and we have proved \eqref{eq:etrememevoisinage}. 

Since $\mathcal{M}_{1,j}(z_{k})$ is by construction a connected component of $\{z_{k}\}\times q_{z_{k}}^{-1}(\mathbb{D}_R)$, we deduce that there exists $g_{k}\in\Gamma(z_{k})$ of the form $g_{k}(z_{k},w)=(z_{k},\widehat{g}_{z_{k}}(w))$ such that $Q^{k}_{z_0}(w_l)\in\widehat{g}_{z_{k}}(\mathbb{D}_R)$ (see Section \ref{sec:admissibleparameters} and Proposition \ref{lemma:contractiondanslesfibres}). Since $q_{z_k}\circ\widehat{g}_{z_k}=\mathrm{Id}$, we deduce $Q^{k}_{z_0}(w_l)=\widehat{g}_{z_{k}}(Q^{k+1}_{z_0}(w_l))$. In particular, we deduce by induction that $w_l=\widehat{g}_{z_0}(Q^1_{z_0}(w_l))=\widehat{g}_{z_0}\circ\widehat{g}_{z_1}(Q^2_{z_0}(w_l))=\cdots=\widehat{g}_{z_0}\circ\cdots\circ\widehat{g}_{z_{n-1}}\left(Q^n_{z_0}(w_l)\right)$. So, for each $n\geq1$ and each $l\geq n$:
\begin{equation}\label{eq:w_l=gcircg(Q(w_l))}
    w_l=\widehat{g}_{z_0}\circ\cdots\circ\widehat{g}_{z_{n-1}}\left(Q^n_{z_0}(w_l)\right).
\end{equation}
Applying \eqref{eq:w_l=gcircg(Q(w_l))} to $l=n$ and $l=n+m$ with $m\geq0$, we get:
$$w_{n}=\widehat{g}_{z_0}\circ\cdots\circ\widehat{g}_{z_{n-1}}\left(Q^n_{z_0}(w_n)\right)\ \mathrm{and}\ w_{n+m}=\widehat{g}_{z_0}\circ\cdots\circ\widehat{g}_{z_{n-1}}\left(Q^n_{z_0}(w_{n+m})\right).$$
By Proposition \ref{lemma:contractiondanslesfibres} we deduce $\rho_R(w_{n},w_{n+m})\leq k(f)^n\rho_R(Q^{n}_{z_0}(w_{n}),Q^n_{z_0}(w_{n+m}))\leq \rho k(f)^{n}$, where $\rho=\sup\{\rho_R(a,b),\ a,b\in\overline{D}\}$. Observe that all the preceding arguments are valid for any $t\in[0,1]$ and any $n\geq1$, $m\geq 0$. So, we finally deduce :
$$\forall n\geq1,\forall m\geq0,\ \sup_{t\in[0,1]}\rho_R\left(W_{j,n}(t),W_{j,n+m}(t)\right)\leq \rho k(f)^n.$$
Since $k(f)<1$ we finally deduce that $(W_{j,n})_{n\geq0}$ is a Cauchy sequence which converge uniformly on $[0,1]$ to a continuous function $W_j:[0,1]\to\overline{D}$. 

Let $t\in[0,1]$, $n>m\geq0$ and denote $z_k:=e^{2i\pi td^km_j}$. From \eqref{eq:w_l=gcircg(Q(w_l))} we observe that $Q^m_{z_0}(W_{j,n}(t))\in\widehat{g}_{z_m}(\mathbb{D}_R)\subset\overline{D}$, for some $g_{m}\in\Gamma(z_m)$. Taking the limit as $n$ goes to infinity, we get $Q^m_{z_0}(W_j(t))\in\overline{D}$. So, we have $\mathcal{P}_j(t)=(z_0,W_j(t))\in\bigcap_{m\geq0}f^{-m}(\mathcal{M})=J(f)$. Since this true for any $t\in[0,1]$, we can conclude that $\mathcal{P}_j$ takes its values in $J(f)$.\\

\noindent\textit{2.} Let us fix $k\in\{0,\cdots,m_j-1\}$ and let $P_k:=\mathcal{P}_{j}(k/m_j)$ and $P_{n,k}:=\mathcal{P}_{j,n}(k/m_j)$, we recall that $P_{n,k}\to P_k$. Observe that $t\in[0,1]\mapsto f^{-n}(\gamma_{d^nm_j})_{L(nj)}(tk/m_j)$ and $f^{-n}(\gamma_{d^nk})_{L(nj)}$ are two lifts by $f^n$ of $\gamma_{d^nk}$. By uniqueness of lifts by $f^n$ these two lifts are equal, so we deduce that $P_{n,k}=\mathcal{P}_{j,n}(k/m_j)=f^{-n}(\gamma_{d^nm_j})_{L(nj)}(k/m_j)=f^{-n}(\gamma_{d^nk})[L(nj)]$. By definition of the monodromy action \eqref{eq:morphismedegroupesX}, we have thus $P_{n,k}=L(nj)^{\gamma_{d^nk}}$. So, by applying Corollary \ref{cor:lafomulederecurrencedeNekrashevych} we deduce:
\begin{align*}
    P_{n,k} & = L(j\cdots j)^{\gamma_{d^nk}} = L\left(S^{d^{n-1}k}(j)\cdots S^{k}(j)\right).
\end{align*}
In particular, using Lemma \ref{lemma:Lemma1} it gives :
\begin{align*}
    P_{k} = \lim_{n\to+\infty}P_{n,k} & = \lim_{n\to+\infty}L\left(S^{d^{n-1}k}(j)\cdots S^{k}(j)\right)\\
          & = \lim_{n\to+\infty} g_{S^k(j)}\circ\cdots\circ\cdots g_{S^{d^{n-1}k}}(t_0)\ \mathrm{by\ Lemma\ \ref{lemma:Lemma1}},\\
          & = \omega^{-1}\left(S^k(j)\cdots S^{d^{n-1}k}(j)\cdots\right)\ \mathrm{by\ \eqref{eq:inversedeomega}}.
\end{align*}
Since the elements of $\{S^{k}(j),\ k=0,\cdots,m_j-1\}$ are distinct two by two by definition of $m_j$, we deduce that the points $(P_k)_{k=0,\cdots,m_j-1}$ do not have the same image by $\omega$, thus they are distinct two by two. The conclusion then follows.\\

\noindent\textit{3.} At last we have $\mathcal{P}_{j,n}(1)=L(nj)^{\gamma_{d^nm_j}}=L(S^{d^{n-1}m_j}(j)\cdots S^{d^{0}m_j}(j))=L(nj)$
by Corollary \ref{cor:lafomulederecurrencedeNekrashevych} and by $S^{m_j}(j)=j$. Thus $\mathcal{P}_{n,j}(1)=L(nj)\to p_j$ and $\mathcal{P}_j(1)=p_j$. 
Similarly, $\mathcal{P}_{j,n}(0)=L(nj)\to p_j$ and thus $\mathcal{P}_j(0)=p_j$ as well. So we have proved that $\mathcal{P}_j$ is a loop based at $p_j$ inside $J(f)$.\\

\noindent\textit{4.} Let $\mathcal{M}_j$ be the connected component of $J(f)$ containing $p_j$. To finish the proof, it remains to prove that $\mathcal{C}_j=\mathcal{M}_j$, where $\mathcal{C}_j=\mathcal{P}_j([0,1])$. According to Lemma \ref{lemma:nestedintersectionandconnectedcomponent}, $\mathcal{M}_j=\bigcap_{n\geq1}\mathcal{M}_{n,j}$ where $\mathcal{M}_{n,j}$ is the connected component of $L(nj)$ in $\mathcal{M}_n$. Since the curve $\mathcal{P}_j$ takes its values in $J(f)$, and since $\mathcal{P}_j(0)=p_j$, we have that $\mathcal{C}_j\subset\mathcal{M}_j$. Let us prove the reverse implication. Let us fix $p:=(z_0,w_0)\in\mathcal{M}_j$, and let $t\in[0,1]$ such that $z_0=e^{2i\pi t m_j}$. Observe that $f^{n}(p)\in\mathcal{M}_{1,j}(z_n)$ with $z_n=e^{2i\pi td^nm_j}$ for each $n$. So, similarly has we have proved \eqref{eq:w_l=gcircg(Q(w_l))}, we can prove that there exists $g_n\in\Gamma(z_n)$ of the form $g_n(z_n,w)=(z_n,\widehat{g}_{z_n}(w))$ such that $w_0=\widehat{g}_{z_0}\circ\cdots\circ\widehat{g}_{z_{n-1}}(Q^n_{z_0}(w_0))$ for each $n$. Similarly, with the same sections $g_n\in\Gamma(z_n)$, we will have $W_{j,n}(t)=\widehat{g}_{z_0}\circ\cdots\circ\widehat{g}_{z_{n-1}}(Q^n_{z_0}(W_{j,n}(t)))=\widehat{g}_{z_0}\circ\cdots\circ\widehat{g}_{z_{n-1}}(0)$ for each $n$. It implies by Proposition \ref{lemma:contractiondanslesfibres} that ${d}_{\mathcal{M}}(p,\mathcal{P}_{j,n}(t))\leq \rho k(f)^n$ for each $n$, and thus $p=\mathcal{P}_j(t)$ proving that $\mathcal{M}_j\subset\mathcal{C}_j$.\qed

\subsection{Description of all connected components}\label{sec:topologicaldescriptionIII}

To finish the proof of Theorem \ref{thm:CompConnexesJ(f)}, we would like to generalize the preceding description of connected components of periodic points for any point $p\in J(f)$. To do so, we introduce:
\begin{equation}\label{eq:ordredeS}
    m_S:=\min\{n\geq 1\ :\ S^n=\mathrm{Id}_{\mathfrak{S}_d}\}.
\end{equation}

\begin{prop}\label{prop:connectedcomponentpointC_A}
    Let $p\in\{1\}\times J_1$ and denote $\omega(p)=(a_1\cdots a_n\cdots)$. We then denote for any $n\geq1$ and for any $t\in[0,1]$:
    $$\mathcal{P}_{p,n}(t):=f^{-n}(\gamma_{d^{n}m_S})_{L\left(a_n\cdots a_1\right)}(t),$$
    where $m_S$ is defined by \eqref{eq:ordredeS}. We can write $\mathcal{P}_{p,n}(t)=(e^{2i\pi tm_S},W_{p,n}(t))$. Then the sequence of continuous functions $W_{p,n}:[0,1]\to\overline{D}$ converges uniformly to a continuous function $W_p$. We then denote $\mathcal{P}_p(t):=(e^{2i \pi tm_S},W_p(t)):[0,1]\to\mathcal{M}$ which takes its values in $J(f)$ and which satisfies:
    $$\mathcal{P}_{p}(0)=\mathcal{P}_{p}(1)=p.$$
    In other words, $\mathcal{P}_{p}$ is a loop based at $p$ winding at least one time above $\mathbb{S}^1$. Moreover, $\mathcal{C}_p:=\mathcal{P}_p([0,1])$ is the connected component of $J(f)$ containing $p$, and we have:
    \begin{equation}\label{eq:P(k/m_S)}
    \mathcal{P}_{p}\left(\frac{k}{m_S}\right)=\lim_{n\to+\infty}L\left(S^{d^{n-1}k}(a_n)\cdots S^{d^0k}(a_1)\right),\ 0\leq k\leq m_S.
\end{equation}
\end{prop}
\noindent\textbf{\underline{Proof$:$}}
The proof is similar to the one of Proposition \ref{prop:Lemma2}. The reader should only notice that for $k\in\{0,\cdots,m_S\}$, we have 
$\mathcal{P}_{p}\left(\frac{k}{m_S}\right)=\lim_{n}L(a_n\cdots a_1)^{\gamma_{d^nk}}$ which is equal to $\lim_{n}L\left(S^{d^{n-1}k}(a_n)\cdots S^{d^0k}(a_1)\right)$ by Corollary \ref{cor:lafomulederecurrencedeNekrashevych}, thus we have \eqref{eq:P(k/m_S)}. In particular, for $k=0$ and for $k=m_S$ we get $\mathcal{P}_p(0)=\mathcal{P}_p(1)=p$.\qed\\ 

We can extend this result for every point in the Julia set as follows:

\begin{thm}\label{thm:theoremedescourbes}
    Let $p\in J(f)$. Then there exists $\widetilde{p}\in\{1\}\times J_1$ such that the connected component $\mathcal{C}_p$ of $p$ inside $J(f)$ is equal to $\mathcal{C}_p=\mathcal{P}_{\widetilde{p}}([0,1])$ (where $\mathcal{P}_{\widetilde{p}}$ is given by Prop. \ref{prop:connectedcomponentpointC_A}). One can reparametrize the loop $\mathcal{P}_{\widetilde{p}}$ into a loop $\mathcal{P}_p$ based at $p$, and winding above $\mathbb{S}^1$ in the trigonometric sense.     
\end{thm}
\noindent\textbf{\underline{Proof$:$}} Let us denote $\mathcal{M}_n'$ the connected component of $p$ in $\mathcal{M}_n$, for each $n\geq1$. As in Lemma \ref{lemma:nestedintersectionandconnectedcomponent}, one can show that $\mathcal{C}_p=\mathcal{M}'$, where $\mathcal{M}'=\bigcap_{n\geq1}\overline{\mathcal{M}}_n'$, by showing that $(\overline{\mathcal{M}}_n')_n$ is a nested sequence of connected compact sets. 
Since $\pi:\mathcal{M}\to\mathbb{S}^1$ is surjective, and since $f^n:\mathcal{M}_n\to\mathcal{M}$ is a (surjective) covering map, we have $\pi(\mathcal{M}_n')=\mathbb{S}^1$ for each $n\geq1$. It is a well known topological property that, since $(\overline{\mathcal{M}}_n')_n$ is a nested sequence of compact subsets of the solid torus $\mathcal{M}$, we have $\pi\left(\bigcap_{n\geq1}\overline{\mathcal{M}}'_n\right)=\bigcap_{n\geq1}\pi\left(\overline{\mathcal{M}}_n'\right)=\mathbb{S}^1$.
So there exists $\widetilde{p}\in\mathcal{C}_p\cap\pi^{-1}(1)$ and thus $\mathcal{C}_p=\mathcal{C}_{\widetilde{p}}=\mathcal{P}_{\widetilde{p}}([0,1])$ by Proposition \ref{prop:connectedcomponentpointC_A}.\qed

\subsection{Technical results : equicontinuity and convergences}\label{sec:uniform_equicontinuity} 

In the following next subsection  our aim is to give a topological description of the Julia set $J(f)$ using the fact that the connected components of $J(f)$ are continuous loops winding above $\mathbb{S}^1$. But to fulfill a rigorous proof of this topological description we need to have additional results concerning the convergence of the curves $(\mathcal{P}_{p,n})_{n}$ to $\mathcal{P}_p$ for each $p\in\{1\}\times J_1$, see Proposition \ref{prop:connectedcomponentpointC_A}. To do so, in what follows we consider the following family 
\begin{equation}\label{eq:familleH}
    \mathcal{H}:=\left\{\mathcal{P}_{p},\ p\in\{1\}\times J_1\right\}
\end{equation}
of elements of the space $C^0([0,1],\mathcal{M})$. 

Let us recall that there exists a covering of $J_p=\mathbb{S}^1$ given by \eqref{eq:recouvrementduJuliaJ_p}:
\begin{equation}\label{eq:recouvrementcompactS^1}
    \mathbb{S}^1=\bigcup_{j=1}^{N}I_j,
\end{equation}
and according to Proposition \ref{lemma:contractiondanslesfibres}, there exists a family $\underset{\ \ \ \ \ 1\leq k\leq d}{(\mathcal{G}_{jk})_{1\leq j\leq N}}$ of sections of the covering $f:\mathcal{M}_1\to\mathcal{M}$ such that each section $\mathcal{G}_{jk}$ is defined on $I_j\times\mathbb{D}_R$, and the Lipschitz constant for $d_{\mathcal{M}}$ of $\mathcal{G}_{jk}$ on $I_j\times\overline{D}$ (recall $\overline{D}\subset\mathbb{D}_R$) is equal to $l_{jk}\leq k(f)<1$. 

\begin{prop}\label{eq:lemmeequi1}
    The family $\mathcal{H}$ of $C^0([0,1],\mathcal{M})$ defined by \eqref{eq:familleH} is (uniformly) equicontinuous. 
\end{prop}
\noindent\textbf{\underline{Proof$:$}} Let $\delta>0$ be a Lebesgue number of the covering \eqref{eq:recouvrementcompactS^1} i.e. for any $A\subset\mathbb{S}^1$ of diameter $\leq\delta$, there exists $j$ such that $A\subset I_j$. Let us fix $\varepsilon>0$ very small. There exists $n_{\varepsilon}\geq1$ large enough such that (with $\rho=\sup\{d_{\mathcal{M}}(a,b),\ a,b\in\mathcal{M}_1\}$):
\begin{equation}\label{eq:deltacontrolepsilon3}
    \rho k(f)^{n_{\varepsilon}+1}<\varepsilon.
\end{equation}
Let $\delta_{\varepsilon}>0$ be sufficiently small so that:
\begin{equation}\label{eq:deltacontrolepsilon2}
    \underset{0\leq k\leq n_{\varepsilon}}{\sup_{|x|\leq\delta_{\varepsilon}}}\left|e^{2i\pi xd^km_S}-1\right|\leq\delta.
\end{equation}
Let us fix $t_1,t_2\in[0,1]$ such that $|t_1-t_2|\leq\delta_{\varepsilon}$, let us fix $n\geq n_{\varepsilon}$, and let us fix $p\in\{1\}\times J_1$. Denote $\omega(p)=(a_1\cdots a_n\cdots)$, and denote $z_k(t_l)=e^{2i\pi t_ld^km_S}$, $l=1,2$. For each $k\in\{0,\cdots,n_{\varepsilon}\}$, we have $|z_k(t_1)-z_k(t_2)|=|e^{2i\pi(t_1-t_2)d^km_S}-1|\leq\delta$ by \eqref{eq:deltacontrolepsilon2}, thus since $\delta$ is a Lebesgue number of the covering \eqref{eq:recouvrementcompactS^1}, there exists $j(k)$ such that $z_k(t_1)\in I_{j(k)}$ and $z_k(t_2)\in I_{j(k)}$. Recall that by definition $\mathcal{P}_{p,n}(t)=f^{-n}(\gamma_{d^nm_S})_{L(a_n\cdots a_1)}(t)$, so similarly as we did in the proof of Proposition \ref{prop:Lemma2}, by using Proposition \ref{lemma:contractiondanslesfibres} we deduce that there exist inverse branches $\mathcal{G}_{j(0)k(0)},\cdots,\mathcal{G}_{j(n_{\varepsilon})k(n_{\varepsilon})}$ of $f$ such that:
$$\mathcal{P}_{p,n}(t_l)=\mathcal{G}_{j(0)k(0)}\circ\cdots\circ\mathcal{G}_{j(n_{\varepsilon})k(n_{\varepsilon})}(P_l),$$
where $P_l:=f^{-(n-n_{\varepsilon}-1)}(\gamma_{d^nm_S})_{L(a_n\cdots a_{n_{\varepsilon}+2})}(t_l)\in\overline{D}$. In consequence, we have using \eqref{eq:deltacontrolepsilon3}:
\begin{align*}
    d_{\mathcal{M}}\left(\mathcal{P}_{p,n}(t_1),\mathcal{P}_{p,n}(t_2)\right) & \leq k(f)^{n_{\varepsilon}+1}d_{\mathcal{M}}(P_1,P_2) \leq \varepsilon,\ \forall n\geq n_{\varepsilon}.
\end{align*} 
According to Proposition \ref{prop:connectedcomponentpointC_A}, $(\mathcal{P}_{p,n})_{n}$ converges uniformly to $\mathcal{P}_{p}$. So, taking the limit as $n$ goes to infinity in the last inequality, we obtain that $d_{\mathcal{M}}\left(\mathcal{P}_{p}(t_1),\mathcal{P}_{p}(t_2)\right)\leq \varepsilon$ for all $\mathcal{P}_p\in\mathcal{H}$. Since this is true for any $t_1,t_2\in[0,1]$ such that $|t_1-t_2|\leq\delta_{\varepsilon}$, and for any small $\varepsilon>0$, we deduce that $\mathcal{H}$ is uniformly equicontinuous.\qed\\

From Proposition \ref{eq:lemmeequi1} we deduce immediately the following result:

\begin{cor}\label{cor:corollairelemmeequicontinuite}
    Let $(t_n)_n$ be a sequence of $[0,1]$ which converges to $t$. Then we have:
    $$\lim_{n\to+\infty}\sup_{p\in\{1\}\times J_1}d_{\mathcal{M}}\left(\mathcal{P}_{p}(t),\mathcal{P}_{p}(t_n)\right) = 0.$$
\end{cor}

Let us now prove the following:

\begin{prop}\label{prop:lemmeequi2}
    Let $(p_n)$ be a sequence of $\{1\}\times J_1$ which converges to $p$. Then we have:
    $$\lim_{n\to+\infty}\sup_{t\in[0,1]}d_{\mathcal{M}}\left(\mathcal{P}_{p_n}(t),\mathcal{P}_{p}(t)\right)=0.$$
\end{prop}
\noindent\textbf{\underline{Proof$:$}} Let us write for every $n\geq0$ :
$$\omega(p_n)=(a_1^n\cdots a_m^n\cdots)\ \ \mathrm{and}\ \ \omega(p)=(a_1\cdots a_m\cdots).$$
Let us fix $\varepsilon>0$ and let $l_{\varepsilon}$ such that (with $\rho=\sup\{d_{\mathcal{M}}(a,b),\ a,b\in\mathcal{M}_1\}$):
\begin{equation}\label{eq:diam(D)k(f)leps<eps}
    \rho k(f)^{l_{\varepsilon}+1}<\varepsilon.
\end{equation}
Since $p_n\to p$ there exists $\varphi:\mathbb{N}\to\mathbb{N}$ a subsequence, and there exists $n_{\varepsilon}\geq1$, such that $\varphi(n)\geq n$ for all $n\geq1$, and such that for all $n\geq n_{\varepsilon}$: 
\begin{equation}\label{eq:extractricetrick2}
    a_{1}^{\varphi(n)}=a_{1}^{n}=a_1,\cdots,a_{l_{\varepsilon}}^{\varphi(n)}=a_{l_{\varepsilon}}^n=a_{l_{\varepsilon}}.
\end{equation}

Let us fix now $n\geq n_{\varepsilon}$ and $m\geq l_{\varepsilon}$. Let us also fix $t\in[0,1]$. As in the proof of Proposition \ref{prop:Lemma2}, if we denote $z_k(t)=e^{2i\pi td^km_S}$, there exists $g_{0}\in\Gamma(z_0(t)),\cdots,g_{m-1}\in\Gamma(z_{m-1}(t))$ such that:
$$W_{p_n,m}(t)=\widehat{g}_{z_0}\circ\cdots\circ\widehat{g}_{z_{l_{\varepsilon}}}\circ\widehat{g}_{z_{l_{\varepsilon}+1}}\circ\cdots\circ\widehat{g}_{z_{m-1}}(0).$$
Denoting $w_1:=\widehat{g}_{z_{l_{\varepsilon}+1}}\circ\cdots\circ\widehat{g}_{z_{m-1}}(0)\in\overline{D}$ we have:
$$W_{p_n,m}(t)=\widehat{g}_{z_0}\circ\cdots\circ\widehat{g}_{z_{l_{\varepsilon}}}(w_1).$$
Since by \eqref{eq:extractricetrick2} the $l_{\varepsilon}-$first elements of the sequences $\omega(p_n)$ and $\omega(p_{\varphi(n)})$ are equal, we can express $W_{p_{\varphi(n)},m}(t)$ using the same sections $\{g_{k}\in\Gamma(z_k(t))\}_{0\leq k\leq l_{\varepsilon}}$ than the ones used to express $W_{p_n,m}(t)$. So we have:
$$W_{p_{\varphi(n)},m}(t)=\widehat{g}_{z_0}\circ\cdots\circ\widehat{g}_{z_{l_{\varepsilon}}}(w_2),$$
where $w_2\in\overline{D}$ is defined similarly as $w_1$ but using sections $g\in\bigcup_{z\in\mathbb{S}^1}\Gamma(z)$ which are not necessarily the same than the one used to define $w_1$, since the sequences $(a^n_k)_{l_{\varepsilon}<k<m}$ and $(a_{k}^{\varphi(n)})_{l_{\varepsilon}<k<m}$ may be different. 
From these expressions we deduce that for every $m\geq l_{\varepsilon}$:
\begin{align*}
    d_{\mathcal{M}}\left(\mathcal{P}_{p_n,m}(t),\mathcal{P}_{p_{\varphi(n)},m}(t)\right) & = \rho_R\left(W_{p_n,m},W_{p_{\varphi(n)},m}\right)\\
                                               & = \rho_R\left(\widehat{g}_{z_0}\circ\cdots\circ\widehat{g}_{z_{l_{\varepsilon}}}(w_1),\widehat{g}_{z_0}\circ\cdots\circ\widehat{g}_{z_{l_{\varepsilon}}}(w_2)\right)\\
                                               & \leq k(f)^{l_{\varepsilon}+1}\rho_R(w_1,w_2) \leq \varepsilon,
\end{align*}
where the last inequality comes from \eqref{eq:diam(D)k(f)leps<eps}. 

Similarly, by using \eqref{eq:extractricetrick2}, we can express $W_{p,m}(t)$ as follows:
$$W_{p,m}(t)=\widehat{g}_{z_0}\circ\cdots\circ\widehat{g}_{z_{l_{\varepsilon}}}(w_3),$$
where $w_3\in\overline{D}$. So we also have:
\begin{align*}
    d_{\mathcal{M}}\left(\mathcal{P}_{p_{\varphi(n)},m}(t),\mathcal{P}_{p,m}(t)\right) & = \rho_R\left(W_{p_{\varphi(n)},m}(t),W_{p,m}(t)\right) \leq k(f)^{l_{\varepsilon}+1}\rho_R(w_2,w_3) \leq \varepsilon.
\end{align*}
Let $m\geq l_{\varepsilon}$, we proved
$d_{\mathcal{M}}\left(\mathcal{P}_{p_n,m}(t),\mathcal{P}_{p_{\varphi(n)},m}(t)\right)\leq\varepsilon\ \mathrm{and}\ d_{\mathcal{M}}\left(\mathcal{P}_{p_{\varphi(n)},m}(t),\mathcal{P}_{p,m}(t)\right)\leq\varepsilon$, so we have:
{\small \begin{align*}
    d_{\mathcal{M}}\left(\mathcal{P}_{p_n}(t),\mathcal{P}_{p}(t)\right) & \leq d_{\mathcal{M}}\left(\mathcal{P}_{p_n}(t),\mathcal{P}_{p_n,m}(t)\right) + d_{\mathcal{M}}\left(\mathcal{P}_{p_n,m}(t),\mathcal{P}_{p_{\varphi(n)},m}(t)\right)
    + d_{\mathcal{M}}\left(\mathcal{P}_{p_{\varphi(n)},m}(t),\mathcal{P}_{p}(t)\right)\\
    & \leq d_{\mathcal{M}}\left(\mathcal{P}_{p_n}(t),\mathcal{P}_{p_n,m}(t)\right) + \varepsilon + d_{\mathcal{M}}\left(\mathcal{P}_{p_{\varphi(n)},m}(t),\mathcal{P}_{p}(t)\right)\\
    & \leq d_{\mathcal{M}}\left(\mathcal{P}_{p_n}(t),\mathcal{P}_{p_n,m}(t)\right) + \varepsilon + d_{\mathcal{M}}\left(\mathcal{P}_{p_{\varphi(n)},m}(t),\mathcal{P}_{p,m}(t)\right) + d_{\mathcal{M}}\left(\mathcal{P}_{p,m}(t),\mathcal{P}_{p}(t)\right)\\
    & \leq d_{\mathcal{M}}\left(\mathcal{P}_{p_n}(t),\mathcal{P}_{p_n,m}(t)\right) + \varepsilon + \varepsilon + d_{\mathcal{M}}\left(\mathcal{P}_{p,m}(t),\mathcal{P}_{p}(t)\right).
\end{align*}}
Since by Proposition \ref{prop:connectedcomponentpointC_A}, $\lim_md_{\mathcal{M}}\left(\mathcal{P}_{p_n}(t),\mathcal{P}_{p_n,m}(t)\right)=0$ and $\lim_md_{\mathcal{M}}\left(\mathcal{P}_{p,m}(t),\mathcal{P}_{p}(t)\right)=0$ (and these limits are uniform in $t$), there exists $N(n,\varepsilon)\geq l_{\varepsilon}$ (depending on $n\geq n_{\varepsilon}$ that has been fixed earlier) such that for all $m\geq N(n,\varepsilon)$ we have $d_{\mathcal{M}}\left(\mathcal{P}_{p_n}(t),\mathcal{P}_{p_n,m}(t)\right)\leq\varepsilon\ \mathrm{and}\ d_{\mathcal{M}}\left(\mathcal{P}_{p,m}(t),\mathcal{P}_{p}(t)\right)\leq\varepsilon$. Thus we deduce $d_{\mathcal{M}}\left(\mathcal{P}_{p_n}(t),\mathcal{P}_{p}(t)\right)\leq4\varepsilon$. Since all the computations made above are valid for any $t\in[0,1]$ and for any $n\geq n_{\varepsilon}$, we can conclude that:
$$\sup_{t\in[0,1],\ n\geq n_{\varepsilon}}d_{\mathcal{M}}\left(\mathcal{P}_{p_n}(t),\mathcal{P}_{p}(t)\right)\leq 4\varepsilon.$$
Since that for all $\varepsilon>0$ there exists $n_{\varepsilon}$ satisfying this inequality, the proof is complete.\qed

\subsection{The Julia set is homeomorphic to a suspension}\label{sec:topologicaldescriptionIV}

In the previous two sections, we described the connected components of $J(f)$. We would like now to use this description to prove Theorem \ref{thm:Returnmapisahomeo} i.e. to
construct a homeomorphism between $J(f)$ and the suspension $\mathcal{S}_{\mathcal{A}}$ defined by \eqref{eq:suspensionS_A}. We mention that a result of this type has already been obtained in degree $d=2$ by Astorg-Bianchi \cite[Thm. 5.14]{AstBian23}.

To prove Theorem \ref{thm:Returnmapisahomeo}, we are going to construct to maps $\mathcal{P}$ and $\mathcal{L}$, see Equation \eqref{eq:applicationPandL}. We will show that $\mathcal{P}$ is continuous and bijective with $\mathcal{P}^{-1}=\mathcal{L}$, thus $\mathcal{P}$ will be a homemorphism between $J(f)$ and $\mathcal{S}_{\mathcal{A}}$. The map $\mathcal{P}$ is the one already defined in \eqref{eq:Returnmap}, thus we will obtain a proof of Theorem \ref{thm:Returnmapisahomeo}.
The mappings $\mathcal{P}$ and $\mathcal{L}$ are defined with two maps $t\in[0,1]\mapsto\mathcal{P}_t$ and $t\in[0,1]\mapsto\mathcal{L}_t$ that we should introduce first. Note that $\mathcal{P}_t$ has already be defined in \eqref{eq:mapsP_t}. From a geometrical point of view, this two maps are the first return map (forward and backward) of the monodromy above $\mathbb{S}^1$ of the connected components $\mathcal{C}_p$ of the points $p\in J(f)$. 

We denote for each $t\in[0,1]$, $I_t:=[0,1]$ if $t<1$ and $I_t:=]0,1]$ if $t=1$. 

\begin{defn}\label{defn:returnmaps} 
    The forward/backward first return maps are defined by:
    \begin{enumerate}
        \item
        Let $t\in[0,1]$, let $p\in J(f)\cap\pi^{-1}(e^{2i\pi t})$ and let $\mathcal{P}_p$ be the loop given by $\mathrm{Theorem\ \ref{thm:theoremedescourbes}}$.
        We then define the forward return map as:
        $$\mathcal{P}_t(p):=\mathcal{P}_p\left(\min\{s\in I_t:\mathcal{P}_p(s)\in\pi^{-1}(1)\}\right).$$
        \item Let $t\in[0,1]$, $p\in J(f)\cap\pi^{-1}(1)$ and let $\mathcal{P}_p$ be the loop given by $\mathrm{Theorem\ \ref{thm:theoremedescourbes}}$, we denote $\mathcal{P}_p^{-1}(s):=\mathcal{P}_p(1-s)$ the inverse loop.
        We then define the backward first return map as:
        $$\mathcal{L}_t(p):=\mathcal{P}_p^{-1}\left(\min\{s\in I_t:\mathcal{P}_p^{-1}(s)\in\pi^{-1}(e^{2i\pi t})\}\right).$$
    \end{enumerate}
\end{defn}

\begin{lemme}\label{lemma:htour}
    Let $p_1$ and $p_2$ be two points of $\{1\}\times J_1$. Let $\mathcal{P}_{p_1}$ and $\mathcal{P}_{p_2}$ be two loops based at $p_1$ and $p_2$ which are given by $\mathrm{Proposition\ \ref{prop:connectedcomponentpointC_A}}$. Denote $\mathcal{C}_{p_i}:=\mathcal{P}_{i}([0,1])$ for each $i\in\{1,2\}$. Then the following points hold:
    \begin{enumerate}
        \item $\mathcal{P}_t(p_1)=\mathcal{P}_{p_1}(t/m_S)$ for $t\in\{0,1\}$, and $\mathcal{L}_t(p_1)=\mathcal{P}_{p_1}\left(1-\frac{1-t}{m_S}\right)$ for $t\in[0,1]$.
        \item $\omega(p_2)=h(\omega(p_1))$ if and only if $\mathcal{P}_1(p_1)=p_2$. 
        \item $\mathcal{P}_1(p_1)=p_2$ if and only if $\mathcal{L}_1(p_2)=p_1$.
        \item We have $\mathcal{P}_1^{\circ n}=\omega^{-1}\circ h^n\circ \omega$ and $\mathcal{L}_1^{\circ n}=\omega^{-1}\circ h^{-n}\circ \omega$ for each $n\geq0$.
        \item We have $\mathcal{C}_{p_1}=\mathcal{C}_{p_2}$ if and only if $\omega(p_2)\in\mathrm{orb}_h(\omega(p_1))$, where $\mathrm{orb}_h(\omega(p_1))$ is the orbit of $\omega(p_2)$ under $h$.  
    \end{enumerate}
\end{lemme}
\noindent\textbf{\underline{Proof$:$}}\ \\ 
\noindent\textit{1.} $\mathcal{P}_0(p_1)=p_1$ by definition, and $\mathcal{P}_{p_1}(0)=p_1$ by construction of $\mathcal{P}_{p_1}$, thus $\mathcal{P}_0(p_1)=\mathcal{P}_{p_1}(0/m_S)$. Since the path $\mathcal{P}_{p_1}(t)$ as the form $\mathcal{P}_{p_1}(t)=(e^{2i\pi t m_S},W_{p_1}(t))$ we deduce:
$$\mathcal{P}_1(p_1)=\mathcal{P}_{p_1}\left(\min\{s\in]0,1]\ :\ \mathcal{P}_{p_1}(s)\in\pi^{-1}(1)\}\right)=\mathcal{P}_{p_1}\left(\frac{1}{m_S}\right).$$
Similarly, we observe that $\mathcal{L}_t(p_1)=\mathcal{P}_{p_1}^{-1}((1-t)/m_S)=\mathcal{P}_{p_1}(1-(1-t)/m_S)$ for $t\in[0,1]$.\\

\noindent\textit{2.} By \eqref{eq:P(k/m_S)} if $\omega(p_1)=(a_1\cdots a_n\cdots)$
then $\mathcal{P}_{p_1}\left(\frac{1}{m_S}\right) = \lim_nL\left(S^{d^{n-1}}(a_n)\cdots S^{d^0}(a_1)\right)$. By using Lemma \ref{lemma:Lemma1}, it implies $\mathcal{P}_{p_1}(\frac{1}{m_S})=\omega^{-1}\left(S^{d^0}(a_1)\cdots S^{d^{n-1}}(a_n)\cdots\right)$ and thus $\omega\circ\mathcal{P}_{p_1}\left(\frac{1}{m_S}\right) = h\circ \omega(p_1)$. So, by using the first item, we have $\omega(\mathcal{P}_1(p_1))=h(\omega(p_1))$. It follows that $\mathcal{P}_{1}(p_1)=p_2$ if and only if $\omega(p_2)=h(\omega(p_1))$. We also deduce that $h(\omega(p_1))=(\omega\circ\mathcal{P}_1\circ\omega^{-1})(\omega(p_1))$, which implies Item \textit{4}.\\

\noindent\textit{3.} If $\mathcal{P}_{1}(p_1)=p_2$, then $\mathcal{C}_{p_1}\cap\mathcal{C}_{p_2}\neq\emptyset$, thus $\mathcal{C}_{p_1}=\mathcal{C}_{p_2}$ since these are connected components of $J(f)$. In this case, $\mathcal{P}_{p_1}$ and $\mathcal{P}_{p_2}$ parametrize the same curve. So, $\mathcal{P}_{1}(p_1)=p_2$ implies that by traveling (in the non-trigonometric sense) along the path $\mathcal{P}_{p_2}$ from the point \textcolor{black}{$p_2$}, we land in $\pi^{-1}(1)$ at the point $p_1$. This is exactly the meaning of the equality $\mathcal{L}_1(p_2)=p_1$. Conversely, by the same arguments, if $\mathcal{L}_1(p_2)=p_1$ it means that $\mathcal{C}_{p_1}=\mathcal{C}_{p_2}$ and thus by traveling (in the trigonometric sense) along the path $\mathcal{P}_{p_1}$ from the point \textcolor{black}{$p_1$}, we land in $\pi^{-1}(1)$ exactly at the point $p_2$ i.e. $\mathcal{P}_1(p_1)=p_2$.\\

\noindent\textit{5.} Let us first assume that $\mathcal{C}_{p_1}=\mathcal{C}_{p_2}$. Since $\mathcal{C}_{p_1}$ is a simple closed curve winding above $\mathbb{S}^1$ a finite number of times, by traveling along this curve (in the trigonometric sense) from the point $p_1$ we will return in the fiber $\pi^{-1}(1)$ a certain number of times, say $n\geq1$, until we will reach the point $p_2$ (if $p_1=p_2$ we travel along all the curve until we return to the starting point $p_1$, hence $n\geq1$). Then by definition of $\mathcal{P}_1$ we have $\mathcal{P}_1^{\circ n}(p_1)=p_2$. Thus by the fourth item of this lemma we deduce $\omega(p_2)=\omega\circ\mathcal{P}_1^{\circ n}(p_1)=h^n\circ\omega(p_1)$, thus $\omega(p_2)\in\mathrm{orb}_h(\omega(p_1))$.

Let us assume now that $\omega(p_2)\in\mathrm{orb}_h(\omega(p_1))$. Since $h$ is of finite order by construction (from \eqref{eq:homeoh} we observe that $h$ is of order $\leq$ the order of $S$), there exists $n\geq 0$ such that $\omega(p_2)=h^n\circ\omega(p_1)$. By using again the preceding point of this lemma we deduce $p_2=\mathcal{P}_1^{\circ n}(p_1)\in\mathcal{C}_{p_1}$, and thus $\mathcal{C}_{p_1}=\mathcal{C}_{p_2}$ since these are connected components of $J(f)$.\qed\\

We can now introduce the mappings $\mathcal{P}$ and $\mathcal{L}$ as follows:

\begin{equation}\label{eq:applicationPandL}
\left|
    \begin{array}{ll}
      J(f)  & \overset{\mathcal{P}}{\longrightarrow}\ \mathcal{S}_{\mathcal{A}}\\
      p=(e^{2i\pi t},w),t<1 & \longmapsto\ \left(t,\omega(\mathcal{P}_t(p))\right)
    \end{array}
\right.
\mathrm{;}\
\left|
    \begin{array}{ll}
      \mathcal{S}_{\mathcal{A}} & \overset{\mathcal{L}}{\longrightarrow}\ J(f)\\
      (t,\omega(p)) & \longmapsto\ \mathcal{L}_t(p)
    \end{array}
\right.
\end{equation}

\begin{lemme}\label{lemma:applicationsPetL}
    The following points hold:
    \begin{enumerate}
        \item The maps $\mathcal{P}$ and $\mathcal{L}$ are well defined.
        \item The map $\mathcal{L}$ is continuous.
        \item We have $\mathcal{P}\circ\mathcal{L}=\mathrm{Id}_{\mathcal{S}_{\mathcal{A}}}$ and $\mathcal{L}\circ\mathcal{P}=\mathrm{Id}_{J(f)}$.\\
    \end{enumerate}
\end{lemme}

\noindent\textbf{\underline{Proof of Theorem \ref{thm:Returnmapisahomeo}$:$}} The third item of Lemma \ref{lemma:applicationsPetL} shows that $\mathcal{L}:\mathcal{S}_{\mathcal{A}}\to J(f)$ is bijective. Since the spaces $\mathcal{S}_{\mathcal{A}}$ and $J(f)$ are compact, and since $\mathcal{L}$ is continuous by the second item of Lemma \ref{lemma:applicationsPetL},  we can conclude that $\mathcal{L}:\mathcal{S}_{\mathcal{A}}\to J(f)$ is a homeomorphism.\qed\\

\noindent\textbf{\underline{Proof of Lemma \ref{lemma:applicationsPetL}$:$}}\ \\
\noindent\textit{1.} The map $(t,p)\mapsto \mathcal{P}_t(p)$ is well defined and $\mathcal{P}$ is the composition of this map with the quotient map $[0,1]\times\mathcal{C}_{\mathcal{A}}\to \mathcal{S}_{\mathcal{A}}$. Thus $\mathcal{P}$ is well defined. Let us prove now that $\mathcal{L}$ is well defined. Let $\omega(p_1)$ and $\omega(p_2)$ be two elements of $\mathcal{C}_{\mathcal{A}}$ such that $h(\omega(p_1))=\omega(p_2)$. By definition of the equivalent relation $\sim_h$, we only need to check that $\mathcal{L}_1(p_2)=\mathcal{L}_0(p_1)$. According to Lemma \ref{lemma:htour} we have $\mathcal{L}_1(p_2)=p_1$, and by Definition \ref{defn:returnmaps} of $\mathcal{L}_0$ we have $\mathcal{L}_0(p_1)=p_1$. So finally we have $\mathcal{L}_1(p_2)=\mathcal{L}_0(p_1)$ and $\mathcal{L}$ is well defined.\\

\noindent\textit{2.} Let us define $\widehat{\mathcal{L}}:[0,1]\times\mathcal{C}_{\mathcal{A}}\to J(f)$ by $\widehat{\mathcal{L}}:(t,p)\mapsto \mathcal{L}_t(p)$. Since $\mathcal{S}_{\mathcal{A}}$ is equipped with the quotient topology, it is sufficient to show that $\widehat{\mathcal{L}}$ is continuous. So let us fix $(t_n,x_n)$ a sequence of $[0,1]\times \mathcal{C}_{\mathcal{A}}$ which converges to $(t,x)$. We denote $p_n:=\omega^{-1}(x_n)$ and $p:=\omega^{-1}(x)$. Let $\mathcal{P}_p$ (resp. $\mathcal{P}_{p_n}$) be the curve parameterizing the connected component of $p$ (resp. $p_n$) in $J(f)$ given by Theorem \ref{thm:theoremedescourbes}. According to the first item of Lemma \ref{lemma:htour}, we have $\mathcal{L}_t(p)=\mathcal{P}_p\left(1-\frac{1-t}{m_S}\right)\ \mathrm{and}\ \mathcal{L}_{t_n}(p)=\mathcal{P}_{p_n}\left(1-\frac{1-t_n}{m_S}\right)$.
Thus we have:
\begin{align*}
    d_{\mathcal{M}}\left(\widehat{\mathcal{L}}(t,x),\widehat{\mathcal{L}}(t_n,x_n)\right) & \leq d_{\mathcal{M}}\left(\mathcal{P}_{p}\left(\frac{t}{m_S}\right),q_n\right) + d_{\mathcal{M}}\left(q_n,\mathcal{P}_{p_n}\left(1-\frac{1-t_n}{m_S}\right)\right) =:u_n+v_n,
\end{align*}
where $q_n:=\mathcal{P}_{p_n}\left(1-\frac{1-t}{m_S}\right)$. We have $\lim_nu_n=0$ by Proposition \ref{prop:lemmeequi2}, and we have $\lim_nv_n=0$ by Corollary \ref{cor:corollairelemmeequicontinuite}. Finally, $\lim_n\widehat{\mathcal{L}}(t_n,x_n)=\widehat{\mathcal{L}}(t,x)$ and thus $\mathcal{L}$ is continuous.\\

\noindent\textit{3.} Let $(t,p)$ be an element of $[0,1]\times(\{1\}\times J_1)$ and denote $\widetilde{p}:=\mathcal{L}_t(p)$. Let $\mathcal{P}_p$ and $\mathcal{P}_{\widetilde{p}}$ be the curves given by Theorem \ref{thm:theoremedescourbes}. By construction of $\widetilde{p}$, $\mathcal{P}_p$ and $\mathcal{P}_{\widetilde{p}}$ parametrize the same curve.
\begin{itemize}
    \item[\textbullet] If $t\in]0,1[$, by definition of $\mathcal{L}_t$, we deduce that by traveling (in the trigonometric sense) along $\mathcal{P}_{\widetilde{p}}$ from $\widetilde{p}$, we will land in $\pi^{-1}(1)$ at the point $p$. Moreover, by definition of $\mathcal{P}_t$, since $\mathcal{P}_p$ and $\mathcal{P}_{\widetilde{p}}$ parametrize the same curve we also have $\mathcal{P}_t(\widetilde{p})=p$. So we have $\mathcal{P}_t\circ\mathcal{L}_t(p)=p$ and thus $\mathcal{P}\circ\mathcal{L}(t,\omega(p))=(t,\omega(p))$.
    \item[\textbullet] If $t=0$ we have by definition $\mathcal{L}(0,\omega(p))=\mathcal{L}_0(p)=p$ and $\mathcal{P}_0(p)=p$, thus $\mathcal{P}\circ\mathcal{L}(0,\omega(p))=(0,\omega(p))$. 
    \item[\textbullet] If $t=1$ then $\mathcal{L}(1,\omega(p))=\mathcal{L}_1(p)$, but according to the fourth point of Lemma \ref{lemma:toursdanslabase} we have $\mathcal{L}_1(p)=\omega^{-1}\circ h^{-1}\circ\omega(p)$ and $\mathcal{P}_1=\omega^{-1}\circ h\circ \omega$, so we deduce that $\mathcal{P}_1\circ\mathcal{L}_1(p)=p$, and thus $\mathcal{P}\circ\mathcal{L}(1,\omega(p))=(1,\omega(p))$. 
\end{itemize}
We have then proved that $\mathcal{P}\circ\mathcal{L}=\mathrm{Id}_{\mathcal{S}_{\mathcal{A}}}$. By following these arguments we can also prove that $\mathcal{L}\circ\mathcal{P}=\mathrm{Id}_{J(f)}$.\qed   

\subsection{Isotopies of solid torus and proof of Theorem \ref{thm:Kisotopicf^-1(S^1times0)}}

\begin{defn}\label{defn:nsolidtorus}
    A $n-$solid torus with $n\geq1$, is an open subset $\mathbb{T}_n$ of $\mathbb{S}^1\times\cmplex$ such that there exists a strong deformation retract of $\mathbb{T}_n$ on a simple closed curve $C_n\subset\mathbb{T}_n$, such that $\pi:C_n\to\mathbb{S}^1$ is an unramified covering map of degree $n$.
\end{defn}

\begin{prop}\label{prop:isotopiessolidtorus}
    Let $n\geq 1$ and $m\geq 1$ be two integers, and let $\mathbb{T}_n$ be a $n-$solid torus. 
    Let $\alpha,\beta:[0,1]\to\mathbb{T}_n$ be two disjoint simple closed curves inside $\mathbb{T}_n$, both winding exactly $m$ times above $\mathbb{S}^1$. Then there exists $(h_t)_{t\in[0,1]}$ an ambient isotopy of $\overline{\mathbb{T}}_n$ deforming $\alpha$ on $\beta$ such that $h_t$ has the form $h_t(z,w)=(z,k_t(z,w))$, and such that $h_t=\mathrm{Id}$ on $\partial\mathbb{T}_n$.
\end{prop}
\noindent\textbf{\underline{Proof$:$}} It is a standard result that $\mathbb{T}_n$ is actually homeomorphic to the standard solid torus $\mathbb{T}:=\mathbb{S}^1\times\mathbb{D}$, thus it is sufficient to prove the statement for $\mathbb{T}$. Let $D\Subset\mathbb{D}$ be an open disc such that the curves $\alpha$ and $\beta$ are contained in $\mathbb{S}^1\times D$. Let $s\in[0,1]$ and denote $z=e^{2i\pi sm}\in\mathbb{S}^1$. We denote $\alpha(s)=:(z,\widehat{\alpha}(s))$ and $\beta(s)=:(z,\widehat{\beta}(s))$. Let $M_s:=\varphi_{\widehat{\beta}(s)}\circ\varphi_{\widehat{\alpha}(s)}$ where $\varphi_a$ is the Möbius transformation
\begin{equation}\label{eq:Mobius}
    \varphi_a(w)=\frac{a-w}{1-\overline{a}w},\ w\in\mathbb{D},\ a\in\mathbb{D}.
\end{equation}
The map $M_s$ sends $\widehat{\alpha}(s)$ to $\widehat{\beta}(s)$, and $(s,w)\in[0,1]\times\overline{D}\mapsto(s,M_s(w))$ is a homeomorphism on its image. Then there exists $\widehat{M}$ a self homeomorphism of $[0,1]\times\mathbb{D}_2$ of the form $\widehat{M}(s,w)=(s,g_s(w))$ such that $\widehat{M}=M$ on $[0,1]\times\overline{D}$ and $\widehat{M}=\mathrm{Id}$ on $[0,1]\times\mathbb{D}_2\backslash\mathbb{D}$, see \cite[Thm. 1.8]{GGH24}. Then observe that $H(z,w):=(z,g_s(w)),\ w\in\overline{\mathbb{D}}\ \mathrm{and}\ z=e^{2i\pi s},\ s\in[0,1]$, defines a homeomorphism of $\mathbb{T}$ such that $H(\alpha(s))=\beta(s)$, $s\in[0,1]$.

Let us now define a isotopy $(h_t)_{t\in[0,1]}$ of $\mathbb{T}$ such that $h_0=\mathrm{Id}$ and $h_1=H$ using the standard Alexander's trick: 
$$h_t(e^{2i\pi s},w):=\left(e^{2i\pi s},tg_s\left(\frac{w}{t}\right)\right),\ \mathrm{if}\ 0\leq|w|< t,$$
and
$$h_t(e^{2i\pi s},w):=(e^{2i\pi s},w),\ \mathrm{if}\ t\leq|w|\leq 1.$$
The isotopy $(h_t)_{t\in[0,1]}$ finally deforms $\alpha$ on $\beta$ and satisfies all the required property of the statement.\qed\\

\noindent\textbf{\underline{Proof of Theorem \ref{thm:Kisotopicf^-1(S^1times0)}$:$}} According to Lemma \ref{lemma:htour}, two fixed points $\omega(p_{j_1})$ and $\omega(p_{j_2})$ are in the same orbit under $h$ if and only if the points $p_{j_1}$ and $p_{j_2}$ are contained in the same connected component inside $J(f)$ i.e. $\mathcal{C}_{j_1}=\mathcal{C}_{j_2}$. Up to re-index, let us assume that $p_1,\cdots ,p_m$ (with $1\leq m\leq d$) are exactly the different fixed points such that the curves $(\mathcal{C}_{j})_{1\leq j\leq m}$ are disjoint two by two. Recall that our aim is to prove that the class of $\bigcup_{j=1}^{m}\mathcal{C}_j$ is equal to the class of the compact set $f^{-1}(\mathbb{S}^1\times\{0\})$ for the equivalent relation \eqref{eq:equivrelaambientisotopy}. 

Let us fix $j\in\{1,\cdots,m\}$. We denote $\mathcal{M}_{1j}$ the connected component of $\mathcal{M}_1$ which contains the loop $\mathcal{P}_{j,1}(t)$. According to Lemma \ref{lemma:Lemma1} and Proposition \ref{prop:Lemma2}, $\mathcal{M}_{1j}$ contains also the loop $\mathcal{P}_j(t)$. These two loops winds inside $\mathcal{M}_{1j}$ exactly $m_j$ times above $\mathbb{S}^1$, see again Proposition \ref{prop:Lemma2}, and we can observe that there exists a strong deformation retract of $\mathcal{M}_{1j}$ on $C_{j,1}:=\mathcal{P}_{j,1}([0,1])$ and on $\mathcal{C}_j$. So, $\mathcal{M}_{1j}$ is a $m_j-$solid torus in the sense of Definition \ref{defn:nsolidtorus}. Thus we can apply Proposition \ref{prop:isotopiessolidtorus} to the curves $\alpha(t):=\mathcal{P}_{j}(t)$ and $\beta(t):=\mathcal{P}_{1,j}(t)$ (with $n=m:=m_j$) to conclude that there exists an isotopy $h_t:\overline{\mathcal{M}}_{1j}\to\overline{\mathcal{M}}_{1j}$ of the form $h_t(z,w)=(z,k_t(z,w))$ deforming $\mathcal{P}_{j}$ into $\mathcal{P}_{j,1}$ such that $h_{j,t}=\mathrm{Id}$ on $\partial\mathcal{M}_{1j}$. We can then extend $h_{j,t}$ into an ambient isotopy of $\mathbb{S}^1\times\cmplex$ by putting for each $t\in[0,1]$:
\begin{equation}\label{eq:extensionisoh_jt}
    h_{j,t}:=\mathrm{Id}\ \mathrm{on}\ (\mathbb{S}^1\times\cmplex)\backslash\mathcal{M}_{1j}.
\end{equation}

Repeating these arguments for each $j$, we can define for each $t\in[0,1]$ a homeomorphism $H_t:=h_{1,t}\circ\cdots\circ h_{m,t}$ of $\mathbb{S}^1\times\cmplex$ of the form $H_t(z,w)=(z,K_t(z,w))$. At last, observe that $(H_{t})_{t\in[0,1]}$ is an ambient isotopy of $\mathbb{S}^1\times\cmplex$ which deform $\bigcup_{j=1}^{m}\mathcal{C}_j$ into $\bigcup_{j=1}^{m}\mathcal{C}_{j,1}=f^{-1}(\mathbb{S}^1\times\{0\})$ by deforming each curve $\mathcal{C}_j$ into the curve $\mathcal{C}_{1,j}$ independently from the other curves.\qed

\section{Algebraic braids and Julia sets}\label{sec:AlgBraidsandJuliasets}

We still consider the case $p(z)=z^d$ and we use the notations and the results obtained in Section \ref{sec:IMG} and Section \ref{sec:topologicaldescriptionofJ(f)}. Let $\Lambda\in\pi_0(\mathcal{D}')$ be fixed. Our goal in this section is to associate to each $f_{\lambda},\ \lambda\in\Lambda$, a canonical isotopy class of curves lying in $J(f)$ which corresponds to an isotopy class of an \textit{algebraic braid} of degree $d$ (see Definition \ref{defn:algbraids}) denoted $\mathrm{ab}(f_{\lambda})$. By structural stability of Julia sets in $\Lambda$ (see Theorem \ref{thm:isotopyoverJ_p}), we deduce that $\lambda\in\Lambda\mapsto \mathrm{ab}(f_{\lambda})$ is constant, see Lemma \ref{lemma:braidab(f_lambda)constantonLambda}. It yields a natural map from $\pi_0(\mathcal{D}')$ to $AB_d$, the set of algebraic braids of degree $d$. We also recall in this paragraph that there exists a canonical map $AB_d\to\mathrm{Con}(\mathfrak{S}_d)$ and we show that the composition $\pi_0(\mathcal{D}')\to\mathrm{Con}(\mathfrak{S}_d)$ is surjective, it gives a proof of Corollary \ref{cor:abcmap} announced in introduction. 

{\subsection{Closed braids and algebraic braids}\label{sec:closedbraids}}

We start by giving a definition of (closed) braids embedded in $\mathbb{S}^1\times\cmplex$. 

\begin{defn}\label{defn:closedbraids}
    A closed braid of $\mathbb{S}^1\times\cmplex$ of degree $d$ is a compact $K\subset\mathbb{S}^1\times\cmplex$ which has a structure of topological manifold of dimension $1$ such that $\pi:K\to\mathbb{S}^1$ is a (surjective) non ramified covering map of degree $d$, where $\pi:\mathbb{S}^1\times\cmplex\to\mathbb{S}^1$ is the first projection. 
\end{defn}

If $K\subset\mathbb{S}^1\times\cmplex$ is a closed braid of degree $d$, and if $K(1)=K\cap\pi^{-1}(1)$, then observe that there exists a group morphism (use that $\pi:K\to\mathbb{S}^1$ is a covering map):
\begin{equation*}
    {\rho_K}:\left\{
    \begin{array}{ll}
          \mathbb{Z} & \longrightarrow\ \mathrm{Aut}(K(1))\simeq\mathfrak{S}_d\\
          m & \longmapsto\ (x\mapsto\pi^{-1}(\gamma_m)[x])
    \end{array}
    \right.
\end{equation*}
where $\gamma_m(t)=e^{2i\pi tm}$ and where $\pi^{-1}(\gamma_m)[x]$ is the end point of the lift of $\gamma_m$ by $\pi$ starting at $x$. The connected component of $K$ are given by smooth paths of the form $t\in[0,1]\mapsto\pi^{-1}(\gamma_m)_x(t)$ with $m\in\mathbb{N}^*$ and $x\in K(1)$. Since $\mathfrak{S}_d$ is a finite group, for each connected component $C$ of $K$ one can chose $x\in K(1)$ and $m\in\mathbb{N}^*$ such that the path $t\in[0,1]\mapsto \pi^{-1}(\gamma_m)_x(t)$ is a loop based at $x$ parametrizing $C$. Finally, $K$ is a union of disjoint closed loops winding above $\mathbb{S}^1$ which may interlace, there are at most $d$ different loops ; it justifies the terminology of closed braid of degree $d$. We associate to $K$ a permutation $S_K$:
\begin{equation}\label{eq:permutationS_K}
    S_K:=\rho_K(d)\in\mathfrak{S}_d
\end{equation}
which we identify to an element of $\mathfrak{S}_d$ by identifying $K(1)$ to $\{1,\cdots,d\}$. In the sequel we will be interested in the class of $S_K$ modulo conjugacy and thus $S_K$ can be seen as a element of $\mathrm{Conj}(\mathfrak{S}_d)$, which does not depend on the choice made to identify $K(1)$ to $\{1,\cdots,d\}$. 

It is known in the literature that the notion of (closed) braid of degree $d$ is linked to the group $B_d$ of Artin. It is the group defined by generators and relations by:
\begin{equation}\label{eq:GroupeArtin}
    B_d:=\left\langle\sigma_1,\cdots,\sigma_{d-1}|\sigma_i\sigma_{i+1}\sigma_i=\sigma_{i+1}\sigma_i\sigma_{i+1}\ ;\ \sigma_i\sigma_j=\sigma_j\sigma_i\ \mathrm{if}\ |i-j|\geq2\right\rangle.
\end{equation}
To each closed braid $K\subset\mathbb{S}^1\times\cmplex$ of degree $d$ there exists a procedure to associate a word $\sigma(K)\in B_d$ to $K$. The definition of closed braids given above is invariant under ambient isotopies of $\mathbb{S}^1\times\cmplex$ respecting the product structure of $\mathbb{S}^1\times\cmplex$ i.e. isotopies of the form $H_t(z,w)=(z,h_t(z,w)),\ t\in[0,1]$. Moreover, the closed braid $\widetilde{K}$ obtained after a modification of $K$ by such an isotopy is associated with a word of $B_d$ of the form $\sigma(\widetilde{K})=w\sigma(K)w^{-1},\ w\in B_d$. Moreover, each element of $B_d$ can be viewed as the word associated with a closed braid in $\mathbb{S}^1\times\cmplex$. In other words, one can identify the set $CB_d$ of ambient isotopy classes (respecting the product structure of $\mathbb{S}^1\times\cmplex$) of (closed) braids of degree $d$ with the set $\mathrm{Conj}(B_d)$ of conjugacy classes of $B_d$, we refer to the article by Rudolph \cite{RL83} for an exposition of (closed) braids in $\mathbb{S}^1\times\cmplex$.

Each word $w\in B_d$ induce a permutation $S_w\in\mathfrak{S}_d$, for $w=\sigma_i$ a generator of $B_d$ the induced permutation is $S_{\sigma_i}=(i,i+1)$ a transposition, for a general word $w\in B_d$ we express $w$ as a product of generators and the composition of the associated transpositions give $S_w$. If $w=\sigma(K)$ comes from a closed braid in $\mathbb{S}^1\times\cmplex$, the permutation $S_K$ defined by \eqref{eq:permutationS_K} coincide with $S_{\sigma(K)}$ up to conjugacy i.e. $S_{\sigma(K)}=S_K$ in $\mathrm{Conj}(B_d)$. This gives us a map (which coincide with the map \eqref{eq:themapc:AB_dtoCB_d} announced in introduction):
\begin{equation}\label{eq:themapccorpsdutexte}
    \mathrm{c}:CB_d\longrightarrow\mathrm{Conj}(\mathfrak{S}_d).
\end{equation}

Let us now introduce \textit{algebraic braids} of degree $d$. Let $q(z,w)$ be a polynomial map of the form:
\begin{equation}\label{eq:monicpolynomialsq(z,w)versionadmissible}
    q(z,w)=a_d(z)w^d+a_{d-1}(z)w^{d-1}+a_{d-2}w^{d-2}+\cdots+a_0(z)\in\cmplex[z,w],
\end{equation}
where $a_j(z)\in\cmplex[z]$ and $a_d(z)\not\equiv0$. If $q(z,w)$ does not contains any factor of the form $(z-c)$ and no square factors in $\cmplex[z,w]$, then the first projection $\pi:\{q=0\}\to\cmplex$ is a covering map of degree $d$ ramified only above a finite number of points. Denoting 
$$K=\{q=0\}\cap(\mathbb{S}^1\times\cmplex),$$
if $\pi:K\to\mathbb{S}^1$ is not ramified then $K$ is a closed braid of degree $d$. Since this braid is obtained by using a polynomial mapping we call it an \textit{algebraic braid}. 

\begin{defn}\label{defn:algbraids}
    For any polynomial $q$ of the form \eqref{eq:monicpolynomialsq(z,w)versionadmissible} such that $\pi:\{q=0\}\cap(\mathbb{S}^1\times\cmplex)\to\mathbb{S}^1$ is unramified of degree $d$, we denote $\mathrm{ab}(q)$ the isotopy class of $K:=\{q=0\}\cap(\mathbb{S}^1\times\cmplex)$ (modulo ambient isotopies preserving the product structure of $\mathbb{S}^1\times\cmplex$). We denote $AB_d\subset CB_d$ the set of all these isotopy classes. We also denote $S_q:=\mathrm{c}(K)$ the permutation given by \eqref{eq:themapccorpsdutexte}.
\end{defn}

If $K=\{q=0\}\cap(\mathbb{S}^1\times\cmplex)$ is an algebraic braid of degree $d$, according to Rudolph \cite{RL83} the word $\sigma(K)$ associated to $K$ in the Artin group $B_d$ is \textit{quasipositive}. It means that $\sigma(K)$ is conjugated in $B_d$ to a word of the form $\prod_{k=1}^{n}w_k\sigma_{i(k)}^{+1}w_k^{-1}$ where the indices $1\leq i(k)\leq d-1$ are distinct two by two. In particular, the set $AB_d$ is not equal to $CB_d$ i.e. the algebraic braids do not represent all the possible closed braids given by Definition \ref{defn:closedbraids}. For $d=2$, we have $B_2=\{\sigma_1^n,\ n\in\mathbb{Z}\}$ and $AB_2=\{\sigma_1^{+n},\ n\geq0\}$, see Section \ref{sec:comparaison}.

{\subsection{Braids and Julia sets: proof of Theorem \ref{thm:algebraicbraidsinJuliasets} and Corollary \ref{cor:abcmap}}\label{sec:proofofcorollary1.3}}

Let us fix $\Lambda\in\pi_0(\mathcal{D}')$. Our goal in this paragraph is to introduce a canonical isotopy class $\mathcal{K}(f_{\lambda})\in AB_d$ for each $\lambda\in\Lambda$. First observe that if $\lambda_0\in\Lambda$ is an admissible parameter (see Definition \ref{defn:admissileparameters}) then the set 
$$K_{\lambda_0}:=\{q_{\lambda_0}=0\}\cap(\mathbb{S}^1\times\cmplex)=f^{-1}_{\lambda_0}(\mathbb{S}^1\times\{0\})$$
may be not an algebraic braid if the first projection $\pi:K_{\lambda_0}\to\mathbb{S}^1$ is ramified. But if $[\lambda_0]\not\in E$ with $E$ defined by \eqref{defn:EetlesE_z_version_intro}, then $\pi:K_{\lambda_0}\to\mathbb{S}^1$ is unramified and the isotopy class of $K_{\lambda_0}$ gives an algebraic braid $\mathrm{ab}(f_{\lambda_0}):=\mathrm{ab}(q_{\lambda_0})\in AB_d$. We also recall that $\mathrm{ab}(f_{\lambda_0})$ is the isotopy class of the union $\mathcal{C}(f_{\lambda_0})=\bigcup_{j=1}^{d}\mathcal{C}_j$ of connected components of fixed points of $f_{\lambda_0}|_{\{1\}\times J_1(f_{\lambda_0})}$ according to Theorem \ref{thm:Kisotopicf^-1(S^1times0)}. We also recall that there exists an admissible parameter $\lambda_0\in\Lambda$ satisfying $[\lambda_0]\not\in E$ according to Proposition \ref{prop:lambdaisadmissible}.\\ 

Let us recall that according to Theorem \ref{thm:isotopyoverJ_p}, we have the structural stability of Julias sets in the hyperbolic component $\Lambda$. In particular, there exists $U\subset\Lambda$ a neighborhood of $\lambda_0$ and there exists $h:U\times J(f_{\lambda_0})\to\mathbb{S}^1\times\cmplex$ a continuous map such that:
\begin{enumerate}
    \item[a)] $h$ has the form $h(\lambda,z,w)=(\lambda,z,g_{\lambda}(z,w))=:(\lambda,h_{\lambda}(z,w))$,
    \item[b)] $h_{\lambda}:J(f_{\lambda_0})\to J(f_{\lambda})$ is a homeomorphism and $h_{\lambda}\circ f_{\lambda_0}=f_{\lambda}\circ h_{\lambda}$ on $J(f_{\lambda_0})$.
\end{enumerate}

\begin{lemme}\label{lemme:imagedetresses}
    Let $\lambda_1\in\Lambda$ and let $t\in[0,1]\mapsto\lambda_t$ be a continuous path inside $U$ joining $\lambda_0$ to $\lambda_1$. For each $t\in[0,1]$, denote $\mathcal{K}_t$ the ambient isotopy class (in the sense of \eqref{eq:equivrelaambientisotopy}) of $h_{\lambda_t}(\mathcal{C}(f_{\lambda_0}))$. Then for each $t\in[0,1]$, $h_t(\mathcal{C}(f_{\lambda_0}))$ is a closed braid of degree $d$ and $\mathcal{K}_0=\mathcal{K}_t$.
\end{lemme}
\noindent\textbf{\underline{Proof$:$}} Let $(\alpha_j)_{1\leq j\leq m}$ be the different disjoint closed curves which compose $\mathcal{C}(f_{\lambda})$. For $1\leq j\leq m$, we denote $\alpha_j(s)=:(e^{2i\pi sn_j},\widehat{\alpha}_j(s))$, where $n_j$ is the number of times $\alpha_j$ winds above $\mathbb{S}^1$. Since $h_t$ has the form $h_t(z,w)=(z,g_{t}(z,w))$, the first assertion of the lemma is straightforward. Let us prove the second assertion of the lemma i.e. $\mathcal{C}(f_{\lambda_0})$ and $h_t(\mathcal{C}(f_{\lambda_0}))$ are ambient isotopic. Since $\bigcup_{t\in[0,1]}h_t(\mathcal{C}(f_{\lambda_0}))$ is a compact subset of $\mathbb{S}^1\times\cmplex$, there exists $R'>0$ such that $\bigcup_{t\in[0,1]}h_t(\mathcal{C}(f_{\lambda_0}))\subset\mathbb{S}^1\times\mathbb{D}_{R'}$. Without loss of generality, we can assume that $R'=1$ in what follows. 

There exists $D\Subset\mathbb{D}$ such that for every $t\in[0,1]$ we have $\pi'\circ h_t(\mathcal{C}(f_{\lambda_0}))\subset D$, where $\pi'$ denotes the second projection of $\mathbb{S}^1\times\cmplex$. For $\varepsilon>0$ and for $1\leq j\leq m$, we define $V_{j}(\varepsilon):=\bigcup_{s\in[0,1]}\{e^{2i\pi sn_j}\}\times D_s(\alpha_j)$, where $D_s(\alpha_j)\Subset D$ is the open disc of radius $\varepsilon>0$ centered at $\widehat{\alpha}_j(s)$. We can choose $\varepsilon>0$ small enough such that the sets $(\overline{V}_{j}(\varepsilon))_{1\leq j\leq m}$ are disjoint two by two. 
 
Let $1\leq j\leq m$ be fixed. Denote $\alpha_{j,t}(s):=h_t\circ \alpha_j(s)$, which has the form $\alpha_{j,t}(s)=(e^{2i\pi s n_j},\widehat{\alpha}_{j,t}(s))$. For each $(t,s)\in[0,1]^2$, let $M_{s}^{j,t}:=\varphi_{\widehat{\alpha}_{j,t}(s)}\circ\varphi_{\widehat{\alpha}_j(s)}$, where $\varphi_a$ is for each $a\in\mathbb{D}$ the Möbius transformation defined by \eqref{eq:Mobius}. Then $M_j(t,s,w):=(t,s,M_{s}^{j,t}(w))$ is a self homeomorphism of $[0,1]^2\times\overline{\mathbb{D}}$. According to \cite[Thm. 1.8]{GGH24}, there exists $\widehat{M}_j$ a homeomorphism of $[0,1]^2\times\overline{\mathbb{D}}$ such that $\widehat{M}_j=M_j$ on $[0,1]^2\times\overline{D}$, such that $\widehat{M}_j=\mathrm{Id}$ on $[0,1]^2\times\overline{\mathbb{D}}\backslash\overline{D}$, and such that $\widehat{M}_j$ has the fibered form $\widehat{M}_j(t,s,w)=(t,s,L_s^{j,t}(w))$. We can then define for each $t\in[0,1]$, $h_{j,t}$ a self homeomorphism of $\mathbb{S}^1\times\overline{\mathbb{D}}$ by:
$$h_{j,t}(e^{2i\pi s},w):=(e^{2i\pi s},L_{s}^{j,t}(w)).$$
Observe that $(h_{j,t'})_{t'\in[0,t]}$ is an isotopy deforming $\alpha_j$ on $\alpha_{j,t}$ for each $t\in[0,1]$, which is equal to $\mathrm{Id}$ on $\mathbb{S}^1\times\partial\mathbb{D}$. In particular, $h_{j,t}$ extends to $\mathbb{S}^1\times\cmplex$ by putting $h_{j,t}=\mathrm{Id}$ on $\mathbb{S}^1\times(\cmplex\backslash\mathbb{D})$. Let $\{\tau_{j,s}\}_{s\in[0,1]}$ be a family of homeomorphisms of $\cmplex$, depending continuously on $s$, such that $\tau_{j,s}(D_s(\alpha_j))=\mathbb{D}$ and $\tau_{j,s}(\widehat{\alpha}_j(s))=\widehat{\alpha}_j(s)$. Then define a new isotopy $\widetilde{h}_{j,t}$ of $\mathbb{S}^1\times\cmplex$ as follows:
$$\widetilde{h}_{j,t}(e^{2i\pi s},w):=h_{j,t}(e^{2i\pi s},\tau_{j,s}(w)).$$
It satisfies $\widetilde{h}_{j,t}(e^{2i\pi s},w)=(e^{2i\pi s},w)$ if $w\notin D_{s}(\alpha_j)$, thus $\widetilde{h}_{j,t}=\mathrm{Id}$ on $(\mathbb{S}^1\times\cmplex)\backslash V_{j}(\varepsilon)$. Moreover, since $\tau_{j,s}$ fix $\widehat{\alpha}_j(s)$, we have $\widetilde{h}_{j,t}\circ\alpha_j(s)=h_{j,t}\circ\alpha_j(s)=\alpha_{j,t}(s)$ by construction of $h_{j,t}$.

Since the ${\overline{V}_j(\varepsilon)}'$s are disjoint two by two, and since each isotopy $\widetilde{h}_{j,t}$ is equal to $\mathrm{Id}$ on the complement of $V_j(\varepsilon)$, we obtain a well define ambient isotopy $\widetilde{H}_t:\mathbb{S}^1\times\cmplex\to\mathbb{S}^1\times\cmplex$ by defining:
$$\widetilde{H}_t(z,w):=(z,w),\ \mathrm{if}\ (z,w)\notin\bigcup_{j=1}^nV_j,$$
and 
$$\widetilde{H}_t(z,w):=\widetilde{h}_{j,t}(z,w),\ \mathrm{if}\ (z,w)\in V_j(\varepsilon)\ \mathrm{for\ some}\ j.$$
Moreover, this isotopy has the form $\widetilde{H}_t(z,w)=(z,\widetilde{K}_t(z,w))$, and satisfies $\widetilde{H}_t\circ\alpha_j(s)=\widetilde{h}_{j,t}\circ\alpha_j(s)=\alpha_{j,t}(s)$. It proves that $\mathcal{K}_0=\mathcal{K}_t$ for each $t$.\qed\\ 

\begin{lemme}\label{lemma:braidab(f_lambda)constantonLambda}
    Let $\lambda_0\in\Lambda$ be an admissible parameter such that $[\lambda_0]\not\in E$. Then for every $\lambda\in\Lambda$, $f_{\lambda}:\{1\}\times J_1(f_{\lambda})\to\{1\}\times J_1(f_{\lambda})$ admits exactly $d$ different fixed points and the ambient isotopy class $\mathcal{K}(f_{\lambda})$ (in the sense of \eqref{eq:equivrelaambientisotopy}) of the union $\mathcal{C}(f_{\lambda})$ of the connected components in $J(f_{\lambda})$ of these fixed points, satisfies $\mathcal{K}(f_{\lambda})=\mathrm{ab}(f_{\lambda_0})$ and $\mathcal{C}(f_{\lambda})$ is a closed braid of degree $d$. In the sequel we denote $\mathrm{ab}(f_{\lambda}):=\mathcal{K}(f_{\lambda})=\mathrm{ab}(f_{\lambda_0})$. 
\end{lemme}
\noindent\textbf{\underline{Proof$:$}} Let $U$ and let $(h_\lambda)_{\lambda\in U}$ be the neighborhood of $\lambda_0$ inside $U$ and the associate structural isotopy of Julia sets satisfying the two points a) and b) given above. Let us fix $\lambda_1\in U$ and let $(\lambda_t)_{t\in[0,1]}$ be a continuous path joining $\lambda_0$ to $\lambda_1$ inside $U$. Observe that because $h_{\lambda_1}:J(f_{\lambda_0})\to J(f_{\lambda_1})$ is a homeomorphism such that $h_{\lambda_1}\circ f_{\lambda_0}=f_{\lambda_1}\circ h_{\lambda_1}$, the map $f|_{\{1\}\times J_1(f_{\lambda_1})}$ has also exactly $d$ different fixed points, and $\mathcal{C}(f_{\lambda_1}):=h_{\lambda_1}(\mathcal{C}(f_{\lambda_0}))$ is the union of the connected components of these fixed points inside $J(f_{\lambda_1})$. According to Lemma \ref{lemme:imagedetresses}, $\mathcal{C}(f_{\lambda_1})$ is a closed braid of degree $d$ whose ambient isotopy class $\mathcal{K}(f_{\lambda_1})$ is equal to the ambient isotopy class $\mathcal{K}(f_{\lambda_0})$ of $\mathcal{C}(f_{\lambda_0})$. Moreover, According to Theorem \ref{thm:Kisotopicf^-1(S^1times0)}, $\mathcal{K}(f_{\lambda_0})=\mathrm{ab}(f_{\lambda_0})$ thus we have $\mathcal{K}(f_{\lambda_1})=\mathrm{ab}(f_{\lambda_0})$.

For any over continuous path $(\widetilde{\lambda}_t)_{t\in[0,1]}$ joining $\widetilde{\lambda}_0=\lambda_0$ to some $\widetilde{\lambda}_1$ inside $\Lambda$, one can cover this path by a finite number of open sets $(U_{i})_{0\leq i\leq N}$ such that each of them satisfies the content of Theorem \ref{thm:isotopyoverJ_p}, and by applying the preceding arguments to the path restricted to each $U_{i}$, one can deduce that $\mathcal{K}(f_{\widetilde{\lambda}_1})=\mathrm{ab}(f_{\lambda_0})$.\qed\\

In consequence, we have proved the third item of Theorem \ref{thm:algebraicbraidsinJuliasets} and the existence of the map \eqref{eq:themapab:pi_0(D')_dtoAB_d} $\mathrm{ab}:\pi_0(\mathcal{D}')\to{AB}_d$ announced in Corollary \ref{cor:abcmap}. The two first item of Theorem \ref{thm:algebraicbraidsinJuliasets} are proved in Proposition \ref{prop:lambdaisadmissible} and Theorem \ref{thm:CompConnexesJ(f)}, the proof of Theorem \ref{thm:algebraicbraidsinJuliasets} is then complete. 

Our goal is now to prove Corollary \ref{cor:abcmap}.
We can compose $\mathrm{ab}:\pi_0(\mathcal{D}')\to{AB}_d$ with the map \eqref{eq:themapccorpsdutexte} $\mathrm{c}:CB_d\to\mathrm{Conj}(\mathfrak{S}_d)$ introduced just above, to get the map \eqref{eq:themapab:AB_dtoCB_d} $\mathrm{cab}:\pi_0(\mathcal{D}')\to\mathrm{Conj}(\mathfrak{S}_d)$. To prove Corollary \ref{cor:abcmap}, it remains to prove that $\mathrm{cab}:\pi_0(\mathcal{D}')\to\mathrm{Conj}(\mathfrak{S}_d)$ is surjective. To do so, it is enough to prove that for each product $C_1\circ\cdots\circ C_m$ of cycles with disjoint supports in $\mathfrak{S}_d$ can be realized as a permutation of the roots $\{w\in\cmplex :q_{\lambda,1}(w)=0\}=\{x_1\cdots x_d\}$ (with an arbitrary labeling) for a certain parameter $\lambda\in\mathcal{D}'$. We prove this fact in Proposition \ref{lemma:particularexamples} below. Indeed, it implies Corollary \ref{cor:abcmap} since the products of cycles with disjoint supports generate the conjugacy classes of $\mathfrak{S}_d$.

\begin{prop}\label{lemma:particularexamples}
    Let $d=\delta+d_1+\cdots+d_{m},\ m\in\mathbb{N},$ be a decomposition of $d$ with $\delta\geq0$ and
    $d_j\geq2$ being integers. Let $f(z,w)=(z^d,q_z(w))$ be a regular skew-product of degree $d$ defined by
    $$q_z(w):=\left(w^{\delta}-R_0^{\delta}z^{\delta}\right)\prod_{j=1}^{m}(w^{d_j}-R_{j}^{d_j}z),\ R_0, R_{j}>0.$$
    If $\delta=0$ remove the first factor, if $m=0$ remove the second factor. Then the ${R_{j}}'$s can be chosen appropriately such that there exists $\Lambda\in\pi_0(\mathcal{D}')$ and $\lambda\in\Lambda$ satisfying $[\lambda]\not\in E$, $\lambda$ is admissible, $f=f_{\lambda}$ and $:$ 
    \begin{enumerate}
        \item[\textbullet] If $m\geq1$, $\mathrm{cab}(\Lambda)=S_{q_{\lambda}}=C_1\circ\cdots\circ C_{m}\ \mathrm{mod.\ conj.}$ where
        \begin{equation}\label{eq:C_jcycle}
            C_j=(\delta+d_1+\cdots+d_{j-1}+1,\cdots,\delta+d_1+\cdots+d_{j})
        \end{equation}
        is a cycle of length $d_j$ in $\mathfrak{S}_d$.
        \item[\textbullet] If $m=0$, $\mathrm{cab}(\Lambda)=S_{q_{\lambda}}=\mathrm{Id}$.
    \end{enumerate}
    In particular, $\mathrm{cab}:\pi_0(\mathcal{D}')\to\mathrm{Conj}(\mathfrak{S}_d)$ is surjective.
\end{prop}

\noindent\textbf{\underline{Proof$:$}} Denote $\mathcal{R}:=(R_0,R_{1},\cdots,R_{m})\in\mathbb{C}^{m+1}$. Since each $d_j\geq 2$, observe that the coefficient of $w^{d-\delta-1}$ in $\prod_{j=1}^{m}(w^{d_j}-R_{j}^{d_j}z)$ is equal to $0$, and thus the coefficient of $w^{d-1}$ in $q_z(w)$ is also equal to $0$. So, there exists $\lambda=\lambda(\mathcal{R})\in\cmplex^{D_d}$ such that $f=f_{\lambda}$ since $q_z$ is a monic polynomial map. 
Observe that if $R_{m}>\cdots> R_{0}$, then the polynomial mappings $(w^{\delta}-R_0^{\delta}z^{\delta})$
and $(w^{d_j}-R_{j}^{d_j}z)_{1\leq j\leq m}$ do not have any common roots (for any $z\in\mathbb{S}^1$), thus we deduce $P(\lambda,z)\neq0$ for any $z\in\mathbb{S}^1$ i.e. $[\lambda]\not\in E$. Moreover, $\lambda(\mathcal{R})\in\mathcal{D}$ if $||\mathcal{R}||_{\cmplex^{m+1}}$ is large enough. From these remarks, we deduce that $\lambda=\lambda(\mathcal{R})\in\mathcal{D}'$ with $[\lambda]\not\in E$ if $\mathcal{R}$ is chosen appropriately. 

According to Proposition \ref{prop:lambdaisadmissible}, up to increase $||\mathcal{R}||_{\cmplex^{m+1}}$, we can assume that $\lambda=\lambda(\mathcal{R})\in\mathcal{D}'$ and and that $\lambda$ is admissible in the sense of Definition \ref{defn:admissileparameters}. In particular, the theory developed in Section \ref{sec:IMG} applies to $f$, which is a partial self-cover on $\mathbb{S}^1\times\mathbb{D}_R$ for some $R>0$, and we are allowed to compute its iterated monodromy and the associated permutation $S_{q}$. We should assume $\delta\geq1$ and $m\geq1$, the cases $\delta=0$ and $\delta=d$ can be treated similarly. 

We compute for $z=1$ the roots of $q_z$:
\begin{enumerate}
    \item[\textbullet] The roots of $w^\delta-R_{0}^{\delta
    }$ \ \ \ : $x^{1,j}:=R_0e^{2i\pi\frac{k-1}{\delta}},\ k\in\{1,\cdots,\delta\}$.
    \item[\textbullet] The roots of $w^{d_j}-R_{2,j}^{d_j}$ : $x_k^{2,j}:=R_{j}e^{2i\pi\frac{k-1}{d_j}},\ k\in\{1,\cdots,d_j\}$.
\end{enumerate}
The roots of $q_1$ are thus given by:
$$x_1:=x^{1,1},\cdots,x_{\delta}:=x^{1,\delta}\ \mathrm{and}$$
$$x_{\delta+1}=x_{1}^{2,1},\cdots,x_{\delta+d_1}:=x_{d_1}^{2,1}\ ;\ \cdots\ ;\ x_{d-d_m+1}:=x_1^{2,m},\cdots,x_d:=x_{d_m}^{2,m}.$$
This labeling $\{x_1,\cdots,x_d\}$ of $\{q_{1}=0\}$ then allows us to define $L:\mathcal{A}^*\to X$ as in Definition \ref{defn:codageL}, with $X=\bigcup_{n=0}^{+\infty}f_1^{-n}(1,0)$ where $f_1$ is the restriction of $f$ on $\{1\}\times\mathbb{D}_R$. In particular, the permutation $S_{q}$ associated to $f$ is equal (modulo conjugacy) to $S:j\in\{1,\cdots,d\}\mapsto L^{-1}\left(x_j^{\gamma_d}\right)$, see Definition \ref{defn:definitiondeS}. Let us know compute explicitly $S$. For each $t\in\reels$ and for each $1\leq j\leq d_j$ let us denote:
$$\alpha_{1,j}(t):=\left(e^{2i\pi t}, R_{0}e^{2i\pi t}\right)\ \mathrm{and}\ \alpha_{2,j}(t):=\left(e^{2i\pi t}, R_{j}e^{2i\pi\frac{t}{d_j}}\right).$$
Recall that $\gamma_d(t)=(e^{2i\pi dt},0)$ in Section \ref{sec:IMG}, and observe that $f\circ\alpha_{l,j}(t)=\gamma_d(t)$, $\alpha_{1,j}(0)=x^{1,j}$ and $\alpha_{2,j}(k-1)=x_{k}^{2,j}$. Thus by uniqueness of lifts starting at the same point by the covering map $f:f^{-1}(\mathbb{S}^1\times\mathbb{D}_R)\to\mathbb{S}^1\times\mathbb{D}_R$, we deduce that for each $1\leq j\leq\delta$:
$$\alpha_{1,j}(t)=f^{-1}(\gamma_d)_{x^{1,j}}(t),\ \forall t\in[0,1],$$
and that for each $1\leq j\leq m$ and for each $1\leq k\leq d_j$:
$$\alpha_{2,j}(t)=f^{-1}(\gamma_d)_{x_k^{2,j}}(t),\ \forall t\in[k-1,k].$$
In particular, its implies that: 
$$\{q=0\}\cap(\mathbb{S}^1\times\cmplex)=\left(\bigsqcup_{j=1}^{\delta}\alpha_{1,j}([0,1])\right)\bigsqcup\left(\bigsqcup_{j=1}^{m}\alpha_{2,j}([0,d_j])\right),$$
and that $S$ is given by:
$$S(1)=1,\cdots, S(\delta)=\delta,\ \mathrm{if}\ \delta>0,$$
and for each $1\leq j\leq m$ by:
$$S(\delta+d_1+\cdots+d_{j-1}+k)=\delta+d_1+\cdots+d_{j-1}+k+1,\ \forall k\in\{1,\cdots,d_j-1\}.$$
We deduce that $S=C_1\circ\cdots\circ C_m$ where each $C_j$ is defined by \eqref{eq:C_jcycle}. At last, we obtain $\mathrm{cab}(\Lambda)=S_q=S=C_1\circ\cdots\circ C_m\ \mathrm{mod.\ conj.}$\qed

{\section{The set \texorpdfstring{$AB_2$}{TEXT} and Astorg-Bianchi classification}\label{sec:comparaison}}

We would like to compare here our result to the one of Astorg-Bianchi \cite{AstBian23} in the case $d=2$ and $p(z)=z^2$. To do so, we give some more information concerning algebraic braids. For $d=2$, the braid group $B_d=B_2$ (see \eqref{eq:GroupeArtin}) becomes a free group generating by one element:
$$B_2=\{\sigma_1^n,\ n\in\mathbb{Z}\}.$$
As explain above, the main result of Rudolph \cite{RL83} said that if $K=\{q=0\}\cap(\mathbb{S}^1\times\cmplex)$ is an algebraic braid of degree $d$, then the word $\sigma(K)$ associated to $K$ (up to conjugacy) in $B_d$ is \textit{quasipositive} i.e. it is a word of the form $\prod_{k=1}^{n}w_k\sigma_{i(k)}^{+1}w_k^{-1}$ where the indices $1\leq i(k)\leq d-1$ are distinct two by two, and where $w_k\in B_d$. Rudolph also proved in \cite{RL83} that for any quasipositive word $w\in B_d$, there exists a polynomial map $q$ such that $\{q=0\}\cap(\mathbb{S}^1\times\cmplex)$ is represented by $w$ in $B_d$. Applying these results to $d=2$ it gives:
$$AB_2=\{\sigma_1^{+n},\ n\geq0\}.$$

Let us recall the results obtained by Astorg-Bianchi \cite[Theorem B]{AstBian23} mentioned in introduction for $d=2$ and $p(z)=z^2$. For each $\lambda=(a,b,c)\in\cmplex^{D_2}=\cmplex^3$, define the number:
$$s(\lambda):=\#\{z\in\mathbb{D}\ :\ a^2z^2+b^2z+c^2=0\},$$
it is not difficult to see that $\{s(\lambda),\ \lambda\in\mathcal{D}'\}=\{0,1,2\}$, and thus the set \eqref{eq:shift_p} mentioned above is identified with the set $\{0,1,2\}$. In other words, $s$ induces a bijective map:
$$
s:\left|
\begin{array}{ll}
  \pi_0(\mathcal{D}') & \overset{\simeq}{\longrightarrow}\ \{0,1,2\}\\
  \textcolor{white}{} & \textcolor{white}{}  \textcolor{white}{}\\
  \Lambda\ni\lambda   & \longmapsto\ s(\lambda)\\
\end{array}
\right.
$$
The link between braids and the mappings $\mathrm{ab}:\pi_0(\mathcal{D}')\to AB_2$ and $\mathrm{c}:\pi_0(\mathcal{D}')\to\mathfrak{S}_2$ is the following. Let us state the following lemma from \cite{AstBian23}:\\


\begin{lemme}[Astorg-Bianchi {\cite[Lem. 5.11]{AstBian23}}]\label{lemma:lemmeAB}
    Let $\lambda=(a,b,c)\in\mathcal{D}'$ such that $[\lambda]\not\in E$, and denote $K=\{q_{\lambda}=0\}\cap(\mathbb{S}^1\times\cmplex)$. If the norm of $\lambda$ is large enough, then $K$ is an algebraic braid of degree $2$ and the the following points hold:
    \begin{enumerate}
        \item If $s(\lambda)\in\{0,2\}$ then $K$ is composed of two connected components parametrized by closed loops winding above $\mathbb{S}^1$ exactly one time, and the linking number of these two loops is $\mathrm{Wind}(K)=\frac{s(\lambda)}{2}$.
        \item If $s(\lambda)=1$ then $K$ is connected and is parametrized by a loop winding exactly two times above $\mathbb{S}^1$. 
    \end{enumerate}
\end{lemme}
Since $AB_2=\{\sigma_1^n,\ n\geq0\}$ we deduce that all the informations given by Astorg-Bianchi results and by the maps $\mathrm{ab}:\pi_0(\mathcal{D}')\mapsto AB_2$ and $\mathrm{c}:\pi_0(\mathcal{D}')\to\mathfrak{S}_2$ can be gathered in the following table. This table should be read with $\lambda\in\mathcal{D}'$ of large norm and such that $[\lambda]\not\in E$, so that Lemma \ref{lemma:lemmeAB} applies.

\begin{center} 
\footnotesize
\begin{tabular}{|>{\centering\arraybackslash}p{1.5cm}||>{\centering\arraybackslash}p{7.5cm}|>{\centering\arraybackslash}p{2.5cm}|>{\centering\arraybackslash}p{2.5cm}|}  
\hline
\multicolumn{4}{|c|}{$\phantom{0}^{\phantom{0}}$Comparison with Astorg-Bianchi classification for \( d = 2 \) and \( p(z) = {z^{2^{\phantom{0}^{\phantom{0}}}}}_{\phantom{0}} \)$\phantom{0}^{\phantom{0}}$} \\
\hline
\begin{center}\( s(\lambda) \)\end{center} & \begin{center} \( \{q_{\lambda} = 0\} \cap (\mathbb{S}^1 \times \mathbb{C}) \) \end{center} & \begin{center} \( \mathrm{ab}(f_{\lambda}) \) \end{center} & \begin{center} \( \mathrm{c}(f_{\lambda}) \) \end{center} \\
\hline
\begin{center} 0 \end{center} & \begin{center} 2 connected components which do not interlace and wind 1 time above \( \mathbb{S}^1 \) \end{center} & \begin{center} \( \mathrm{Id} \) \end{center} & \begin{center} \( \mathrm{Id} \) \end{center} \\
\hline
\begin{center} 1 \end{center} & \begin{center} 1 connected component winding 2 times above \( \mathbb{S}^1 \) \end{center} & \begin{center} \( \sigma_1 \) \end{center} & \begin{center} \( (1,2) \) \end{center} \\
\hline
\begin{center} 2 \end{center} & \begin{center} 2 connected components which interlace 1 time and wind 1 time above \( \mathbb{S}^1 \) \end{center} & \begin{center} \( \sigma_1^2 \) \end{center} & \begin{center} \( \mathrm{Id} \) \end{center} \\
\hline
\end{tabular}
\end{center}

\begin{cor}[$d=2$, $p(z)=z^2$]\label{eq:abford=2}
    For each $\lambda\in\mathcal{D}'$, $\mathrm{ab}(f_{\lambda})=\sigma_1^{s(\lambda)}$. 
\end{cor}
We deduce that the map $\lambda\mapsto s(\lambda)$ of Astorg-Bianchi and our map $\lambda\mapsto\mathrm{ab}(f_{\lambda})$ carry exactly the same information and give the same classification of hyperbolic components in $\pi_0(\mathcal{D}')$.
 
\begin{otherlanguage}{english}

\end{otherlanguage}

$ $ \\
\noindent {\footnotesize Virgile Tapiero}\\
{\footnotesize Universit\'e d'Orl\'eans}\\
{\footnotesize CNRS, IDP - UMR 7013}\\
{\footnotesize F-45000 Orl\'eans, France}\\
{\footnotesize virgile.tapiero@univ-orleans.fr}\\

\end{document}